\documentclass[12pt,twoside,a4paper]{article}
\usepackage[utf8]{inputenc}
\usepackage{amsmath, amssymb, mathtools, enumerate}
\usepackage{mathrsfs} 
\usepackage{amsfonts, amsthm, stmaryrd, fullpage}
\usepackage{authblk}
\usepackage{tikz-cd}
\usepackage[german, english]{babel}
\usepackage[nottoc,numbib]{tocbibind}
\usepackage{needspace}

\usepackage{color}
\usepackage{hyperref}
\hypersetup{colorlinks=true,
linktoc=all,
citecolor=blue,
filecolor=black,
linkcolor=blue,
urlcolor=blue}

\newtheorem{thm}{Theorem}[section]
\newtheorem{lemma}[thm]{Lemma}
\newtheorem{cor}[thm]{Corollary}
\newtheorem{prop}[thm]{Proposition}

\theoremstyle{definition} 
\newtheorem{mydef}[thm]{Definition}
  
\newtheorem{example}[thm]{Example}

\theoremstyle{remark}
\newtheorem{rmk}[thm]{Remark}

\newcommand\ad{{\rm ad}}

\newcommand\alg{{\rm alg}}

\newcommand\gr{{\rm gr}}

\newcommand\Nil{{\rm Nil}}

\newcommand\perfd{{\rm perfd}}

\newcommand\Proj{{\rm Proj}}

\newcommand\spc{{\rm sp}}

\newcommand\Spa{{\rm Spa}}

\newcommand\Spec{{\rm Spec}} 
\newcommand\Spf{{\rm Spf}}
\newcommand\supp{{\rm supp}}

\begin{document}

\title{The $\Proj$ construction for pro-graded Huber pairs}
\author{Dimitri Dine}
\date{}
\maketitle

\begin{abstract}Inspired by the theory of pro-graded adic rings developed recently by Ishizuka and Yoshikawa, we introduce and study the notions of a pro-graded Huber ring and a pro-graded Huber pair \begin{equation*}((\mathcal{A}, \mathcal{A}_{\gr}), \mathcal{A}^{+}),\end{equation*}which capture the intuition of a complete Huber ring (or Huber pair) equipped with a dense graded subring. Using these notions, we define and study an analog \begin{equation*}\Proj(\mathcal{A}, \mathcal{A}_{\gr}, \mathcal{A}^{+})\end{equation*}of the $\Proj$ construction from algebraic geometry in the realm of adic spaces, without restricting our attention to adic spaces satisfying any Noetherian or finite-type assumptions. We consider gradings both by $\mathbb{Z}$ and by $\mathbb{Z}[1/p]$ for a fixed prime number $p$, so our $\Proj$ construction recovers both classical examples such as usual adic projective space $\mathbb{P}_{K}^{n,\ad}$ and infinite-type, perfectoid examples such as perfectoid projective space $\mathbb{P}_{K}^{n,\perfd}$. Using GAGA, we prove that the adic analytification of any projective variety $X$ over a nonarchimedean field $K$ can be described as the adic $\Proj$ of a suitable pro-graded affinoid algebra $\mathcal{A}$. Moreover, we prove that in this case $\dim X=\dim \mathcal{A}-1$. We also relate our adic $\Proj$ construction to formal models. \end{abstract} 

\tableofcontents

\section{Introduction}\label{sec:introduction}

\subsection{Motivation}

The $\Proj$ construction, which describes projective algebraic varieties in terms of graded algebras, is an indispensable tool of modern scheme-theoretic algebraic geometry. In nonarchimedean-analytic geometry, in the language of Berkovich spaces, there is the Berkovich $\Proj$ construction, introduced by Berkovich in \cite{Berkovich95Drinfeld}, proof of Theorem~1. The goal of this paper is to provide a general analog \begin{equation*}\Proj(\mathcal{A}, \mathcal{A}_{\gr}, \mathcal{A}^{+})\end{equation*}of these constructions in the setting of adic spaces which would recover classical examples (analytifications of projective varieties), but could also describe more exotic, ``infinite-dimensional'' geometric objects such as Zariski-closed subsets of perfectoid projective space. Moreover, we want our construction to be naturally related to formal models of the adic spaces in question. 

\subsection{Definition of pro-graded Huber rings}

One difficulty in finding such a general $\Proj$ construction is that, unlike the coordinate rings of algebraic varieties, the algebras of functions considered in $p$-adic geometry (complete Huber rings) are usually not equipped with any natural grading. Indeed, the most basic example of a complete Huber ring is the Tate algebra $T_{n, K}$ of convergent power series over the nonarchimedean base field $K$, and it is not a graded algebra in any obvious way. However, it does have a dense subalgebra which is equipped with a natural grading: The polynomial ring $K[X_1,\dots, X_n]$, with its usual grading by total degree. 

Starting from this elementary example and inspired by the theory of ``pro-graded adic rings'' developed recently by Ishizuka and Yoshikawa \cite{Ishizuka-Yoshikawa25}, we develop a formalism of pro-graded Huber rings and Huber pairs, which are complete Huber rings (or pairs) equipped with a dense graded subring. 
\begin{mydef}[Definition~\ref{Graded Huber ring}, Definition~\ref{Pro-graded Huber ring}]Let $G$ be an abelian group.~\begin{enumerate}[(1)]\item A \textit{$G$-graded Huber ring} $\mathcal{A}$ is a Huber ring which is also equipped with a structure of $G$-graded ring \begin{equation*}\mathcal{A}=\bigoplus_{g\in G}\mathcal{A}_{g}\end{equation*}and which has a pair of definition $(A, I)$ such that $A$ is a $G$-graded subring of $\mathcal{A}$ and $I$ is generated by finitely many homogeneous elements of degree~$0$. \item A \textit{pro-$G$-graded Huber ring} is a pair of Huber rings $(\mathcal{A}, \mathcal{A}_{\gr})$ where $\mathcal{A}_{\gr}$ is a $G$-graded Huber ring all of whose graded pieces are complete for the subspace topology and where $\mathcal{A}$ is the completion of $\mathcal{A}_{\gr}$. \end{enumerate} \end{mydef}
When the group $G$ in question is $\mathbb{Z}$, we simply speak of graded and pro-graded Huber rings. The most basic example is the one we already alluded to above. 
\begin{example}[Example~\ref{Pro-graded Tate algebra}]Let $K$ be a nonarchimedean field and let $n\geq1$. The polynomial ring $K[X_1,\dots, X_n]$, with topology defined by the Gauss norm, is a graded Huber ring and \begin{equation*}(T_{n, K}, K[X_1,\dots, X_n])\end{equation*}is a pro-graded Huber ring. When $p$ is a fixed prime number and $G=\mathbb{Z}[1/p]$, there is also an analog of the above example involving the perfectoid Tate algebra over a perfectoid field $K$, giving a rise to a pro-$\mathbb{Z}[1/p]$-graded Huber ring \begin{equation*}(T_{n, K}^{\perfd}, (T_{n, K}^{\perfd})_{\gr}),\end{equation*}see Example~\ref{Pro-graded perfectoid Tate algebra}.\end{example}
There is also a related notion of \textit{pro-$G$-graded Huber pair}, which is a triple \begin{equation*}((\mathcal{A}, \mathcal{A}_{\gr}), \mathcal{A}^{+})\end{equation*}such that $(\mathcal{A}, \mathcal{A}_{\gr})$ is a pro-graded Huber ring and $\mathcal{A}^{+}\subseteq \mathcal{A}$ is a ring of integral elements of $\mathcal{A}$ compatible with the graded subring $\mathcal{A}_{\gr}$ in a suitable sense (see Definition~\ref{Pro-graded Huber pair, definition}). Every classical affinoid algebra (i.e., every quotient of the Tate algebra $T_{n, K}$) is naturally an example of a pro-graded Huber ring. 
\begin{mydef}[Definition~\ref{Pro-graded affinoid algebras, definition}]Let $K$ be a nonarchimedean field. A pro-graded Huber ring $(\mathcal{A}, \mathcal{A}_{\gr})$ with a continuous map $K\to \mathcal{A}_{\gr,0}$ is called a \textit{pro-graded affinoid $K$-algebra} if there exists an integer $n\geq0$ and a morphism of pro-$\mathbb{Z}$-graded Tate rings \begin{equation*}(\varphi, \varphi_{\gr}): (T_{n, K}, (T_{n, K})_{\gr}))\to (\mathcal{A}, \mathcal{A}_{\gr}),\end{equation*}such that $\varphi$ is surjective.\end{mydef}

\subsection{Main results}

When $G=\mathbb{Z}$ or $G=\mathbb{Z}[1/p]$ for a fixed prime number $p$, we also define a notion of \textit{graded-sheafy} pro-graded Huber rings and pairs which is an analog of the notion of sheafy Huber rings in the pro-graded context, and every pro-graded affinoid $K$-algebra is pro-graded sheafy. For a projective variety $X$ over a nonarchimedean field $K$, we denote by $X^{\ad}$ the adic analytification of $X$. In particular, we denote by $\mathbb{P}_{K}^{n,\ad}$ the adic projective space over $K$. Similarly, we denote by $\mathbb{P}_{K}^{n,\perfd}$ the \textit{perfectoid projective space} (see \cite{Scholze}), which can be constructed by gluing multiple copies of the perfectoid closed unit polydisk along their common (perfectoid) polyannuli. Our results on the adic $\Proj$ construction can then be summarized as follows.  
\begin{thm}[Definition~\ref{Adic Proj, part 2}, Theorem~\ref{Dimension theorem}, Theorem~\ref{Graded GAGA}]\label{Main theorem}Suppose that $G=\mathbb{Z}$ or $G=\mathbb{Z}[1/p]$ for some fixed prime number $p$. For every graded-sheafy pro-$G$-graded Huber pair $((\mathcal{A}, \mathcal{A}_{\gr}), \mathcal{A}^{+})$, there is an adic space \begin{equation*}\Proj(\mathcal{A}, \mathcal{A}_{\gr}, \mathcal{A}^{+})\end{equation*}such that, for any nonarchimedean field $K$ (respectively, any perfectoid field $K$) and integer $n\geq1$, \begin{equation*}\mathbb{P}_{K}^{n,\ad}=\Proj(T_{n+1, K}, (T_{n+1, K})_{\gr}, T_{n+1, K}^{\circ})\end{equation*}and \begin{equation*}\mathbb{P}_{K}^{n,\perfd}=\Proj(T_{n+1, K}^{\perfd}, (T_{n+1, K}^{\perfd})_{\gr}, T_{n+1, K}^{\perfd\circ}),\end{equation*}and such that for every projective variety $X$ over $K$ there exists a pro-graded affinoid $K$-algebra $(\mathcal{A}, \mathcal{A}_{\gr})$ with \begin{equation*}X^{\ad}=\Proj(\mathcal{A}, \mathcal{A}_{\gr}, \mathcal{A}^{\circ}).\end{equation*}Moreover, for any pro-graded affinoid $K$-algebra $(\mathcal{A}, \mathcal{A}_{\gr})$, the dimension of $\Proj(\mathcal{A}, \mathcal{A}_{\gr}, \mathcal{A}^{\circ})$ is $\dim(\mathcal{A})-1$. \end{thm}
Our definition of $\Proj(\mathcal{A}, \mathcal{A}_{\gr}, \mathcal{A}^{+})$ rests on the notion of topologically homogeneous valuations on a pro-graded Huber ring. These valuations take on the role played by homogeneous prime ideals in the classical algebraic geometry picture. 
\begin{mydef}[Definition~\ref{Topologically homogeneous valuations}]Let $G$ be an abelian group and let $(\mathcal{A}, \mathcal{A}_{\gr})$ be a pro-$G$-graded Huber ring. A continuous valuation $v$ on $\mathcal{A}$ is called topologically homogeneous if  \begin{equation*}v(f)=\max\{\, v(f_{g})\mid g\in G\,\}\end{equation*}for every element $f\in \mathcal{A}_{\gr}$ with homogeneous decomposition $f=\sum_{g\in G}f_{g}$. \end{mydef}
For a pro-graded Huber pair $((\mathcal{A}, \mathcal{A}_{\gr}), \mathcal{A}^{+})$, kernels of topologically homogeneous valuations $v$ on $(\mathcal{A}, \mathcal{A}_{\gr})$ are homogeneous prime ideals of $\mathcal{A}$ and centers of topologically homogeneous valuations on the subring of integral elements $\mathcal{A}^{+}$ are homogeneous open prime ideals of $\mathcal{A}^{+}$. In this way, the description of the adic analytification $X^{\ad}$ of a projective variety $X$ as an adic $\Proj$ is compatible both with GAGA and with formal models of $X^{\ad}$. Regarding the comparison with formal models, we prove the following result. 
\begin{thm}[Definition~\ref{Formal Proj, part 2}, Theorem~\ref{Two specialization maps}, Corollary~\ref{Dimension estimate for the Proj}]\label{Main theorem 2}Suppose that $G=\mathbb{Z}$ or $G=\mathbb{Z}[1/p]$ for a fixed prime number $p$. There is a formal $\Proj$ construction \begin{equation*}\Proj(R, R_{\gr})\end{equation*}for pro-$G$-graded adic rings $(R, R_{\gr})$ in the sense of Ishizuka-Yoshukawa \cite{Ishizuka-Yoshikawa25} such that for every complete adic ring $A$, \begin{equation*}\widehat{\mathbb{P}_{A}^{n}}=\Proj(A\langle X_1,\dots, X_n\rangle, A[X_1,\dots, X_n]),\end{equation*}the formal projective space over $A$, and such that for every graded sheafy pro-$G$-graded uniform Tate Huber pair $((\mathcal{A}, \mathcal{A}_{\gr}), \mathcal{A}^{+})$ with $\mathcal{A}_{\gr}^{+}=\mathcal{A}^{+}\cap \mathcal{A}_{\gr}$, the adic projective spectrum is the adic analytic generic fiber \begin{equation*}\Proj(\mathcal{A}, \mathcal{A}_{\gr}, \mathcal{A}^{+})=\Proj(\mathcal{A}^{+}, \mathcal{A}_{\gr}^{+})_{\eta}^{\ad}\end{equation*}of the adic formal scheme $\Proj(\mathcal{A}^{+}, \mathcal{A}_{\gr}^{+})$. In particular, \begin{equation*}\dim(\Proj(\mathcal{A}, \mathcal{A}_{\gr}, \mathcal{A}^{+}))\geq \dim(\Proj(\mathcal{A}^{+}, \mathcal{A}_{\gr}^{+})).\end{equation*} \end{thm}
In our future work, we hope to generalize the theory begun in this paper and include a relative $\Proj$ construction for certain sheaves of pro-graded algebras on adic spaces and formal schemes. 

\subsection{Outline}

We conclude this introduction with an outline of the individual sections of the paper. Section~\ref{sec:Huber rings} recalls the notion of pro-graded adic rings from the work of Ishizuka and Yoshikawa~\cite{Ishizuka-Yoshikawa25} and then introduces graded and pro-graded Huber rings. Various analogs of notions from the theory of Huber rings are studied in the pro-graded setting: Pro-graded pairs of definition, pro-graded subrings, pro-graded rational localizations, etc. We also study (in some detail) the notion of topological homogeneous decomposition (suggested by the analogous notion for pro-graded adic rings introduced by Ishizuka and Yoshikawa), which turns out to be a very useful tool when proving results in the pro-graded formalism. Section~\ref{sec:Huber pairs} proceeds with the study of pro-graded Huber pairs, their rational localizations and quotients. Finally, topologically homogeneous valuations are introduced towards the end of Section~\ref{sec:Huber pairs}. In Section~\ref{sec:adic Proj} we use the notions of pro-graded Huber pairs and topologically homogeneous valuations to construct the adic projective spectrum of a graded-sheafy pro-graded Huber pair. We also discuss the examples of adic projective space and perfectoid projective space. In Section~\ref{sec:formal}, we introduce the formal $\Proj$ and prove Theorem~\ref{Main theorem 2}. Finally, in Section~\ref{sec:GAGA}, we prove the main assertions of Theorem~\ref{Main theorem} by relating our adic $\Proj$ construction from Section~\ref{sec:adic Proj} to projective varieties and GAGA. 

\subsection{Acknowledgements}

First of all, I would like to express my sincere gratitude to my advisor, Kiran Kedlaya, for his valuable advice and feedback on this work, as well as for his continuous encouragement and support. 

I would like to thank Kazuma Shimomoto for helpful conversations and for his insightful comments on my work, including on a previous version of this paper. I am also grateful to Linquan Ma for numerous enlightening conversations on commutative algebra. 

This paper owes an obvious intellectual debt to the ideas of Ishizuka and Yoshikawa on (pro-)graded perfectoid rings \cite{Ishizuka-Yoshikawa25}, which were the starting point and initial inspiration for the present work. I heartily thank Ryo Ishizuka for his very detailed comments on an early draft of my paper, which helped me catch some mistakes and improved this paper's exposition considerably, for kindly giving me permission to include his proof of Proposition~\ref{Pro-graded perfectoid Tate algebra and tilting}, as well as for all our stimulating correspondence over the years. I am also grateful to Shou Yoshikawa for helpful conversations around the notion of pro-graded perfectoid rings and the $\Proj$ construction. Finally, I would like to thank Jack J Garzella for numerous interesting discussions on the foundations of analytic geometry which continue to inform my perspective on the subject. 

The author gratefully acknowledges the support from Jane Street via the Jane Street Graduate Research Fellowship 2026. 

\subsection{AI disclosure}All results and proofs in the paper were obtained by the author without the use of generative AI, with the following exceptions: The reduction to finitely generated subgroups in the proof of Lemma~\ref{Integral closure and grading} was suggested to the author by Claude Sonnet (the author was lazy that day), and a literature search on the notion of graded valuations by Claude Fable~5.1 led to Fable explaining to the author that the two notions of graded valuations and homogeneous valuations appear in the literature in various special cases and that there is a presumably well-known connection between the two which can be proved by a ``leading-term argument''. This information ultimately led to the author's formulation and proof of Lemma~\ref{Gauss Lemma}. 

Other than on these occasions, Claude Fable 5.1 was used to proofread the paper and point out mistakes, which were subsequently corrected by the author.    

\section{Pro-graded Huber rings}\label{sec:Huber rings}

In order to study the projective geometry of adic spaces, we consider graded structures on complete Huber rings. More precisely, we observe that many examples of complete Huber rings occurring in adic geometry, such as the Tate algebra and the perfectoid Tate algebra, while not themselves equipped with a natural grading, contain a dense subring which has a natural structure of a graded ring. It is to the study of such complete Huber rings that this section is dedicated.

Throughout this section, we fix an abelian group $G$ and we are going to consider gradings by elements of $G$. The main examples we would like the reader to keep in mind are the usual Tate algebra $T_{n, K}=K\langle X_1,\dots, X_n\rangle$ over a nonarchimedean field $K$, where $G=\mathbb{Z}$ and the graded dense subring is the usual polynomial ring with its grading by total degree, and the perfectoid Tate algebra $T_{n, K}^{\perfd}=K\langle X_{1}^{1/p^{\infty}},\dots, X_{n}^{1/p^{\infty}}\rangle$, where $G=\mathbb{Z}[1/p]$ and the graded dense subring in question is a so-called ``gradedwise completion'' of the perfectly polynomial ring $K[X_{1}^{1/p^{\infty}},\dots, X_{n}^{1/p^{\infty}}]$ with the grading being once again the grading by total degree. Our story begins with the notions of graded and pro-graded adic rings introduced by Ishizuka and Yoshikawa \cite{Ishizuka-Yoshikawa25}, which we recall below. 
\begin{mydef}[Ishizuka-Yoshikawa, \cite{Ishizuka-Yoshikawa25}, Definition~3.1, Definition~3.9]\begin{enumerate}[(1)] \item A $G$-graded adic ring $R$ is a $G$-graded ring \begin{equation*}R=\bigoplus_{g\in G}R_{g}\end{equation*}whose underlying ring is equipped with an adic topology admitting a finitely generated homogeneous ideal of definition. A $G$-graded adic ring $R$ is called gradedwise complete if, for every $g\in G$, the subgroup $R_{g}$ is complete with respect to the subspace topology inherited from $R$. \item A pro-$G$-graded adic ring is a pair $(R, R_{\gr})$, where $R_{\gr}$ is a gradedwise complete $G$-graded adic ring and $R$ is the completion of $R_{\gr}$. \end{enumerate}For a pro-$G$-graded adic ring $(R, R_{\gr})$, an ideal $I$ of $R$ is called homogeneous if $I=JR$ for some homogeneous ideal $J$ of $R_{\gr}$.\end{mydef}
Graded and pro-graded adic rings form categories under the following notion of morphisms.
\begin{mydef}[Ishizuka-Yoshikawa \cite{Ishizuka-Yoshikawa25}]\begin{enumerate}[(1)]\item A morphism of $G$-graded adic rings $R\to S$ is a continuous ring homomorphism (for the chosen adic topologies) which is also a homomorphism of $G$-graded rings. \item A morphism of pro-$G$-graded adic rings, or a pro-$G$-graded homomorphism, \begin{equation*}(R, R_{\gr})\to (S, S_{\gr})\end{equation*}is a pair $(\varphi, \varphi_{\gr})$ consisting of a continuous ring homomorphism $R\to S$ and a morphism of $G$-graded adic rings $\varphi_{\gr}: R_{\gr}\to S_{\gr}$ such that $\varphi_{\gr}$ is the restriction of $\varphi$ to $R_{\gr}$ (equivalently, $\varphi=\widehat{\varphi_{\gr}}$). \end{enumerate}\end{mydef}
Analogously, we introduce a notion of (pro-)graded Huber rings. 
\begin{mydef}[Graded Huber ring]\label{Graded Huber ring}\begin{enumerate}[(1)]\item A $G$-graded Huber ring $\mathcal{A}$ is a (not necessarily complete) Huber ring with a $G$-grading \begin{equation*}\mathcal{A}=\bigoplus_{g\in G}\mathcal{A}_{g}\end{equation*}on its underlying ring such that there exists a $G$-graded subring $A$ of $\mathcal{A}$ which is a ring of definition of $\mathcal{A}$ and which has an ideal of definition generated by finitely many homogeneous elements of degree $0$. Such a ring of definition is called a $G$-graded ring of definition of $\mathcal{A}$. \item A pair of definition $(A, I)$ of a $G$-graded Huber ring $\mathcal{A}$ is called a $G$-graded pair of definition if $A$ is a $G$-graded ring of definition of $\mathcal{A}$ and $I$ is a homogeneous ideal of the $G$-graded ring $A$ generated by finitely many homogeneous elements of degree $0$. \item The $G$-graded Huber ring $\mathcal{A}$ is called gradedwise complete if, for every $g\in G$, the degree-$g$ part $\mathcal{A}_{g}$ of $\mathcal{A}$ is complete with respect to the subspace topology. \item A morphism of $G$-graded Huber rings is a continuous homomorphism of $G$-graded rings. \end{enumerate}\end{mydef}
For an ideal $I$ of a $G$-graded ring $A=\bigoplus_{g\in G}A_{g}$ and for any $g\in G$, we denote by $IA_{g}$ the image of $I$ under the projection map $A\to A_{g}$. 
\begin{lemma}\label{Ideals and gradings}Let $A=\bigoplus_{g\in G}A_{g}$ be a $G$-graded ring and let $I$ be a homogeneous ideal. Then $IA_{g}=I\cap A_{g}$ for every $g\in G$. \end{lemma}
\begin{proof}Every $f\in I\cap A_{g}$ is the degree-$g$ part of an element of $I$, namely, of $f$ itself. Thus $I\cap A_{g}\subseteq IA_{g}$. The other inclusion follows from the ideal $I$ being homogeneous.  \end{proof}
\begin{lemma}\label{Ideals and gradings 2}Let $\mathcal{A}$ be a $G$-graded Huber ring with $G$-graded pair of definition $(A, I)$. For every $g\in G$ and every integer $n\geq1$, we have $I^{n}A_{g}=I^{n}\cap \mathcal{A}_{g}$.\end{lemma}
\begin{proof}Since the ideal $I$ of the $G$-graded ring $A$ is homogeneous, we have $I^{n}A_{g}=I^{n}\cap A_{g}$ by Lemma~\ref{Ideals and gradings}. But $A_{g}=A\cap \mathcal{A}_{g}$ since $A$ is a $G$-graded subring of $\mathcal{A}$. \end{proof}
\begin{lemma}\label{Degree-0 part}Let $\mathcal{A}$ be a $G$-graded Huber ring with pair of definition $(A, I)$. For every $g\in G$, the family $\{I^{n}A_{g}\}_{n\geq0}$ is a fundamental system of open neighborhoods of $0$ for the subspace topology on $\mathcal{A}_{g}$. 
In particular, the degree-$0$ part $\mathcal{A}_{0}$, with the subspace topology induced from $\mathcal{A}$, is a Huber ring, with pair of definition $(A_{0}, IA_{0})$. \end{lemma}
\begin{proof}By Lemma~\ref{Ideals and gradings 2}, we see that $I^{n}A_{g}$ is an open neighborhood of $0$ with respect to the subspace topology for every $n\geq0$. Let $V\subseteq \mathcal{A}_{g}$ be an arbitrary open neighborhood of $0$ with respect to the subspace topology. Then there exists an open neighborhood $U$ of $0$ in $\mathcal{A}$ such that $V=U\cap \mathcal{A}_{g}$. Then there exists an integer $n\geq1$ such that $I^{n}\subseteq U$, so, by Lemma~\ref{Ideals and gradings 2}, $I^{n}A_{g}=I^{n}\cap \mathcal{A}_{g}\subseteq V$. \end{proof}
Recall from \cite{Ishizuka-Yoshikawa25}, Construction~3.3, the notion of gradedwise completion of a $G$-graded adic ring. For a $G$-graded adic ring $R$ with homogeneous ideal of definition $I$, this is the gradedwise complete $G$-graded adic ring \begin{equation*}R^{\wedge\gr}=\bigoplus_{g\in G}\varprojlim_{n\geq1}R_{g}/(I^{n}\cap R_{g}).\end{equation*}This gives rise to an analogous construction for graded Huber rings. 
\begin{mydef}[Gradedwise completion]\label{Gradedwise completion, definition}For any $G$-graded Huber ring $\mathcal{A}$, consider the $G$-graded ring \begin{equation*}\mathcal{A}^{\wedge\gr}=\bigoplus_{g\in G}\widehat{\mathcal{A}_{g}},\end{equation*}where $\widehat{\mathcal{A}_{g}}$ is the completion of $\mathcal{A}_{g}$ with respect to the subspace topology induced from $\mathcal{A}$. Choosing a $G$-graded pair of definition $(A, I)$ of $\mathcal{A}$, we endow $\mathcal{A}^{\wedge\gr}$ with the structure of a $G$-graded Huber ring by letting $(A^{\wedge\gr}, IA^{\wedge\gr})$ be a $G$-graded pair of definition. \end{mydef}
\begin{lemma}\label{Gradedwise completion and pairs of definition}The definition of gradedwise completion does not depend on the choice of the $G$-graded pair of definition $(A, I)$. \end{lemma} 
\begin{proof}Let $(A, I)$ and $(A^{\prime}, J)$ be two $G$-graded pairs of definition of $\mathcal{A}$. In view of Lemma~\ref{Degree-0 part} it suffices to prove that there exist integers $m, n\geq1$ such that \begin{equation*}\bigoplus_{g\in G}J^{m}A_{g}^{\prime}\subseteq \bigoplus_{g\in G}IA_{g}\end{equation*}and \begin{equation*}\bigoplus_{g\in G}I^{n}A_{g}\subseteq \bigoplus_{g\in G}JA_{g}^{\prime}.\end{equation*}By symmetry, it suffices to prove only one of these statements. Since $(A, I)$ and $(A^{\prime}, J)$ are pairs of definition of the same Huber ring, there exists an integer $m\geq1$ such that $J^{m}\subseteq I$. We check: \begin{align*}\bigoplus_{g\in G}J^{m}A_{g}^{\prime}=\bigoplus_{g\in G}(J^{m}\cap \mathcal{A}_{g}) \\ \subseteq \bigoplus_{g\in G}(I\cap \mathcal{A}_{g})=\bigoplus_{g\in G}IA_{g},\end{align*}where we applied Lemma~\ref{Ideals and gradings 2} twice. \end{proof} 
The grading on a graded Huber ring is automatically compatible with the topology, in the following sense.
\begin{lemma}\label{Grading and topology}Let $(f_{n})_{n}$ be a sequence in a $G$-graded Huber ring $\mathcal{A}$ and let $f\in \mathcal{A}$. Let $f=\sum_{g\in G}f_{g}$ and, for every $n$, let $f_{n}=\sum_{g\in G}f_{n, g}$ be the homogeneous decomposition of the element $f_{n}$. If $f_{n}\to f$ in $\mathcal{A}$, then $f_{n, g}\to f_{g}$ in $\mathcal{A}_{g}$ for every $g\in G$. More precisely, if $(A, I)$ is a $G$-graded pair of definition of $\mathcal{A}$, then for every integer $k\geq1$ there exists an integer $N_{k}$ (not depending on $g\in G$) such that for all $n\geq N_{k}$, \begin{equation*}f_{g}-f_{n, g}\in I^{k}\end{equation*}for all $g\in G$. \end{lemma}
\begin{proof}Let $(A, I)$ be a $G$-graded pair of definition. For every $k\geq1$, let $N_{k}$ be such that $f-f_{n}\in I^{k}$ for all $n\geq N_{k}$. Then $f_{g}-f_{n, g}\in I^{k}$ for all $g\in G$ and all $n\geq N_{k}$ since $I^{k}$ is a homogeneous ideal of $A$. \end{proof}
\begin{cor}\label{Completeness}For every $G$-graded Huber ring $\mathcal{A}=\bigoplus_{g\in G}\mathcal{A}_{g}$ and every $g\in G$, the degree-$g$ part $\mathcal{A}_{g}$ is closed in $\mathcal{A}$. In particular, if $\mathcal{A}$ is complete, then it is also gradedwise complete. \end{cor}
The converse to Corollary~\ref{Completeness} fails: There exist many $G$-graded Huber rings which are gradedwise complete, but not complete in the usual sense. 
\begin{example}[Polynomial ring]\label{Polynomial ring}Let $n\geq1$ be an integer, let $K$ be a nonarchimedean field with pseudo-uniformizer $\varpi$ and consider the polynomial ring $K[X_1,\dots, X_n]$, regarded as a Tate ring with pair of definition $(K^{\circ}[X_1,\dots, X_n], \varpi)$. With its usual grading by total degree, $K[X_1,\dots, X_n]$ is a $\mathbb{Z}$-graded Huber ring, and it is gradedwise complete: For every $d\in \mathbb{Z}_{\geq0}$, the degree-$d$ piece $K[X_1,\dots, X_n]_{d}$ is the finite-dimensional $K$-vector space of homogeneous polynomials of total degree $d$, which is complete for its unique norm topology as a finite-dimensional $K$-vector space. But the Tate ring $K[X_1,\dots, X_n]$ is obviously not complete, its completion being the Tate algebra $T_{n, K}=K\langle X_1,\dots, X_n\rangle$. \end{example}
This observation naturally leads us to the notion of a pro-graded Huber ring, analogous to Ishizuka's and Yoshikawa's notion of pro-graded adic ring which we have just recalled.  
\begin{mydef}[Pro-graded Huber ring]\label{Pro-graded Huber ring}\begin{enumerate}[(1)]\item A pro-$G$-graded Huber ring is a pair $(\mathcal{A}, \mathcal{A}_{\gr})$, where $\mathcal{A}$ is a Huber ring and $\mathcal{A}_{\gr}$ is a gradedwise complete $G$-graded Huber ring such that $\mathcal{A}$ is the completion of the Huber ring $\mathcal{A}_{\gr}$. \item A morphism of pro-$G$-graded Huber rings \begin{equation*}(\varphi, \varphi_{\gr}): (\mathcal{A}, \mathcal{A}_{\gr})\to (\mathcal{B}, \mathcal{B}_{\gr})\end{equation*}is a pair consisting of a continuous homomorphism $\varphi: \mathcal{A}\to \mathcal{B}$ and a continuous homomorphism of $G$-graded rings $\mathcal{A}_{\gr}\to \mathcal{B}_{\gr}$ with $\varphi\vert_{\mathcal{A}_{\gr}}=\varphi_{\gr}$. \item If $(\mathcal{A}, \mathcal{A}_{\gr})$ is a pro-$G$-graded Huber ring, an ideal $I$ of $\mathcal{A}$ is called homogeneous if there exists a homogeneous ideal $J$ of $\mathcal{A}_{\gr}$ such that $J\mathcal{A}=I$.\end{enumerate}\end{mydef}
We also have a corresponding notion of pro-graded pair of definition.
\begin{mydef}[Pro-graded pair of definition]\label{Pro-graded pair of definition}\begin{enumerate}[(1)]Let $(\mathcal{A}, \mathcal{A}_{\gr})$ be a pro-$G$-graded Huber ring. \item A ring of definition $A$ of $\mathcal{A}$ is called a pro-$G$-graded ring of definition of $(\mathcal{A}, \mathcal{A}_{\gr})$ if there exists a $G$-graded ring of definition $A_{\gr}$ of $\mathcal{A}_{\gr}$ such that $A$ is the completion of the adic ring $A_{\gr}$. \item A pair of definition $(A, I)$ of $\mathcal{A}$ is called a pro-$G$-graded pair of definition of $(\mathcal{A}, \mathcal{A}_{\gr})$ if there exists a $G$-graded pair of definition $(A_{\gr}, J)$ of $\mathcal{A}_{\gr}$ with \begin{equation*}(A, I)=(\widehat{A_{\gr}}, J\widehat{A_{\gr}}),\end{equation*}i.e., $A$ is the completion of the adic ring $A_{\gr}$ and $I=JA$.\end{enumerate}\end{mydef} 
Our main examples of pro-graded Huber rings are going to be Tate.
\begin{mydef}[Graded Tate ring, pro-graded Tate ring]A $G$-graded Huber ring $\mathcal{A}$ is called a $G$-graded Tate ring if the Huber ring $\mathcal{A}$ is Tate and there exists a $G$-graded pair of definition $(A, (\varpi)_{A})$ such that $\varpi\in A_{0}\cap \mathcal{A}^{\times}$, i.e., $\varpi$ is a unit in $\mathcal{A}$ and is homogeneous of degree $0$ in the $G$-graded ring $A$. A pro-$G$-graded Huber ring $(\mathcal{A}, \mathcal{A}_{\gr})$ is called a pro-$G$-graded Tate ring if $\mathcal{A}_{\gr}$ is a $G$-graded Tate ring.\end{mydef}
We now introduce our two main examples of the above concepts.  
\begin{example}[Pro-graded Tate algebra]\label{Pro-graded Tate algebra}Consider \begin{equation*}T_{n, K}=K\langle X_1,\dots, X_n\rangle,\end{equation*}the Tate algebra of convergent power series in $n\geq1$ variables over some nonarchimedean field $K$. Let $\varpi$ be a pseudo-uniformizer of $K$. For every $d\in\mathbb{Z}$, let $(T_{n, K}^{\circ})_{\gr, d}$ be the $K^{\circ}$-module of homogeneous polynomials $f\in K^{\circ}[X_1,\dots, X_n]$ of (total) degree $d$. This is a finite free $K^{\circ}$-module and thus $\varpi$-adically complete. Set \begin{equation*}(T_{n, K}^{\circ})_{\gr}=\bigoplus_{d\in\mathbb{Z}}(T_{n, K}^{\circ})_{\gr, d}.\end{equation*}Note that this is a gradedwise $\varpi$-adically complete $\mathbb{Z}$-graded adic ring, and that the pair $(T_{n, K}^{\circ}, (T_{n, K}^{\circ})_{\gr})$ is a pro-$\mathbb{Z}$-graded adic ring.

Set \begin{equation*}(T_{n, K})_{\gr, d}=(T_{n, K}^{\circ})_{\gr, d}[\varpi^{-1}]=(T_{n, K}^{\circ})_{\gr, d}\otimes_{K^{\circ}}K\end{equation*}and \begin{equation*}(T_{n, K})_{\gr}=\bigoplus_{d\in\mathbb{Z}}(T_{n, K})_{\gr, d}.\end{equation*}Then $(T_{n, K})_{\gr}$ is a $\mathbb{Z}$-graded Tate ring and the pair \begin{equation*}(T_{n, K}, (T_{n, K})_{\gr})\end{equation*}is a pro-$\mathbb{Z}$-graded Tate ring, which we call the pro-graded Tate algebra. Moreover, we note that the pro-$\mathbb{Z}$-graded adic ring $(T_{n, K}^{\circ}, (T_{n, K}^{\circ})_{\gr})$ is a pro-$\mathbb{Z}$-graded ring of definition of the pro-graded Tate algebra, and we call it the integral pro-graded Tate algebra.  \end{example}
\begin{example}[Pro-graded perfectoid Tate algebra]\label{Pro-graded perfectoid Tate algebra}Analogously, let $K$ be a perfectoid field, with pseudo-uniformizer $\varpi$, and consider the perfectoid Tate algebra \begin{equation*}T_{n, K}^{\perfd}=K\langle X_{1}^{1/p^{\infty}},\dots, X_{n}^{1/p^{\infty}}\rangle\end{equation*}in $n\geq1$ variables over $K$. For every $d\in \mathbb{Z}[1/p]$, the $K^{\circ}$-module of homogeneous ``perfect polynomials'' $f\in K^{\circ}[X_{1}^{1/p^{\infty}},\dots, X_{n}^{1/p^{\infty}}]$ is again a free $K^{\circ}$-module, but for $n>1$ it is of infinite rank and thus no longer $\varpi$-adically complete. Let $(T_{n, K}^{\perfd\circ})_{\gr, d}$ be the $\varpi$-adic completion of this $K^{\circ}$-module. Its elements are elements $f\in T_{n, K}^{\perfd\circ}$ all of whose monomials are homogeneous of degree $d$. For example, \begin{equation*}f=\sum_{k=0}^{\infty}\varpi^{k}X_{1}^{1/p^{k}}X_{2}^{(p^{k}-1)/p^{k}}\in K^{\circ}\langle X_{1}^{1/p^{\infty}}, X_{2}^{1/p^{\infty}}\rangle \end{equation*}is an element of $(T_{2, K}^{\perfd\circ})_{\gr, 1}$ which does not belong to $K^{\circ}[X_{1}^{1/p^{\infty}}, X_{2}^{1/p^{\infty}}]$. Note that this phenomenon only occurs in the perfectoid case and has no analog for the usual polynomial ring and the usual Tate algebra as in Example~\ref{Pro-graded Tate algebra}.

We set \begin{equation*}(T_{n, K}^{\perfd\circ})_{\gr}=\bigoplus_{d\in\mathbb{Z}[1/p]}(T_{n, K}^{\perfd\circ})_{\gr, d}.\end{equation*}This is a gradedwise $\varpi$-adically complete $\mathbb{Z}[1/p]$-graded adic ring and the pair $(T_{n, K}^{\perfd\circ}, (T_{n, K}^{\perfd\circ})_{\gr})$ is a pro-$\mathbb{Z}[1/p]$-graded adic ring. Similarly to the case of the usual Tate algebra, we define $(T_{n, K}^{\perfd})_{\gr, d}=(T_{n, K}^{\perfd\circ})_{\gr, d}\otimes_{K^{\circ}}K$ for every $d\in\mathbb{Z}[1/p]$ and \begin{equation*}(T_{n, K}^{\perfd})_{\gr}=\bigoplus_{d\in\mathbb{Z}[1/p]}(T_{n, K}^{\perfd})_{\gr, d}.\end{equation*}This is a $\mathbb{Z}[1/p]$-graded Tate ring, and the pair \begin{equation*}(T_{n, K}^{\perfd}, (T_{n, K}^{\perfd})_{\gr})\end{equation*}is a pro-$\mathbb{Z}[1/p]$-graded Tate ring, which we call the pro-graded perfectoid Tate algebra. The pair $(T_{n, K}^{\perfd\circ}, (T_{n, K}^{\perfd\circ})_{\gr})$ is a pro-graded ring of definition of $(T_{n, K}^{\perfd}, (T_{n, K}^{\perfd})_{\gr})$ and we call it the integral pro-graded perfectoid Tate algebra.\end{example}
\begin{rmk}\label{Convergent power series, general case}By the same recipe, we also obtain, for every uniform Tate Huber pair $(\mathcal{A}, \mathcal{A}^{+})$, the pro-graded convergent power series algebra \begin{equation*}(\mathcal{A}\langle X_1,\dots, X_n\rangle, (\mathcal{A}\langle X_{1},\dots, X_{n}\rangle)_{\gr})\end{equation*}over $(\mathcal{A}, \mathcal{A}^{+})$, and if $(\mathcal{A}, \mathcal{A}^{+})$ is perfectoid, also the perfectoid convergent power series algebra \begin{equation*}(\mathcal{A}\langle X_{1}^{1/p^{\infty}},\dots, X_{n}^{1/p^{\infty}}\rangle, (\mathcal{A}\langle X_{1}^{1/p^{\infty}},\dots, X_{n}^{1/p^{\infty}}\rangle)_{\gr})\end{equation*}over $(\mathcal{A}, \mathcal{A}^{+})$, where $(\mathcal{A}, \mathcal{A}^{+})$ replaces $(K, K^{\circ})$ in all the constructions of Example~\ref{Pro-graded Tate algebra} and Example \ref{Pro-graded perfectoid Tate algebra}. We also obtain the integral analogs of these constructions. \end{rmk}
Another source of examples of pro-graded Huber rings comes from rational localizations. Recall that for any $G$-graded ring \begin{equation*}\mathcal{A}=\bigoplus_{g\in G}\mathcal{A}_{g}\end{equation*}and any homogeneous element $f\in \mathcal{A}$ of degree $h$, the localization $\mathcal{A}[\frac{1}{f}]$ can be endowed with a canonical structure of a $G$-graded ring as follows. We say that an element \begin{equation*}F\in \mathcal{A}[\frac{1}{f}]\end{equation*}is homogeneous of degree $g\in G$ in $\mathcal{A}[\frac{1}{f}]$ if there exists an integer $k\geq0$ such that \begin{equation*}f^{k}F\in \mathcal{A}_{g+kh}.\end{equation*}This does not depend on the choice of an integer $k\geq0$ with $f^{k}F\in \mathcal{A}$: Indeed, if $f^{k}F\in \mathcal{A}_{g+kh}$ and $k'\geq0$ is another integer such that $f^{k'}F\in \mathcal{A}$, then $f^{k+k'}F\in \mathcal{A}_{g+(k+k')h}$ and thus $f^{k'}F\in \mathcal{A}_{g+k'h}$ (the kernel of $\mathcal{A}\to \mathcal{A}[\frac{1}{f}]$ is homogeneous: If an element is annihilated by $f^{n}$ for some $n\geq1$, so are its homogeneous components). In what follows, we always view the localization $\mathcal{A}[\frac{1}{f}]$ of a $G$-graded ring $\mathcal{A}$ with respect to a homogeneous element $f\in\mathcal{A}$ as a $G$-graded ring with this $G$-graded ring structure. 
\begin{lemma}\label{Localization}Let $\mathcal{A}$ be a $G$-graded Huber ring and let $f\in \mathcal{A}$ be a homogeneous element of degree $0$. Endow the localization $\mathcal{A}[f^{-1}]$ with its natural $G$-grading and with its usual Huber ring topology. Then \begin{equation*}\mathcal{A}[f^{-1}]_{0}=\mathcal{A}_{0}[f^{-1}]\end{equation*}as Huber rings.\end{lemma}
\begin{proof}It is clear that the right-hand side is contained in the left-hand side. On the other hand, an element of $\mathcal{A}[f^{-1}]$ is homogeneous of degree $g\in G$ if and only if it is of the form $\frac{a}{f^{n}}$ with $a$ homogeneous of degree $\deg(a)=g+n\deg(f)$. Since $f$ is homogeneous of degree $0$, this simplifies to $\deg(a)=g$. This shows that $\mathcal{A}[f^{-1}]_{0}$ and $\mathcal{A}_{0}[f^{-1}]$ are equal as abstract subrings of $\mathcal{A}[f^{-1}]$. 

Similarly, if $(A, I)$ is a $G$-graded pair of definition of $\mathcal{A}$, by writing an arbitrary element of $A[f^{-1}]$ as $\sum_{i}a_{i}f^{-k_{i}}$ and using the uniqueness of homogeneous decompositions, we see that $A[f^{-1}]_{0}=A_{0}[f^{-1}]$. 

Set $I_{0}=I\cap A_{0}=I\cap \mathcal{A}_{0}$. By Lemma~\ref{Degree-0 part}, $(A_{0}, I_0)$ is a pair of definition of the Huber ring $\mathcal{A}_{0}$, so \begin{equation*}(A_{0}[f^{-1}], IA_{0}[f^{-1}])\end{equation*}is a pair of definition of the Huber ring $\mathcal{A}_{0}[f^{-1}]$. On the other hand, the Huber ring $\mathcal{A}[f^{-1}]$ has pair of definition $(A[f^{-1}], IA[f^{-1}])$. Applying what we proved in the first and second paragraphs, we see that $(A_{0}[f^{-1}], JA_{0}[f^{-1}])$ is a pair of definition of $\mathcal{A}[f^{-1}]_{0}$, proving equality of the two Huber rings. \end{proof} 
\begin{lemma}\label{Lemma on rational localizations}Let $(\mathcal{A}, \mathcal{A}_{\gr})$ be a pro-$G$-graded Huber ring and let $f_1,\dots, f_{n}, f_{n+1}$ be homogeneous elements of $\mathcal{A}_{\gr}$ generating an open ideal in $\mathcal{A}$. If $A$ is a $G$-graded subring of $\mathcal{A}_{\gr}$, then the subring \begin{equation*}A[\frac{f_1,\dots, f_n}{f_{n+1}}]\subseteq \mathcal{A}_{\gr}[\frac{1}{f_{n+1}}]\end{equation*}of the $G$-graded ring $\mathcal{A}_{\gr}[\frac{1}{f_{n+1}}]$ (with its natural $G$-graded ring structure) is a $G$-graded subring.\end{lemma}
\begin{proof}Let \begin{equation*}f=\sum_{i=(i_{1},\dots, i_{n})}a_{i}(\frac{f_{1}}{f_{n+1}})^{i_{1}}\cdot~\dots~\cdot (\frac{f_{n}}{f_{n+1}})^{i_{n}}\in A[\frac{f_1,\dots, f_{n}}{f_{n+1}}],\end{equation*}with $a_{i}\in A$. Since the elements $f_1,\dots, f_{n}, f_{n+1}$ are homogeneous in $\mathcal{A}_{\gr}$, the elements $\frac{f_{1}}{f_{n+1}},\dots, \frac{f_{n}}{f_{n+1}}$ are homogeneous in $\mathcal{A}_{\gr}[\frac{1}{f_{n+1}}]$. Let $a_{i}=\sum_{g\in G}a_{i, g}$ be the homogeneous decomposition of each $a_{i}$ in $\mathcal{A}_{\gr}$. Since $A$ is a $G$-graded subring of $\mathcal{A}_{\gr}$, we have $a_{i, g}\in A$ for all $g\in G$. Therefore, \begin{equation*}f=\sum_{i=(i_1,\dots, i_{n})}\sum_{g\in G}a_{i, g}(\frac{f_{1}}{f_{n+1}})^{i_{1}}\cdot~\dots~\cdot (\frac{f_{n}}{f_{n+1}})^{i_{n}}\end{equation*}represents $f$ as a finite sum of elements of $A[\frac{f_1,\dots, f_n}{f_{n+1}}]$ which are homogeneous in the $G$-graded ring $\mathcal{A}_{\gr}[\frac{1}{f_{n+1}}]$. By uniqueness of homogeneous decompositions, this representation must then coincide with the homogeneous decomposition of $f$ in the $G$-graded ring $\mathcal{A}_{\gr}[\frac{1}{f_{n+1}}]$.\end{proof}
\begin{mydef}[Pro-graded rational localizations]\label{Pro-graded rational localizations}Let $(\mathcal{A}, \mathcal{A}_{\gr})$ be a pro-$G$-graded Huber ring and let $f_{1},\dots, f_{n}, f_{n+1}\in \mathcal{A}_{\gr}$ be homogeneous elements generating an open ideal of $\mathcal{A}_{\gr}$. Let $h=\deg(f_{n+1})$, and let $(A, I)$ be a $G$-graded pair of definition of the $G$-graded Huber ring $\mathcal{A}_{\gr}$. By Lemma~\ref{Lemma on rational localizations}, the subring \begin{equation*}A[\frac{f_{1},\dots, f_{n}}{f_{n+1}}]\subseteq \mathcal{A}_{\gr}[\frac{1}{f_{n+1}}]\end{equation*}of the $G$-graded ring $\mathcal{A}_{\gr}[\frac{1}{f_{n+1}}]$ (with its natural $G$-graded ring structure) is a $G$-graded subring. This turns the dense subring $\mathcal{A}_{\gr}[\frac{1}{f_{n+1}}]$ into a $G$-graded Huber ring, with $G$-graded pair of definition \begin{equation*}(A[\frac{f_1,\dots, f_n}{f_{n+1}}], IA[\frac{f_1,\dots, f_n}{f_{n+1}}]).\end{equation*}Consequently,\begin{equation*}(\mathcal{A}\langle \frac{f_1,\dots, f_n}{f_{n+1}}\rangle, (\mathcal{A}\langle \frac{f_1,\dots, f_n}{f_{n+1}}\rangle)_{\gr}):=(\mathcal{A}\langle \frac{f_1,\dots, f_n}{f_{n+1}}\rangle, \mathcal{A}_{\gr}[\frac{1}{f_{n+1}}]^{\wedge\gr})\end{equation*}is a pro-$G$-graded Huber ring equipped with a canonical homomorphism of pro-$G$-graded Huber rings from $(\mathcal{A}, \mathcal{A}_{\gr})$, and we call pro-$G$-graded Huber rings constructed from $(\mathcal{A}, \mathcal{A}_{\gr})$ in this way (pro-$G$-graded) rational localizations of the pro-$G$-graded Huber ring $(\mathcal{A}, \mathcal{A}_{\gr})$. \end{mydef}
The pro-graded rational localization satisfies a universal property analogous to that of usual rational localizations. 
\begin{lemma}\label{Universal property of pro-graded rational localization}Let $(\mathcal{A}, \mathcal{A}_{\gr})$ be a pro-$G$-graded Huber ring and let $f_{1}, \dots, f_{n}, f_{n+1}$ be homogeneous elements of $\mathcal{A}_{\gr}$ generating an open ideal in $\mathcal{A}$. Let \begin{equation*}(\varphi, \varphi_{\gr}): (\mathcal{A}, \mathcal{A}_{\gr})\to (\mathcal{B}, \mathcal{B}_{\gr})\end{equation*}be a morphism of pro-$G$-graded Huber rings such that $\varphi(f_{n+1})=\varphi_{\gr}(f_{n+1})$ is invertible in $\mathcal{B}_{\gr}$ and such that the elements $\frac{\varphi(f_{1})}{\varphi(f_{n+1})},\dots, \frac{\varphi(f_{n})}{\varphi(f_{n+1})}$ are power-bounded in $\mathcal{B}$. Then $(\varphi, \varphi_{\gr})$ factors uniquely through the canonical morphism \begin{equation*}(\mathcal{A}, \mathcal{A}_{\gr})\to (\mathcal{A}\langle \frac{f_1,\dots, f_n}{f_{n+1}}\rangle, (\mathcal{A}\langle \frac{f_1,\dots, f_n}{f_{n+1}}\rangle)_{\gr}).\end{equation*}\end{lemma}
\begin{proof}By the usual universal property of rational localizations, $\varphi$ factors uniquely through the canonical map $\mathcal{A}\to \mathcal{A}\langle \frac{f_1,\dots, f_n}{f_{n+1}}\rangle$. Since $f_{n+1}$ is invertible in $\mathcal{B}_{\gr}$, the map $\varphi_{\gr}: \mathcal{A}_{\gr}\to \mathcal{B}_{\gr}$ also factors through $\mathcal{A}_{\gr}[\frac{1}{f_{n+1}}]$, so it remains to show that the resulting continuous map $\psi_{\gr}: \mathcal{A}_{\gr}[\frac{1}{f_{n+1}}]\to \mathcal{B}_{\gr}$ is a map of $G$-graded rings (in which case, it extends by continuity to a continuous map of $G$-graded rings between the gradedwise completions). Let $f\in \mathcal{A}_{\gr}[\frac{1}{f_{n+1}}]$ be homogeneous of some degree $g\in G$. Then $f_{n+1}^{k}f\in \mathcal{A}_{\gr, g+kh}$, where $h$ is the degree of $f_{n+1}$. Since $\varphi_{\gr}$ is a homomorphism of $G$-graded rings, $\varphi_{\gr}(f_{n+1}^{k}f)\in \mathcal{B}_{\gr, g+kh}$. It follows that $\psi_{\gr}(f)\in\mathcal{B}_{\gr, g}$. \end{proof}
\begin{lemma}\label{Pro-graded rational localizations, special case}Let $(\mathcal{A}, \mathcal{A}_{\gr})$ be a pro-$G$-graded Huber ring and let $f\in \mathcal{A}$ be a homogeneous element of degree $0$. Then, for every $g\in G$, we have \begin{equation*}\mathcal{A}_{\gr}[f^{-1}]_{0}=\mathcal{A}_{\gr,0}[f^{-1}]\end{equation*}as Huber rings (where $\mathcal{A}_{\gr,0}$ is viewed as a Huber ring per Lemma~\ref{Degree-0 part}). In particular, \begin{equation*}\mathcal{A}\langle f^{-1}\rangle_{\gr,0}=\mathcal{A}_{\gr,0}\langle f^{-1}\rangle.\end{equation*}\end{lemma}
\begin{proof}The first equality is Lemma~\ref{Localization} applied to the $G$-graded Huber ring $\mathcal{A}_{\gr}$. The second equality follows from the first by taking completions of both Huber rings. \end{proof} 
We also have an analog in the setting of pro-graded Huber rings of the topological homogeneous decomposition of \cite{Ishizuka-Yoshikawa25}, Definition~3.14.
\begin{lemma}\label{Topological homogeneous decomposition}Every element $f\in\mathcal{A}$ in a pro-$G$-graded Huber ring $(\mathcal{A}, \mathcal{A}_{\gr})$ can be uniquely written as a series \begin{equation*}f=\sum_{g\in G}f_{g}\end{equation*}with $f_{g}\in \mathcal{A}_{\gr, g}$ for all $g\in G$, where the support \begin{equation*}\supp(f)=\{\, g\in G\mid f_{g}\neq0\,\}\end{equation*}is countable and the series is convergent in the sense that for every neighborhood $U$ of zero in $\mathcal{A}$ there exist only finitely many $g\in G$ such that $f_{g}\not\in U$. \end{lemma}
\begin{proof}Let $(A_{\gr}, J)$ be a $G$-graded pair of definition of the $G$-graded Huber ring $\mathcal{A}_{\gr}$. Then $(A, I)=(\widehat{A_{\gr}}, J\widehat{A_{\gr}})$, where $\widehat{A_{\gr}}$ is the $J$-adic completion of $A_{\gr}$, is a (pro-$G$-graded) pair of definition of $\mathcal{A}$. For every $f\in \mathcal{A}$ and every integer $n\geq1$, we can find an element $f_{n}\in \mathcal{A}_{\gr}$ such that \begin{equation*}f-f_{n}\in I^{n}.\end{equation*}For every $n$, let \begin{equation*}f_{n}=\sum_{g\in G}f_{n, g}\end{equation*}be the homogeneous decomposition of the element $f_{n}$ of $\mathcal{A}_{\gr}$. Fix $k\geq1$. We want to prove that the set of $g\in G$ with $f_{g}\not\in I^{k}$ is finite. Choose an integer $N_{k}\geq1$ such that $f_{m}-f_{n}\in J^{k}=I^{k}\cap A_{\gr}$ for all $n, m\geq N_{k}$. Since $J^{k}$ is a homogeneous ideal of $A_{\gr}$, this implies \begin{equation*}f_{m, g}-f_{n, g}\in J^{k}\end{equation*}for all $n, m\geq N_{k}$ and for all $g\in G$. Since $A_{\gr}$ is gradedwise complete, there exist elements $f_{g}\in A_{\gr}$, $g\in G$, such that $(f_{n, g})_{n}$ converges to $f_{g}$. Suppose that $f_{g}-f_{n, g}\not\in J^{k}$ for some $g\in G$ and $n\geq N_{k}$. Then, since $f_{m, g}-f_{n, g}\in J^{k}$ for all $g\in G$ and all $m\geq N_{k}$, we see that $f_{g}-f_{m, g}\not\in J^{k}$ for all $m\geq N_{k}$, which contradicts the assumption that $f_{g}=\lim_{n\to\infty}f_{n, g}$. We conclude that \begin{equation*}f_{g}-f_{n, g}\in J^{k}\subseteq I^{k}\end{equation*}for all $g\in G$ and for all $n\geq N_{k}$. In particular, $f_{g}-f_{N_{k}, g}\in I^{k}$ for all $g\in G$. Then $f_{g}\not\in I^{k}$ implies $f_{N_{k}, g}\not\in I^{k}$ and, in particular, $f_{N_{k}, g}\neq0$. Since $f_{N_{k}}\in \mathcal{A}_{\gr}$ has only finitely many non-zero homogeneous components $f_{N_{k}, g}$, we see that $\{\, g\in G\mid f_{g}\not\in I^{k}\,\}$ is finite, as claimed. Since $k\geq1$ was arbitrary, this shows that $\sum_{g\in G}f_{g}$ is a convergent series in $\mathcal{A}$ with countable support, in the sense explained in the statement of the lemma. But then, choosing an increasing sequence of finite subsets $(G_{m})_{m}$ with union $\supp(f)$, we get \begin{equation*}f=\lim_{n\to\infty}f_{n}=\lim_{n\to\infty}\sum_{g\in G}f_{n, g}=\lim_{n\to\infty}\lim_{m\to\infty}\sum_{g\in G_{m}}f_{n, g}=\lim_{m\to\infty}\sum_{g\in G_{m}}\lim_{n\to\infty}f_{n, g}=\sum_{g\in G}f_{g}.\end{equation*}

It remains to prove that the representation $f=\sum_{g\in G}f_{g}$ is unique. Let $f=\sum_{g\in G}f_{g}^{\prime}$ be another such representation. Since $\supp(f)$ is countable, we can choose an increasing sequence $(G_{n})_{n}$ of finite subsets of $\supp(f)$ with \begin{equation*}\supp(f)=\bigcup_{n}G_{n}.\end{equation*}Then $f_{n}^{\prime}=\sum_{g\in G_{n}}f_{g}^{\prime}$ defines a sequence $(f_{n}^{\prime})_{n}$ in $\mathcal{A}_{\gr}$ such that $f_{n}^{\prime}\to f$ in $\mathcal{A}$. Moreover, the degree-$g$ part of $f_{n}^{\prime}$ is precisely $f_{g}^{\prime}$ if $g\in G_{n}$, and zero otherwise. For every $N\geq1$, there exists $k\geq1$ such that for all $n\geq k$, we have \begin{equation*}f_{g}, f_{g}^{\prime}\in I^{N}\end{equation*}for $g\not\in G_{n}$. Then \begin{equation*}\sum_{g\in G_{n}}(f_{g}-f_{g}^{\prime})=(f-f_{n}^{\prime})-\sum_{g\not\in G_{n}}(f_{g}-f_{g}^{\prime})\in I^{N}\end{equation*}for all $n$ large enough. Thus, all $n$ large enough, we have \begin{equation*}\sum_{g\in G_{n}}(f_{g}-f_{g}^{\prime})\in I^{N}\cap A_{\gr}=J^{N}\widehat{A_{\gr}}\cap A_{\gr}=J^{N}.\end{equation*}Since $J^{N}$ is a homogeneous ideal of $A_{\gr}$ and $\supp(f)$ is the increasing union of the finite subsets $G_{n}$, this implies $f_{g}-f_{g}^{\prime}\in J^{N}\subseteq I^{N}$ for all $g\in G$. But $N\geq1$ was arbitrary and $A$ is $I$-adically complete, so $f_{g}=f_{g}^{\prime}$ for all $g\in G$, as desired. \end{proof}
\begin{mydef}[Topological homogeneous decomposition]\label{Topological homogeneous decomposition, definition}As in \cite{Ishizuka-Yoshikawa25}, Definition~3.13, we call the unique convergent series representation \begin{equation*}f=\sum_{g\in G}f_{g}\end{equation*}from Lemma~\ref{Topological homogeneous decomposition} the topological homogeneous decomposition of the element $f\in \mathcal{A}$ of the pro-$G$-graded Huber ring $(\mathcal{A}, \mathcal{A}_{\gr})$.\end{mydef}
We have the following analog of Lemma~\ref{Grading and topology} for topological homogeneous decompositions.
\begin{lemma}\label{Grading and topology 2}Let $(\mathcal{A}, \mathcal{A}_{\gr})$ be a pro-$G$-graded Huber ring and let $(f_{n})_{n}$ be a sequence in $\mathcal{A}$ converging to some element $f\in \mathcal{A}$. Let $f=\sum_{g\in G}f_{g}$ and $f_{n}=\sum_{g\in G}f_{n, g}$ be the topological homogeneous decompositions of $f$ and $f_{n}$, respectively. Then $(f_{n, g})_{n}$ converges to $f_{g}$ for every $g\in G$. More precisely, if $(A, I)$ is a pro-$G$-graded pair of definition of $\mathcal{A}$, then for every integer $k\geq1$ there exists an integer $N_{k}$ (not depending on $g\in G$) such that for all $n\geq N_{k}$, \begin{equation*}f_{g}-f_{n, g}\in I^{k}\end{equation*}for all $g\in G$.\end{lemma}
\begin{proof}Let $(A_{\gr}, J)$ be a $G$-graded pair of definition of $\mathcal{A}_{\gr}$ such that $(A, I)=(\widehat{A_{\gr}}, J\widehat{A_{\gr}})$. Fix an integer $k\geq1$. Choose a finite subset $G_{k}\subseteq G$ such that \begin{equation*}f-\sum_{g\in G_{k}}f_{g}\in I^{k}.\end{equation*}Choose an integer $N_{k}\geq1$ such that $f-f_{n}\in I^{k}$ for all $n\geq N_{k}$. Moreover, for every $n\geq N_{k}$, choose a finite subset $G_{k, n}$ of $G$ containing $G_{k}$ such that \begin{equation*}f_{n}-\sum_{g\in G_{k, n}}f_{n, g}\in I^{k}.\end{equation*}Then for every $n\geq N_{k}$, \begin{equation*}\sum_{g\in G_{k, n}}(f_{g}-f_{n, g})=(\sum_{g\in G_{k, n}}f_{g}-f)+(f-f_{n})+(f_{n}-\sum_{g\in G_{k, n}}f_{n, g})\in I^{k},\end{equation*}so \begin{equation*}\sum_{g\in G_{k, n}}(f_{g}-f_{n, g})\in I^{k}\cap A_{\gr}=J^{k}.\end{equation*}Since $J^{k}$ is a homogeneous ideal of $A_{\gr}$, we deduce that $f_{g}-f_{n, g}\in I^{k}$ for all $n\geq N_{k}$, as desired. \end{proof}
\begin{lemma}\label{Homogeneous decompositions and morphisms}Let \begin{equation*}(\varphi, \varphi_{\gr}): (\mathcal{A}, \mathcal{A}_{\gr})\to (\mathcal{B}, \mathcal{B}_{\gr})\end{equation*}be a morphism of pro-$G$-graded Huber rings and let $f=\sum_{g\in G}f_{g}$ be the topological homogeneous decomposition of some $f\in\mathcal{A}$. Then $\varphi(f)=\sum_{g\in G}\varphi(f_{g})$ is the topological homogeneous decomposition of $\varphi(f)\in \mathcal{B}$.\end{lemma}
\begin{proof}Follows from $\varphi$ being continuous and $\varphi_{\gr}$ being a homomorphism of $G$-graded rings. \end{proof}
\begin{lemma}\label{Surjective morphisms of pro-graded rings}For a morphism of pro-$G$-graded Huber rings \begin{equation*}(\varphi, \varphi_{\gr}): (\mathcal{A}, \mathcal{A}_{\gr})\to (\mathcal{B}, \mathcal{B}_{\gr}),\end{equation*}if $\varphi$ is surjective, so is $\varphi_{\gr}$.\end{lemma}
\begin{proof}Using Lemma~\ref{Homogeneous decompositions and morphisms}, the proof of \cite{Ishizuka-Yoshikawa25}, Lemma~3.18, applies verbatim.\end{proof}
The notion of topological homogeneous decompositions gives rise to the notion of topologically homogeneous ideals.
\begin{mydef}[Topologically homogeneous ideal]\label{Topologically homogeneous ideal, definition}Let $(\mathcal{A}, \mathcal{A}_{\gr})$ be a pro-$G$-graded Huber ring. We call an ideal $I$ of $\mathcal{A}$ topologically homogeneous if for every element $f\in I$ with topological homogeneous decomposition $f=\sum_{g\in G}f_{g}$ in $\mathcal{A}$ all the components $f_{g}$ belong to $I$. \end{mydef}
We have the following faithful analog of \cite{Ishizuka-Yoshikawa25}, Proposition~3.16(1).
\begin{prop}\label{Detecting homogeneous ideals}Let $(\mathcal{A}, \mathcal{A}_{\gr})$ be a pro-$G$-graded Huber ring. If $J$ is a homogeneous ideal of $\mathcal{A}_{\gr}$ and $f=\sum_{g\in G}f_{g}$ is the topological homogeneous decomposition of an element $f\in J\mathcal{A}$, then $f_{g}\in J$ for all $g\in G$. In particular, every homogeneous ideal of $\mathcal{A}$ is also topologically homogeneous, and if $J$ is a homogeneous ideal of $\mathcal{A}_{\gr}$, then $J\mathcal{A}\cap \mathcal{A}_{\gr}=J$.\end{prop}
\begin{proof}The proof of \cite{Ishizuka-Yoshikawa25}, Proposition~3.15(1), applies verbatim. \end{proof}
Topologically homogeneous ideals behave well with respect to morphisms of pro-graded Huber rings.
\begin{lemma}\label{Graded maps and homogeneous ideals}Let $(\varphi, \varphi_{\gr}): (\mathcal{A}, \mathcal{A}_{\gr})\to (\mathcal{B}, \mathcal{B}_{\gr})$ be a morphism of pro-$G$-graded Huber rings and let $I$ be a topologically homogeneous ideal of $\mathcal{B}$. Then the pre-image $\varphi^{-1}(I)$ is a topologically homogeneous ideal of $\mathcal{A}$.\end{lemma}
\begin{proof}Let $f=\sum_{g\in G}f_{g}$ be the topological homogeneous decomposition of an element $f\in \varphi^{-1}(I)$, i.e., of an element $f\in \mathcal{A}$ with $\varphi(f)\in I$. By Lemma~\ref{Homogeneous decompositions and morphisms}, $\varphi(f)=\sum_{g\in G}\varphi(f_{g})$ is the topological homogeneous decomposition of $\varphi(f)$ in $\mathcal{B}$. Since $I$ is topologically homogeneous, $\varphi(f_{g})\in I$ for every $g\in G$. In other words, $f_{g}\in \varphi^{-1}(I)$ for every $g\in G$. This shows that $\varphi^{-1}(I)$ is topologically homogeneous. \end{proof}
If a topologically homogeneous ideal is closed, then it can be described as the closure of a usual homogeneous ideal. In the statement of the next lemma, a bar over an ideal in a Huber ring denotes topological closure.
\begin{lemma}\label{Closed topologically homogeneous ideals}Let $(\mathcal{A}, \mathcal{A}_{\gr})$ be a pro-$G$-graded Huber ring and let $I$ be a closed ideal of $\mathcal{A}$. The following are equivalent: \begin{enumerate}[(1)]\item $I$ is topologically homogeneous. \item There exists a homogeneous ideal $J$ of $\mathcal{A}_{\gr}$ such that \begin{equation*}I=\overline{J\mathcal{A}}.\end{equation*}In other words, $I$ is the (topological) closure of a homogeneous ideal of $\mathcal{A}$. \item The ideal $I\cap \mathcal{A}_{\gr}$ of $\mathcal{A}_{\gr}$ is homogeneous and \begin{equation*}I=\overline{I\cap \mathcal{A}_{\gr}}.\end{equation*}\end{enumerate}\end{lemma}
\begin{proof}Suppose $I$ is topologically homogeneous. Then, for every $f\in I$ with topologically homogeneous decomposition $f=\sum_{g\in G}f_{g}$, we have $f_{g}\in I\cap \mathcal{A}_{\gr}$ for all $g\in G$. This shows both that $I\cap \mathcal{A}_{\gr}$ is homogeneous and that $I\cap \mathcal{A}_{\gr}$ is dense in $I$, so (1)$\Rightarrow$(3)$\Rightarrow$(2).

Conversely, suppose that $I=\overline{J\mathcal{A}}$ for some homogeneous ideal $J$ of $\mathcal{A}_{\gr}$. Let $f\in I$ have topological homogeneous decomposition $f=\sum_{g\in G}f_{g}$ in $\mathcal{A}$ and let $(f_{n})_{n}$ be a sequence in $J$ converging to $f$ in $\mathcal{A}$ (we are using here that $J$ is dense in $J\mathcal{A}$, which follows from $\mathcal{A}_{\gr}$ being dense in $\mathcal{A}$). For every $n$, let $f_{n}=\sum_{g\in G}f_{n, g}$ be the homogeneous decomposition of $f_{n}$. By Lemma~\ref{Grading and topology 2}, $f_{n, g}\to f_{g}$ for all $g\in G$. Since $J$ is homogeneous, $f_{n, g}\in J$ for all $n$ and all $g\in G$. It follows that $f_{g}$ belongs to the closure of $J$ in $\mathcal{A}$ for all $g\in G$, so $f_{g}\in I$ for all $g\in G$. \end{proof}
We obtain an analog of \cite{Ishizuka-Yoshikawa25}, Proposition~3.17, which we formulate in terms of topologically homogeneous ideals. 
\begin{prop}\label{Graded maps and homogeneous ideals 2}The closed topologically homogeneous ideals in a pro-$G$-graded Huber ring $(\mathcal{A}, \mathcal{A}_{\gr})$ are precisely the kernels of morphisms of pro-$G$-graded Huber rings $(\varphi, \varphi_{\gr}): (\mathcal{A}, \mathcal{A}_{\gr})\to (\mathcal{B}, \mathcal{B}_{\gr})$. For every closed topologically homogeneous ideal $I$ of $\mathcal{A}$, the pair \begin{equation*}(\mathcal{A}/I, \mathcal{A}_{\gr}/(I\cap \mathcal{A}_{\gr}))\end{equation*}is a pro-$G$-graded Huber ring.\end{prop}
\begin{proof}By Lemma~\ref{Graded maps and homogeneous ideals}, the kernel $\ker(\varphi)$ of every morphism $(\varphi, \varphi_{\gr}) $ is topologically homogeneous, and it is closed since pro-$G$-graded Huber rings are (by definition) complete. Conversely, if $I$ is closed and topologically homogeneous, then, by Lemma~\ref{Closed topologically homogeneous ideals}, $I\cap \mathcal{A}_{\gr}$ is a homogeneous ideal of $\mathcal{A}_{\gr}$, so $\mathcal{A}_{\gr}/(I\cap \mathcal{A}_{\gr})$ is a graded ring. Moreover, by Lemma~\ref{Closed topologically homogeneous ideals}, we also know that \begin{equation*}I=\overline{I\cap \mathcal{A}_{\gr}},\end{equation*}so $\mathcal{A}/I$ is the completion of $\mathcal{A}_{\gr}/(I\cap \mathcal{A}_{\gr})$ with respect to the quotient topology.

It remains to show that $\mathcal{A}_{\gr}/(I\cap \mathcal{A}_{\gr})$ is gradedwise complete with respect to the quotient topology. By Lemma~\ref{Grading and topology}, the inclusion map \begin{equation*}\mathcal{A}_{\gr, g}\hookrightarrow \mathcal{A}_{\gr}=\bigoplus_{h\in G}\mathcal{A}_{\gr, h}\end{equation*}has a continuous splitting. This induces a continuous splitting of the inclusion map \begin{equation*}\mathcal{A}_{\gr, g}/(I\cap \mathcal{A}_{\gr, g})\hookrightarrow \mathcal{A}_{\gr}/(I\cap \mathcal{A}_{\gr})=\bigoplus_{h\in G}\mathcal{A}_{\gr, h}/(I\cap \mathcal{A}_{\gr, h}).\end{equation*}It follows that the subspace topology on \begin{equation*}(\mathcal{A}_{\gr}/(I\cap \mathcal{A}_{\gr}))_{g}=\mathcal{A}_{\gr, g}/(I\cap \mathcal{A}_{\gr, g})\end{equation*}coincides with the quotient topology induced by the quotient map \begin{equation*}\mathcal{A}_{\gr, g}\to \mathcal{A}_{\gr, g}/(I\cap \mathcal{A}_{\gr, g}).\end{equation*}But this quotient topology is complete, $I\cap \mathcal{A}_{\gr, g}$ being closed in $\mathcal{A}_{\gr, g}$. \end{proof}
In the same vein, we can also define pro-$G$-graded subrings of pro-$G$-graded Huber rings, a notion which occurs naturally in our study of pro-$G$-graded Huber pairs in Section~\ref{sec:Huber pairs}.
\begin{mydef}[Pro-$G$-graded subring]\label{Pro-graded subring, definition}For a pro-$G$-graded Huber ring $(\mathcal{A}, \mathcal{A}_{\gr})$, we call a subring $A$ of $\mathcal{A}$ a pro-$G$-graded subring if for every element $f\in A$ with topological homogeneous decomposition $f=\sum_{g\in G}f_{g}$ in $\mathcal{A}$ we have $f_{g}\in A$ for all $g\in G$.\end{mydef}
We get the following analog of Lemma~\ref{Closed topologically homogeneous ideals} for closed subrings instead of closed ideals.
\begin{lemma}\label{Closed pro-graded subrings}For a pro-$G$-graded Huber ring $(\mathcal{A}, \mathcal{A}_{\gr})$ and a closed subring $A$ of $\mathcal{A}$, the following are equivalent:\begin{enumerate}[(1)]\item $A$ is a pro-$G$-graded subring of $\mathcal{A}$. \item There exists a $G$-graded subring $A_{\gr}$ of $\mathcal{A}_{\gr}$ such that $A$ is the closure of $A_{\gr}$ inside $\mathcal{A}$. \item $\mathcal{A}_{\gr}\cap A$ is a $G$-graded subring of $\mathcal{A}_{\gr}$ and $A$ is the closure of $\mathcal{A}_{\gr}\cap A$ in $\mathcal{A}$.\end{enumerate}In particular, a subring of $\mathcal{A}$ is a pro-$G$-graded ring of definition of $(\mathcal{A}, \mathcal{A}_{\gr})$ if and only if it is both a ring of definition and a pro-$G$-graded subring of $\mathcal{A}$.\end{lemma}
\begin{proof}Suppose that $A$ is a closed pro-$G$-graded subring of $\mathcal{A}$ and let $f$ be an element of $\mathcal{A}_{\gr}\cap A$ with homogeneous decomposition $f=\sum_{g\in G}f_{g}$ in $\mathcal{A}_{\gr}$. Then $f_{g}\in A$ for all $g\in G$ since $A$ is a pro-$G$-graded subring of $\mathcal{A}$. It follows that $f_{g}\in\mathcal{A}_{\gr}\cap A$ for all $g\in G$, showing that $\mathcal{A}_{\gr}\cap A$ is a $G$-graded subring of $\mathcal{A}_{\gr}$. Moreover, if $f\in A$ is any element of $A$ with topological homogeneous decomposition $f=\sum_{g\in G}f_{g}$ in $\mathcal{A}$, then, by the assumption that $A$ is a pro-$G$-graded subring, $f_{g}\in \mathcal{A}_{\gr}\cap A$ for all $g\in G$ and thus all partial sums of the convergent series $\sum_{g\in G}f_{g}$ belong to $\mathcal{A}_{\gr}\cap A$. This shows that $f$ is in the closure of $\mathcal{A}_{\gr}\cap A$. Since $A$ is closed, we conclude that $A$ is precisely the closure of $\mathcal{A}_{\gr}\cap A$ in $\mathcal{A}$. 

Conversely, suppose that $A_{\gr}$ is a $G$-graded subring of $\mathcal{A}_{\gr}$ such that $A$ is the closure of $A_{\gr}$ in $\mathcal{A}$. Let $f=\sum_{g\in G}f_{g}$ be the topological homogeneous decomposition of an element $f\in A$. Let $(f_{n})_{n}$ be a sequence in $A_{\gr}$ converging to $f$ in $\mathcal{A}$ and let $f_{n}=\sum_{g\in G}f_{n, g}$ be the homogeneous decomposition of each $f_{n}$ in $\mathcal{A}_{\gr}$. By Lemma~\ref{Grading and topology 2}, $f_{n, g}\to f_{g}$ for all $g\in G$. Since $A_{\gr}$ is a $G$-graded subring of $\mathcal{A}_{\gr}$, we have $f_{n, g}\in A_{\gr}$ for all $g\in G$. It follows that every component $f_{g}$ belongs to the closure of $A_{\gr}$ in $\mathcal{A}$, i.e., it belongs to $A$. \end{proof}
We conclude this section with some results on pro-graded perfectoid rings in the sense of Ishizuka-Yoshikawa. This class of pro-graded adic rings was introduced by Ishizuka and Yoshikawa in \cite{Ishizuka-Yoshikawa25}, \S4, and the integral pro-graded perfectoid Tate algebra which we discussed in Example~\ref{Pro-graded perfectoid Tate algebra} is a prominent example of this notion.
\begin{mydef}[Pro-$G$-graded perfectoid rings: Ishizuka-Yoshikawa \cite{Ishizuka-Yoshikawa25}, Definition~4.1, Proposition~4.4(1)$\Leftrightarrow$(2)]~

A pro-$G$-graded perfectoid ring is a pro-$G$-graded adic ring $(R, R_{\gr})$ such that $p$ is topologically nilpotent in $R_{\gr}$ and the $p$-adic completion of $R_{\gr}$ is (integral) perfectoid. Equivalently, $p$ is topologically nilpotent in $R_{\gr}$ and both $R$ and the degree-$0$ part $R_{\gr, 0}$ of $R_{\gr}$ are (integral) perfectoid rings. \end{mydef}
\begin{lemma}\label{Perfect}Let $(R, R_{\gr})$ be a pro-$G$-graded adic ring with $R_{\gr}$ of characteristic $p$. Then $(R, R_{\gr})$ is a pro-$G$-graded perfectoid ring if and only if $R$ is perfect. \end{lemma}
\begin{proof}Since $p=0$ on $R_{\gr}$, the condition that $(R, R_{\gr})$ being a pro-$G$-graded perfectoid ring is equivalent to $R_{\gr}$ being perfect. If the adic ring $R_{\gr}$ is perfect, then its completion $R$ is also perfect. To prove the converse, first note that if $R$ is perfect, then $R_{\gr}$ is reduced since it is a subring of $R$. Moreover, the Frobenius on $R$ is clearly continuous with respect to the adic topology and sends homogeneous elements of degree $g\in G$ to homogeneous elements of degree $pg$. Consequently, Frobenius defines a map of pro-$G$-graded adic rings \begin{equation*}(R, R_{\gr})\to (R, R_{\gr}).\end{equation*}Hence, if the Frobenius is surjective on $R$, it must also be surjective on $R_{\gr}$ by \cite{Ishizuka-Yoshikawa25}, Lemma~3.18. \end{proof}
For every abelian group $G$ satisfying $G=G[1/p]$ and every $G$-graded adic ring $R$, Ishizuka and Yoshikawa introduced a $G$-graded adic ring $R^{\flat,\gr}$ of characteristic $p$ whose degree-$g$ part, for $g\in G$, is given by the inverse limit of \begin{equation*}\{((R/pR)_{g/p^{m}}, F)\}_{m\geq0},\end{equation*}with $F$ being the Frobenius $(R/pR)_{g/p^{m+1}}\to (R/pR)_{g/p^{m}}$, see \cite{Ishizuka-Yoshikawa25}, Example~3.8(1). Moreover, if $(R, R_{\gr})$ is a pro-$G$-graded adic ring such that $R_{\gr}$ is a $G$-graded adic ring with homogeneous ideal of definition $J$ containing some power of $p$, then loc.~cit., Lemma~3.21, ensures that \begin{equation*}(R^{\flat}, R_{\gr}^{\flat,\gr}),\end{equation*}where $R^{\flat}=\varprojlim_{f\mapsto f^{p}}R/pR$ is the usual tilt of $R$ equipped with the inverse limit topology, is a pro-$G$-graded adic ring. In particular, this applies to pro-$G$-graded perfectoid rings. 
\begin{mydef}[Graded tilt]In the above situation, we call $(R^{\flat}, R_{\gr}^{\flat,\gr})$ the graded tilt of the pro-$G$-graded adic ring $(R, R_{\gr})$.\end{mydef}
\begin{lemma}\label{Graded tilt is perfectoid}The graded tilt $(R^{\flat}, R_{\gr}^{\flat,\gr})$ of a pro-$G$-graded perfectoid ring $(R, R_{\gr})$ is a pro-$G$-graded perfectoid ring. \end{lemma}
\begin{proof}By \cite{Ishizuka-Yoshikawa25}, Lemma~3.21, $(R^{\flat}, R_{\gr}^{\flat,\gr})$ is a pro-$G$-graded adic ring. Since $R^{\flat}$ is perfect, it follows that $(R^{\flat}, R_{\gr}^{\flat,\gr})$ is a pro-$G$-graded perfectoid ring by Lemma~\ref{Perfect}.  \end{proof}
Similarly, we can define the tilt of a $G$-graded continuous ring homomorphism.
\begin{mydef}[Graded tilt of a morphism]\label{Tilt of a morphism}Suppose that $G=G[1/p]$ and let \begin{equation*}\varphi: R\to S\end{equation*}be a morphism of $G$-graded adic rings. We define the graded tilt of $\varphi$ to be the continuous $G$-graded ring homomorphism \begin{equation*}\varphi^{\flat,\gr}: R^{\flat,\gr}\to S^{\flat,\gr}\end{equation*}whose degree-$g$ part is given by \begin{equation*}(\varphi^{\flat,\gr})_{g}((f_{m})_{m}):=(\overline{\varphi_{g/p^{m}}}(f_{m}))_{m}\end{equation*}for any \begin{equation*}(f_{m})_{m}\in (R^{\flat,\gr})_{g}=\varprojlim_{F}(R/pR)_{g/p^{m}},\end{equation*}where $\overline{\varphi_{h}}$ denotes the mod $p$ reduction of $\varphi_{h}$. \end{mydef}
We have the following analog of \cite{Ishizuka-Yoshikawa25}, Lemma~3.21, for morphisms.
\begin{lemma}\label{Graded morphisms and tilting}Suppose that $G=G[1/p]$ and let \begin{equation*}(\varphi, \varphi_{\gr}): (R, R_{\gr})\to (S, S_{\gr})\end{equation*}be a morphism of pro-$G$-graded adic rings such that $p$ is topologically nilpotent in $R_{\gr}$ (and hence also in $S_{\gr}$). Then the map on completions \begin{equation*}\widehat{\varphi_{\gr}^{\flat,\gr}}: \widehat{R_{\gr}^{\flat,\gr}}\to \widehat{S_{\gr}^{\flat,\gr}},\end{equation*}where $R_{\gr}^{\flat,\gr}$ and $S_{\gr}^{\flat,\gr}$ are endowed with their respective inverse limit topologies, coincides with the usual tilt $\varphi^{\flat}$ of $\varphi$. \end{lemma}
\begin{proof}Recall that the usual tilt \begin{equation*}\varphi^{\flat}: R^{\flat}=\varprojlim_{f\mapsto f^{p}}R/pR\to S^{\flat}=\varprojlim_{f\mapsto f^{p}}S/pS\end{equation*}is given by \begin{equation*}\varphi^{\flat}((f_{m})_{m})=(\overline{\varphi}(f_{m}))_{m}.\end{equation*}In particular, $\varphi^{\flat}$ is continuous for the inverse limit topologies on $R^{\flat}$ and $S^{\flat}$. On the other hand, $R_{\gr}^{\flat,\gr}$ is a dense subring of $R^{\flat}$, by \cite{Ishizuka-Yoshikawa25}, Lemma~3.21. Therefore, it suffices to prove that the restriction of $\varphi^{\flat}$ to $R_{\gr}^{\flat,\gr}$ coincides with $\varphi_{\gr}^{\flat,\gr}$.

By the proof of \cite{Ishizuka-Yoshikawa25}, Lemma~3.21, the completion of $R_{\gr}^{\flat,\gr}$ for its inverse limit topology is canonically identified with $\varprojlim_{f\mapsto f^{p}}R/\overline{pR}$ (where the bar denotes topological closure in $R$), which is in turn canonically identified with $R^{\flat}=\varprojlim_{f\mapsto f^{p}}R/pR$ since $p$ is topologically nilpotent in $R$. The first identification comes from the canonical ring isomorphism (``degree shift'') between \begin{equation*}\bigoplus_{g\in G}R_{g/p^{m}}/pR_{g/p^{m}}\end{equation*}and $R_{\gr}/pR_{\gr}=R/pR$. Since the map $\overline{\varphi}$ (the mod $p$ reduction of $\varphi$) does not take into account degrees, it is invariant under composition with the above canonical isomorphism. Hence the assertion that $\varphi^{\flat}$ coincides with $\varphi_{\gr}^{\flat,\gr}$ on the dense subring $R_{\gr}^{\flat,\gr}$ of $R^{\flat}$ follows directly from the definitions of these two maps.  \end{proof}
The proof of the following proposition was communicated to us by Ryo Ishizuka. 
\begin{prop}[Ishizuka, private communication]\label{Pro-graded perfectoid Tate algebra and tilting}Let $K$ be a perfectoid field. Then the graded tilt of the integral pro-graded perfectoid Tate algebra over $K$ is given by the integral pro-graded perfectoid Tate algebra over $K^{\flat}$.\end{prop}
\begin{proof}Let $(R, R_{\gr})=(T_{n, K}^{\perfd\circ}, (T_{n, K}^{\perfd\circ})_{\gr})$ and let $p^{\flat}\in K^{\flat\circ}$ be a compatible system of $p$-power roots of some unit multiple of $p$ in $K^{\circ}$. Recall that $R^{\flat}$ is $p^{\flat}$-adically complete. Fix a degree $d\in \mathbb{Z}[1/p]$. By \cite{Ishizuka-Yoshikawa25}, Lemma~3.21, $(R^{\flat}, R_{\gr}^{\flat,\gr})$ is a pro-$G$-graded perfectoid ring. In particular, $R_{\gr}^{\flat,\gr}$ is gradedwise complete. By Lemma~\ref{Degree-0 part}, this implies that the degree-$d$ part $(R_{\gr}^{\flat,\gr})_{d}$ of $R_{\gr}^{\flat,\gr}$ is $p^{\flat}$-adically complete. For a multiindex $i=(i_1,\dots, i_n)\in \mathbb{Z}[1/p]^{n}$, we write \begin{equation*}X^{i}=X_{1}^{i_{1}}\cdot~\dots~\cdot X_{n}^{i_{n}}.\end{equation*}We also let \begin{equation*}t_{1}=(X_{1}^{1/p^{m}})_{m\geq0},\dots, t_{n}=(X_{n}^{1/p^{m}})_{m\geq0}\in R^{\flat}\end{equation*}and write \begin{equation*}t^{i}=t_{1}^{i_{1}}\cdot~\dots~\cdot t_{n}^{i_{n}}.\end{equation*}Consider the map \begin{equation*}\iota: \bigoplus_{i_1+\dots+i_{n}=d}K^{\flat\circ}\cdot t^{i}\to (R_{\gr}^{\flat,\gr})_{d}\end{equation*}given by \begin{align*}\bigoplus_{i_1+\dots+i_{n}=d}K^{\flat\circ}\cdot t^{i}\cong \bigoplus_{i_{1}+\dots+i_{n}=d}\varprojlim_{m\geq0} ( \dots \xrightarrow{F} (K^{\circ}/p)\cdot X^{i/p} \xrightarrow{F} (K^{\circ}/p)\cdot X^{i}) \\\to \varprojlim_{m\geq0}(\dots \xrightarrow{F}\bigoplus_{i_1+\dots+i_n=d}(K^{\circ}/p)\cdot X^{i/p}\xrightarrow{F} \bigoplus_{i_1+\dots+i_n=d}(K^{\circ}/p)\cdot X^{i}) \\ \cong \varprojlim_{m\geq0}(\dots \xrightarrow{F} \bigoplus_{i_1+\dots+i_n=\frac{d}{p}}(K^{\circ}/p) \cdot X^{i} \xrightarrow{F} \bigoplus_{i_1+\dots+i_n=d}(K^{\circ}/p) \cdot X^{i})=(R_{\gr}^{\flat,\gr})_{d}, \end{align*}where the arrow in the middle is the map determined by the inclusion maps \begin{equation*}(K^{\circ}/p)\cdot X^{i}\hookrightarrow \bigoplus_{i_{1}+\dots+i_{n}=d}(K^{\circ}/p)\cdot X^{i}.\end{equation*}We want to prove that $\iota$ induces an isomorphism upon taking $p^{\flat}$-adic completions. By topological Nakayama's lemma and a standard dévissage argument for the non-zero-divisor $p^{\flat}$ (see, for example, \cite{Dine25-3}, Lemma~5.2), it suffices to prove that $\iota$ induces an isomorphism modulo $p^{\flat}$. 

To prove this, recall from \cite{Ishizuka-Yoshikawa25}, Lemma~3.21, that the restriction to $R_{\gr}^{\flat,\gr}$ of the multiplicative map \begin{equation*}\sharp: R^{\flat}=\varprojlim_{m\geq0}(\dots \xrightarrow{F} R/pR \xrightarrow{F} R/pR) \to R, (\overline{f_{m}})_{m\geq0}\mapsto \lim_{m\to \infty}f_{m}^{p^{m}},\end{equation*}coincides with the multiplicative map\begin{equation*}\sharp_{\gr}: R_{\gr}^{\flat,\gr}\to R_{\gr}\end{equation*}which is defined by the analogous formula in every degree. It follows that the isomorphism of rings \begin{equation*}R_{\gr}^{\flat,\gr}/p^{\flat}R_{\gr}^{\flat,\gr}=R^{\flat}/p^{\flat}R^{\flat}\xrightarrow{\simeq}R/pR=R_{\gr}/pR_{\gr}\end{equation*}induced by $\sharp$ is actually an isomorphism of graded rings. In our case, this means that the $K^{\circ}/p=K^{\flat\circ}/p^{\flat}$-linear isomorphism \begin{equation*}(K^{\flat\circ}/p^{\flat})[t_{1}^{1/p^{\infty}},\dots, t_{n}^{1/p^{\infty}}]\xrightarrow{\simeq}(K^{\circ}/p)[X_{1}^{1/p^{\infty}},\dots, X_{n}^{1/p^{\infty}}]\end{equation*}which sends $t_{1}^{1/p^{k}},\dots, t_{n}^{1/p^{k}}$ to $X_{1}^{1/p^{k}},\dots, X_{n}^{1/p^{k}}$ for every $k\geq0$ restricts to a $K^{\circ}/p=K^{\flat\circ}/p^{\flat}$-linear isomorphism \begin{equation*}\bigoplus_{i_1+\dots+i_n=d}(K^{\flat\circ}/p^{\flat})\cdot t^{i}\xrightarrow{\simeq} \bigoplus_{i_1+\dots+i_n=d}(K^{\circ}/p)\cdot X^{i}\end{equation*}for every degree $d$. On the other hand, the map $\iota$ modulo $p^{\flat}$ is clearly $K^{\circ}/p$-linear and sends $t_{1}^{1/p^{k}},\dots, t_{n}^{1/p^{k}}$ to $X_{1}^{1/p^{k}},\dots, X_{n}^{1/p^{k}}$ for every $k\geq0$, so it must coincide with this isomorphism in degree $d$. This proves that $\iota$ is an isomorphism mod $p^{\flat}$, finishing the proof of the proposition. \end{proof}
Finally, we conclude this section by checking that the notion of topologically homogeneous ideals is compatible with tilting. Recall that an ideal $I$ in an (integral) perfectoid ring $R$ is called perfectoid if the quotient $R/I$ is again a perfectoid ring. At the level of Tate rings, the perfectoid ideals of $R$ correspond to spectrally reduced ideals of the perfectoid Tate ring $\mathcal{R}=R[\varpi^{-1}]$ (where $\varpi$ is a perfectoid element of $R$ admitting a compatible system of $p$-power roots): These are the ideals of $\mathcal{R}$ which can be described as kernels of continuous power-multiplicative seminorms, see \cite{Dine22}, Definition~2.16, and \cite{Dine-Ishizuka23}, Definition~2.4. By a theorem of Bhatt (\cite{Kedlaya17}, Theorem~2.9.12), these are precisely the ideals $I$ of $\mathcal{R}$ such that the quotient $\mathcal{R}/\mathcal{J}$ is a perfectoid Tate ring. 
\begin{lemma}\label{Quotients by perfectoid ideals}Suppose that $G=G[1/p]$. Let $(R, R_{\gr})$ be a pro-$G$-graded perfectoid ring and let $I$ be a perfectoid ideal of $R$ which is also a topologically homogeneous ideal of the pro-$G$-graded perfectoid ring $(R, R_{\gr})$. Then $(R/I, R_{\gr}/(I\cap R_{\gr}))$ is a pro-$G$-graded perfectoid ring.\end{lemma}
\begin{proof}Since $I$ is a perfectoid ideal, it is in particular closed with respect to the given adic topology on $R$ (as well as with respect to the $p$-adic topology) and thus $I\cap R_{\gr}$ is also closed in $R_{\gr}$. By \cite{Ishizuka-Yoshikawa25}, Proposition~3.17(2), the pair $(R/I, R_{\gr}/(I\cap R_{\gr}))$ is a pro-$G$-graded adic ring. Since $I$ is a perfectoid ideal of $R$, the quotient $R/I$ is perfectoid, so it remains to prove that the degree-$0$ part \begin{equation*}(R_{\gr}/(I\cap R_{\gr}))_{0}=R_{\gr,0}/(I\cap R_{\gr,0})\end{equation*}is also perfectoid. Since $(R, R_{\gr})$ was a pro-$G$-graded perfectoid ring, $R_{\gr,0}$ is a perfectoid ring. Since $I$ is a perfectoid ideal, it is $p$-adically closed in $R$, so $I\cap R_{\gr,0}$ is a $p$-adically closed ideal of $R_{\gr,0}$. It follows that $R_{\gr,0}/(I\cap R_{\gr,0})$ is a semiperfectoid ring. But we have an injective ring map \begin{equation*}R_{\gr,0}/(I\cap R_{\gr,0})\hookrightarrow R/I\end{equation*}into the perfectoid ring $R/I$, so the semiperfectoid ring $R_{\gr,0}/(I\cap R_{\gr,0})$ must be perfectoid by \cite{Dine-Ishizuka23}, Lemma~2.3. \end{proof}
For perfectoid ideals $I$ in a perfectoid ring $R$, we recall the notion of the tilt $I^{\flat}$ of the ideal $I$ introduced in \cite{Dine22}, \S4, and \cite{Dine-Ishizuka23}, Definition~3.2. Let $I$ be a perfectoid ideal in a perfectoid ring $R$, or a spectrally reduced ideal in a perfectoid Tate ring $\mathcal{R}$. Then the tilt $I^{\flat}$ of $I$ is the subset \begin{equation*}I^{\flat}=\{\, f\in R^{\flat}\mid f^{\#}\in I\,\}\end{equation*}of $R^{\flat}$ (respectively, the subset \begin{equation*}I^{\flat}=\{\, f\in \mathcal{R}^{\flat}\mid f^{\#}\in I\,\}\end{equation*}of $\mathcal{R}^{\flat}$). By \cite{Dine-Ishizuka23}, Lemma~3.3 (respectively, by \cite{Dine22}, Lemma~4.8) this is a $\varpi^{\flat}$-adically closed radical ideal of $R^{\flat}$ (respectively, a closed radical ideal of $\mathcal{R}^{\flat}$) such that \begin{equation*}(R/I)^{\flat}=R^{\flat}/I^{\flat} \end{equation*}(respectively, $(\mathcal{R}/I)^{\flat}=\mathcal{R}^{\flat}/I^{\flat}$). Moreover, by \cite{Dine-Ishizuka23}, Theorem~3.5 (respectively, by \cite{Dine22}, Proposition~4.9), the map $I\mapsto I^{\flat}$ defines an inclusion-preserving bijection between the set of perfectoid ideals of $R$ (respectively, of spectrally reduced ideals of $\mathcal{R}$) and the set of closed radical ideals of $R^{\flat}$ or $\mathcal{R}^{\flat}$. We can now check that this tilting operation for ideals is compatible with pro-graded structures. 
\begin{lemma}\label{Homogeneous ideals and tilting}Suppose that $G=G[1/p]$. Let $(R, R_{\gr})$ be a pro-$G$-graded perfectoid ring and let $I$ be a perfectoid ideal of $R$ which is also a topologically homogeneous ideal of the pro-$G$-graded perfectoid ring $(R, R_{\gr})$. Then the tilt $I^{\flat}$ of $I$ is a topologically homogeneous ideal of the pro-$G$-graded perfectoid ring $(R^{\flat}, R_{\gr}^{\flat,\gr})$. \end{lemma}
\begin{proof}By Lemma~\ref{Quotients by perfectoid ideals}, the pair $(R/I, R_{\gr}/(I\cap R_{\gr}))$ is a pro-$G$-graded perfectoid ring. By \cite{Dine-Ishizuka23}, Lemma~3.3, $(R/I)^{\flat}=R^{\flat}/I^{\flat}$, so, in particular, the map $R^{\flat}\to (R/I)^{\flat}$ is surjective. Therefore, by tilting, we obtain a homomorphism of pro-$G$-graded perfectoid rings \begin{equation*}(R^{\flat}, R_{\gr}^{\flat,\gr})\to ((R/I)^{\flat}, (R_{\gr}/(I\cap R_{\gr}))^{\flat,\gr})\end{equation*}whose underlying map of perfectoid rings $R^{\flat}\to (R/I)^{\flat}$ is surjective. By \cite{Ishizuka-Yoshikawa25}, Lemma~3.18, this implies that the map of $G$-graded rings \begin{equation*}R_{\gr}^{\flat,\gr}\to (R_{\gr}/(I\cap R_{\gr}))^{\flat,\gr}\end{equation*}is also surjective. Let $J$ be the kernel of this latter map of $G$-graded rings, a homogeneous ideal of $R_{\gr}^{\flat,\gr}$. Since $R^{\flat}\to (R/I)^{\flat}$ is the $\varpi^{\flat}$-adic completion of the surjective map $R_{\gr}^{\flat,\gr}\to (R_{\gr}/(I\cap R_{\gr}))^{\flat,\gr}$, its kernel $I^{\flat}$ is the (topological) closure of $J$. This shows that $I^{\flat}$ is indeed a topologically homogeneous ideal of $(R^{\flat}, R_{\gr}^{\flat,\gr})$, by Lemma~\ref{Closed topologically homogeneous ideals}. \end{proof}

\section{Pro-graded Huber pairs and topologically homogeneous valuations}\label{sec:Huber pairs}

As in the previous section, we fix an abelian group $G$. In most of the statements in this section we will however put additional assumptions on $G$ such as $G$ being torsion-free or $G$ being a totally ordered abelian group. Our aim in this section is to introduce and study pro-graded analogs of the notions of a ring of integral elements and of a Huber pair.
\begin{mydef}[Pro-graded ring of integral elements]\label{Pro-graded ring of integral elements, definition}Let $(\mathcal{A}, \mathcal{A}_{\gr})$ be a pro-$G$-graded Huber ring. A subring $\mathcal{A}^{+}$ of $\mathcal{A}$ is called a pro-$G$-graded ring of integral elements if it is both a ring of integral elements of $\mathcal{A}$ and a pro-$G$-graded subring of $\mathcal{A}$ (in the sense of Definition~\ref{Pro-graded subring, definition}), and if it contains some pro-$G$-graded ring of definition of $(\mathcal{A}, \mathcal{A}_{\gr})$.\end{mydef} 
\begin{mydef}[Pro-graded Huber pair]\label{Pro-graded Huber pair, definition}A pro-$G$-graded (Tate) Huber pair is a triple \begin{equation*}(\mathcal{A}, \mathcal{A}_{\gr}, \mathcal{A}^{+})=((\mathcal{A}, \mathcal{A}_{\gr}), \mathcal{A}^{+}),\end{equation*}where $(\mathcal{A}, \mathcal{A}_{\gr})$ is a pro-$G$-graded Huber ring (respectively, a pro-$G$-graded Tate ring) and where $\mathcal{A}^{+}$ is a pro-$G$-graded ring of integral elements of $(\mathcal{A}, \mathcal{A}_{\gr})$. 

We call the (Tate) Huber pair $(\mathcal{A}, \mathcal{A}^{+})$ the underlying (Tate) Huber pair of the pro-$G$-graded (Tate) Huber pair $(\mathcal{A}, \mathcal{A}_{\gr}, \mathcal{A}^{+})$. \end{mydef}
\begin{rmk}\label{Two definitions of pro-graded Huber pairs}Note that, by virtue of Lemma~\ref{Closed pro-graded subrings}, for a ring of integral elements $\mathcal{A}^{+}$ of $\mathcal{A}$ to be a pro-$G$-graded ring of integral elements, it is equivalent to require that \begin{equation*}\mathcal{A}_{\gr}^{+}=\mathcal{A}_{\gr}\cap \mathcal{A}^{+}\end{equation*}is a $G$-graded subring of $\mathcal{A}_{\gr}$ which contains some $G$-graded ring of definition of $\mathcal{A}_{\gr}$ (note that, since $\mathcal{A}^{+}$ is an open subring and $\mathcal{A}_{\gr}$ is a dense subring of $\mathcal{A}$, the ring $\mathcal{A}^{+}$ is automatically the closure of $\mathcal{A}_{\gr}^{+}$ in $\mathcal{A}$). We freely use this equivalent characterization of pro-$G$-graded rings of integral elements in what follows.\end{rmk}
\begin{mydef}We call a pro-$G$-graded (Tate) Huber pair uniform (respectively, perfectoid) if the underlying (Tate) Huber pair is uniform (respectively, perfectoid).\end{mydef}
We also have an obvious notion of morphisms between pro-graded Huber pairs. 
\begin{mydef}\label{Morphisms of pro-graded Huber pairs}A morphism of pro-$G$-graded Tate Huber pairs \begin{equation*}(\mathcal{A}, \mathcal{A}_{\gr}, \mathcal{A}^{+})\to (\mathcal{B}, \mathcal{B}_{\gr}, \mathcal{B}^{+})\end{equation*}is a morphism of pro-$G$-graded Huber rings \begin{equation*}(\varphi, \varphi_{\gr}): (\mathcal{A}, \mathcal{A}_{\gr})\to (\mathcal{B}, \mathcal{B}_{\gr})\end{equation*}such that $\varphi(\mathcal{A}^{+})\subseteq \mathcal{B}^{+}$. \end{mydef}
For the purposes of the adic $\Proj$ construction, we will often restrict our attention to pro-graded Huber pairs which are ``positively pro-graded''. Recall that a totally ordered abelian group is an abelian group $G$ equipped with a total order compatible with addition, i.e., $g\leq g^{\prime}$ implies $g+h\leq g^{\prime}+h$ for all $g, g^{\prime}, h\in G$. Note that a totally ordered abelian group is automatically torsion-free.  
\begin{mydef}[Pro-$G_{\geq0}$-graded Huber pair]\begin{enumerate}[(1)]\item For a submonoid $H$ of an abelian group $G$, a $G$-graded ring $\mathcal{A}$ is called $H$-graded if $\mathcal{A}_{g}=0$ for all $g\in G\setminus H$. A pro-$G$-graded Huber ring $(\mathcal{A}, \mathcal{A}_{\gr})$ (respectively, a pro-$G$-graded Huber pair $((\mathcal{A}, \mathcal{A}_{\gr}), \mathcal{A}^{+})$) is called pro-$H$-graded if the $G$-graded ring $\mathcal{A}_{\gr}$ is $H$-graded. \item In particular, if $G$ is a totally ordered abelian group, a pro-$G_{\geq0}$-graded Huber ring $(\mathcal{A}, \mathcal{A}_{\gr})$ is a pro-$G$-graded Huber ring such that $\mathcal{A}_{\gr,g}=0$ for all $g<0$, and a pro-$G_{\geq0}$-graded Huber pair is a pro-$G$-graded Huber pair whose underlying pro-$G$-graded Huber ring is pro-$G_{\geq0}$-graded.\end{enumerate}\end{mydef}  
The following rather elementary lemma is needed to define the notion of rational localization for pro-graded Huber pairs.
\begin{lemma}\label{Integral closure and grading}Let $G$ be a torsion-free abelian group and let $R$ be a $G$-graded subring of a $G$-graded ring $S$. Then the integral closure of $R$ in $S$ is again a $G$-graded subring of $S$.\end{lemma}
\begin{proof}If the torsion-free abelian group $G$ is finitely generated (and hence isomorphic to some $\mathbb{Z}^{e}$), this is proved in \cite{Swanson-Huneke}, Theorem~2.3.2. We reduce the general assertion to this result. Let $f\in S$ be an element of $S$ which is integral over $R$ and let $f=\sum_{g\in G}f_{g}$ be the homogeneous decomposition of $f$. We have to show that $f_{g}$ is integral over $R$ for every $g\in G$. Since $H=\{g\in G\mid f_{g}\neq 0\}$ is finite, it is contained in some finitely generated subgroup $G^{\prime}$ of $G$. Since $G^{\prime}$ is still torsion-free, it is isomorphic to $\mathbb{Z}^{e}$ for some $e\geq0$, so, by \cite{Swanson-Huneke}, Theorem~2.3.2, we know that the integral closure of $R\cap S^{\prime}$ in the $G^{\prime}$-graded ring \begin{equation*}S^{\prime}=\bigoplus_{g\in G^{\prime}}S_{g}\end{equation*}is a $G^{\prime}$-graded subring of $S$. Since $f\in S^{\prime}$, it follows that $f_{g}$ is integral over $R$ for every $g\in G^{\prime}$. This concludes the proof. \end{proof}
Recall from Definition~\ref{Pro-graded rational localizations} the notion of pro-graded rational localizations of pro-graded Huber rings. We now introduce the corresponding notion for pro-graded Huber pairs.
\begin{lemma}\label{Pro-graded rational localizations for Huber pairs}Suppose that $G$ is torsion-free. Let $((\mathcal{A}, \mathcal{A}_{\gr}), \mathcal{A}^{+})$ be a pro-$G$-graded Huber pair and let $f_1,\dots, f_{n}, f_{n+1}\in \mathcal{A}_{\gr}^{+}$ be homogeneous elements generating an open ideal of $\mathcal{A}_{\gr}$. Then \begin{equation*}((\mathcal{B}, \mathcal{B}_{\gr}), \mathcal{B}^{+})=((\mathcal{A}\langle \frac{f_1,\dots, f_n}{f_{n+1}}\rangle, \mathcal{A}\langle \frac{f_1,\dots, f_n}{f_{n+1}}\rangle_{\gr}), \mathcal{A}\langle \frac{f_1,\dots, f_n}{f_{n+1}}\rangle^{+}),\end{equation*}where $\mathcal{A}\langle \frac{f_1,\dots, f_n}{f_{n+1}}\rangle^{+}$ is the integral closure of $\mathcal{A}^{+}\langle \frac{f_1,\dots, f_n}{f_{n+1}}\rangle$ inside $\mathcal{A}\langle \frac{f_1,\dots, f_n}{f_{n+1}}\rangle$, is a pro-$G$-graded Huber pair. \end{lemma}
\begin{proof}In view of Lemma~\ref{Lemma on rational localizations}, it is clear that $\mathcal{B}^{+}$ contains a pro-$G$-graded ring of definition of $(\mathcal{B}, \mathcal{B}_{\gr})$. We have to show that
\begin{equation*}\mathcal{B}_{\gr}^{+}=\mathcal{B}_{\gr}\cap \mathcal{B}^{+}\end{equation*}is a $G$-graded subring of $\mathcal{B}_{\gr}$. Note that $\mathcal{B}_{\gr}^{+}$ is the integral closure of \begin{equation*}B=\mathcal{B}_{\gr}\cap \mathcal{A}^{+}\langle\frac{f_1,\dots, f_{n}}{f_{n+1}}\rangle\end{equation*}(consider an equation of integral dependence of $f\in \mathcal{B}_{\gr}^{+}$ over $\mathcal{A}^{+}\langle\frac{f_1,\dots, f_{n}}{f_{n+1}}\rangle$ and use uniqueness of topological homogeneous decompositions to deduce that the coefficients must be in $\mathcal{B}_{\gr}$).  Hence Lemma~\ref{Integral closure and grading} reduces our task to proving that $B$ is a $G$-graded subring. Let $f\in B$ with homogeneous decomposition $f=\sum_{g\in G}f_{g}$ in $\mathcal{B}_{\gr}$. We have to show that $f_{g}\in B$ for all $g\in G$. If \begin{equation*}(f_{k})_{k}\subseteq \mathcal{A}_{\gr}^{+}[\frac{f_1,\dots, f_n}{f_{n+1}}]\end{equation*}is a sequence converging to $f$ in $B$, then, by Lemma~\ref{Grading and topology}, the homogeneous components $f_{g, k}$ of $f_{k}$ converge to the homogeneous components $f_{g}$. Since $B$ is open (and thus also closed) in $\mathcal{B}_{\gr}$, it then suffices to assume that $f\in \mathcal{A}_{\gr}^{+}[\frac{f_1,\dots, f_n}{f_{n+1}}]$. Thus we are reduced to proving that \begin{equation*}\mathcal{A}_{\gr}^{+}[\frac{f_1,\dots, f_n}{f_{n+1}}]\end{equation*}is a $G$-graded subring of $\mathcal{A}_{\gr}[\frac{1}{f_{n+1}}]$. But, by the assumption that $((\mathcal{A}, \mathcal{A}_{\gr}), \mathcal{A}^{+})$ is a pro-$G$-graded Huber pair, we know that $\mathcal{A}_{\gr}^{+}=\mathcal{A}_{\gr}\cap \mathcal{A}^{+}$ is a $G$-graded subring of $\mathcal{A}_{\gr}$. Hence the desired assertion follows from Lemma~\ref{Lemma on rational localizations}.\end{proof}
\begin{mydef}[Pro-graded rational localizations of pro-graded Huber pairs]For any pro-$G$-graded (Tate) Huber pair $((\mathcal{A}, \mathcal{A}_{\gr}), \mathcal{A}^{+})$, we call pro-$G$-graded (Tate) Huber pairs $((\mathcal{B}, \mathcal{B}_{\gr}), \mathcal{B}^{+})$ constructed from $((\mathcal{A}, \mathcal{A}_{\gr}), \mathcal{A}^{+})$ as in Lemma~\ref{Pro-graded rational localizations for Huber pairs} the (pro-graded) rational localizations of $((\mathcal{A}, \mathcal{A}_{\gr}), \mathcal{A}^{+})$.\end{mydef}
We can also construct quotients of pro-graded Huber pairs, which relies on Lemma~\ref{Integral closure and grading} and the following result on topologically homogeneous ideals. 
\begin{lemma}\label{Saturated homogeneous ideals}Let $((\mathcal{A}, \mathcal{A}_{\gr}), \mathcal{A}^{+})$ be a pro-$G$-graded Huber pair. If $I$ is a closed topologically homogeneous ideal of $\mathcal{A}$ and $A$ is an open pro-$G$-graded subring of $\mathcal{A}$, then $I\cap A$ is a closed topologically homogeneous ideal of the pro-$G$-graded Huber ring $(A, \mathcal{A}_{\gr}\cap A)$. In particular, if $I$ is a closed topologically homogeneous ideal of $\mathcal{A}$, then $I\cap \mathcal{A}^{+}$ is a closed topologically homogeneous ideal of $(\mathcal{A}^{+}, \mathcal{A}_{\gr}^{+})$. \end{lemma}
\begin{proof}Recall from Lemma~\ref{Closed pro-graded subrings} that $A_{\gr}=\mathcal{A}_{\gr}\cap A$ is a $G$-graded subring of $\mathcal{A}_{\gr}$ and $A$ is the closure of $A_{\gr}$, so $(A, A_{\gr})$, with the subspace topology on $A$, is indeed a pro-$G$-graded Huber ring. We use the characterization of closed topologically homogeneous ideals given in Lemma~\ref{Closed topologically homogeneous ideals}. Let $J$ be a homogeneous ideal of $\mathcal{A}_{\gr}$ with $I=\overline{J\mathcal{A}}$. Let $f=\sum_{g\in G}f_{g}$ be the homogeneous decomposition of some element $f\in J\cap A_{\gr}=J\cap A$. Since $J$ is homogeneous, $f_{g}\in J$ for all $g\in G$. Since $A_{\gr}$ is a $G$-graded subring of $\mathcal{A}_{\gr}$, we also have $f_{g}\in A_{\gr}$. This shows that $J\cap A_{\gr}$ is a homogeneous ideal of $A_{\gr}$. On the other hand, if $f\in I\cap A$ and $(f_{n})_{n}$ is a sequence of elements in $J$ with $f_{n}\to f$ (note that $J$ is dense in $J\mathcal{A}$ since $\mathcal{A}_{\gr}$ is dense in $\mathcal{A}$), then, for $n$ large enough, we also have $f_{n}\in A$, since $A$ is an open subring of $\mathcal{A}$. It follows that $I\cap A$ is the topological closure in $A$ of the homogeneous ideal $J\cap A=J\cap A_{\gr}$, proving that $I\cap A$ is topologically homogeneous. \end{proof}
\begin{lemma}\label{Quotients of pro-graded Huber pairs 0}Let $G$ be torsion-free, let $((\mathcal{A}, \mathcal{A}_{\gr}), \mathcal{A}^{+})$ be a pro-$G$-graded Huber pair and let $I$ be a closed topologically homogeneous ideal of $\mathcal{A}$. If $A$ is an open pro-$G$-graded subring of $\mathcal{A}$ and $A_{\gr}=\mathcal{A}_{\gr}\cap A$, then $(A/(I\cap A), A_{\gr}/(I\cap A_{\gr}))$ is a pro-$G$-graded Huber ring and $A/(I\cap A)$ is an open pro-$G$-graded subring of the pro-$G$-graded Huber ring $(\mathcal{A}/I, \mathcal{A}_{\gr}/(I\cap \mathcal{A}_{\gr}))$.\end{lemma}
\begin{proof}We know from Proposition~\ref{Graded maps and homogeneous ideals 2} that $(\mathcal{A}/I, \mathcal{A}_{\gr}/(I\cap \mathcal{A}_{\gr}))$ is indeed a pro-$G$-graded Huber ring. By Lemma~\ref{Closed pro-graded subrings}, $A_{\gr}=\mathcal{A}_{\gr}\cap A$ is a $G$-graded subring of $\mathcal{A}_{\gr}$ and $(A, A_{\gr})$ is a pro-$G$-graded Huber ring. By Lemma~\ref{Saturated homogeneous ideals}, $I\cap A$ is a closed topologically homogeneous ideal of $(A, A_{\gr})$. Applying Proposition~\ref{Graded maps and homogeneous ideals 2} once again, we find that \begin{equation*}(A/(I\cap A), A_{\gr}/(I\cap A_{\gr}))\end{equation*}is a pro-$G$-graded Huber ring (when $A/(I\cap A)$ is equipped with the quotient topology). In particular, $A_{\gr}/(I\cap A_{\gr})$ is a $G$-graded subring of $\mathcal{A}_{\gr}/(I\cap \mathcal{A}_{\gr})$. Since $A$ is an open subring of $\mathcal{A}$, the quotient topology on $A/(I\cap A)$ coincides with the subspace topology induced from $\mathcal{A}/I$ and $A/(I\cap A)$ is an open subring of $\mathcal{A}/I$. We conclude by Lemma~\ref{Closed pro-graded subrings} that $A/(I\cap A)$ is a pro-$G$-graded subring of $\mathcal{A}/I$.  \end{proof}
\begin{lemma}\label{Quotients of pro-graded Huber pairs}Let $((\mathcal{A}, \mathcal{A}_{\gr}), \mathcal{A}^{+})$ be a pro-$G$-graded Huber pair and let $I$ be a closed topologically homogeneous ideal of $\mathcal{A}$. Then \begin{equation*}((\mathcal{A}/I, \mathcal{A}_{\gr}/(I\cap \mathcal{A}_{\gr})), (\mathcal{A}/I)^{+}),\end{equation*}where $(\mathcal{A}/I)^{+}$ is the integral closure of $\mathcal{A}^{+}/(I\cap \mathcal{A}^{+})$ in $\mathcal{A}/I$, is a pro-$G$-graded Huber pair.\end{lemma}
\begin{proof}By Proposition~\ref{Graded maps and homogeneous ideals 2}, $(\mathcal{A}/I, \mathcal{A}_{\gr}/(I\cap \mathcal{A}_{\gr}))$ is a pro-$G$-graded Huber ring. We have to prove that $(\mathcal{A}/I)^{+}$ is a pro-$G$-graded subring of $\mathcal{A}/I$ and contains a pro-$G$-graded ring of definition of $\mathcal{A}/I$. Let $A$ be a pro-$G$-graded ring of definition of $\mathcal{A}$ contained in $\mathcal{A}^{+}$. By Lemma~\ref{Quotients of pro-graded Huber pairs 0}, both $A/(I\cap A)$ and $\mathcal{A}^{+}/(I\cap \mathcal{A}^{+})$ are open pro-$G$-graded subrings of $\mathcal{A}/I$ (in particular, $A/(I\cap A)$ is a pro-$G$-graded ring of definition contained in $(\mathcal{A}/I)^{+}$). Moreover, $(\mathcal{A}^{+}/(I\cap \mathcal{A}^{+}), \mathcal{A}_{\gr}^{+}/(I\cap \mathcal{A}_{\gr}^{+}))$ is a pro-$G$-graded Huber ring.

Let $B_{\gr}$ be the integral closure of $\mathcal{A}_{\gr}^{+}/(I\cap \mathcal{A}_{\gr}^{+})$ in $\mathcal{A}_{\gr}/(I\cap \mathcal{A}_{\gr})$. By Lemma~\ref{Integral closure and grading}, $B_{\gr}$ is a $G$-graded subring of $\mathcal{A}_{\gr}/(I\cap \mathcal{A}_{\gr})$. Since $B_{\gr}$ is an open integrally closed subring of the Huber ring $\mathcal{A}_{\gr}/(I\cap \mathcal{A}_{\gr})$, the topological closure $B$ of $B_{\gr}$ in $\mathcal{A}/I$ is an integrally closed subring of $\mathcal{A}/I$. By fiat, $B$ also contains $\mathcal{A}^{+}/(I\cap \mathcal{A}^{+})$. It follows that $(\mathcal{A}/I)^{+}\subseteq B$. On the other hand, we also have $B\subseteq (\mathcal{A}/I)^{+}$ since $(\mathcal{A}/I)^{+}$ is both topologically and integrally closed in $\mathcal{A}/I$ and contains $B_{\gr}$. Thus, we have exhibited $(\mathcal{A}/I)^{+}$ as the topological closure in $\mathcal{A}/I$ of a $G$-graded subring of $\mathcal{A}_{\gr}/(I\cap \mathcal{A}_{\gr})$. By Lemma~\ref{Closed pro-graded subrings}, this means that $(\mathcal{A}/I)^{+}$ is a pro-$G$-graded subring. \end{proof}
Lemma~\ref{Quotients of pro-graded Huber pairs} allows us to introduce another class of examples of pro-graded Huber pairs, coming from classical affinoid algebras over a nonarchimedean field. 
\begin{mydef}[Pro-graded affinoid algebras]\label{Pro-graded affinoid algebras, definition}Let $K$ be a nonarchimedean field, equipped with the trivial $\mathbb{Z}$-grading. A pro-$\mathbb{Z}$-graded Tate ring $(\mathcal{A}, \mathcal{A}_{\gr})$ over $(K, K)$ is called a pro-graded affinoid $K$-algebra if there exists an integer $n\geq0$ and a morphism of pro-$\mathbb{Z}$-graded Tate rings \begin{equation*}(\varphi, \varphi_{\gr}): (T_{n, K}, (T_{n, K})_{\gr})\to (\mathcal{A}, \mathcal{A}_{\gr}),\end{equation*}where $(T_{n, K}, (T_{n, K})_{\gr})$ is the pro-graded Tate algebra from Example~\ref{Pro-graded Tate algebra}, such that $\varphi$ is surjective.\end{mydef}
Recall from Example~\ref{Pro-graded Tate algebra} that $(T_{n, K})_{\gr}$ is just the polynomial ring in $n$ variables with its usual grading by total degree. Note that, by Lemma~\ref{Surjective morphisms of pro-graded rings}, the map of graded rings $(T_{n, K})_{\gr}\to \mathcal{A}_{\gr}$ in the above definition is always surjective, so $\mathcal{A}_{\gr}$ is an affine algebra over $K$ which is dense in the affinoid algebra $\mathcal{A}$. Moreover, since the polynomial ring is $\mathbb{Z}_{\geq0}$-graded, this also implies that every pro-graded affinoid algebra $(\mathcal{A}, \mathcal{A}_{\gr})$ is a pro-$\mathbb{Z}_{\geq0}$-graded Huber ring. We prove that each pro-graded affinoid algebra naturally defines a pro-graded Huber pair.
\begin{lemma}\label{Pro-graded affinoid algebras define pro-graded Huber pairs}Let $K$ be a nonarchimedean field and let $(\mathcal{A}, \mathcal{A}_{\gr})$ be a pro-graded affinoid $K$-algebra. Then $((\mathcal{A}, \mathcal{A}_{\gr}), \mathcal{A}^{\circ})$ is a pro-$\mathbb{Z}_{\geq0}$-graded Tate Huber pair.\end{lemma}
\begin{proof}We have already observed above that the pro-$\mathbb{Z}$-graded Huber ring $(\mathcal{A}, \mathcal{A}_{\gr})$ is pro-$\mathbb{Z}_{\geq0}$-graded. Let \begin{equation*}(\varphi, \varphi_{\gr}): (T_{n, K}, (T_{n, K})_{\gr})\to (\mathcal{A}, \mathcal{A}_{\gr})\end{equation*}be a morphism of pro-graded Tate rings such that $\varphi$ is surjective. By Proposition~\ref{Graded maps and homogeneous ideals 2}, $I=\ker(\varphi)$ is a topologically homogeneous (closed) ideal of $T_{n, K}$. By \cite{BGR}, Prop.~6.3.4/1, $\mathcal{A}^{\circ}$ is the integral closure of $\varphi(T_{n, K}^{\circ})$ in $\mathcal{A}$. Hence the assertion follows from Lemma~\ref{Quotients of pro-graded Huber pairs}.
\end{proof}
We also observe that for pro-graded affinoid algebras, there is no difference between the notions of homogeneous and topologically homogeneous ideals.
\begin{lemma}\label{Noetherian Banach algebras and homogeneous ideals}If $(\mathcal{A}, \mathcal{A}_{\gr})$ is a pro-$G$-graded Tate ring such that the underlying ring $\mathcal{A}$ is Noetherian, then an ideal of $\mathcal{A}$ is topologically homogeneous if and only if it is homogeneous. In particular, an ideal in a pro-graded affinoid $K$-algebra over a nonarchimedean field $K$ is topologically homogeneous if and only if it is homogeneous.\end{lemma}
\begin{proof}This follows from Lemma~\ref{Closed topologically homogeneous ideals} and the fact that all ideals in any Noetherian complete Tate ring are closed (\cite{Henkel}, Proposition~2.11).\end{proof}
We conclude this section by introducing the notion of topologically homogeneous valuations, which our adic Proj construction in Section~\ref{sec:adic Proj} will rest upon. 
\begin{mydef}[Topologically homogeneous valuation]\label{Topologically homogeneous valuations}Let $G$ be an abelian group and let $\Gamma$ be a totally ordered abelian group. \begin{enumerate}[(1)] \item For a $G$-graded ring $\mathcal{A}$, a valuation \begin{equation*}v: \mathcal{A}\to \Gamma\cup\{0\}\end{equation*}is called homogeneous if, for every $f\in \mathcal{A}$, \begin{equation*}v(f)=\max\{\, v(f_{g})\mid g\in G\,\},\end{equation*}where $f=\sum_{g\in G}f_{g}$ is the decomposition of $f$ into its homogeneous components.\item For a pro-$G$-graded Huber ring $(\mathcal{A}, \mathcal{A}_{\gr})$, a continuous valuation \begin{equation*}v: \mathcal{A}\to \Gamma\cup\{0\}\end{equation*}is called topologically homogeneous if its restriction to $\mathcal{A}_{\gr}$ is a homogeneous valuation on the $G$-graded ring $\mathcal{A}_{\gr}$.\end{enumerate}\end{mydef}
The above notion can be described via topological homogeneous decompositions of elements as in Lemma~\ref{Topological homogeneous decomposition} and Definition~\ref{Topological homogeneous decomposition, definition}
\begin{lemma}\label{Topologically homogeneous valuations and topological homogeneous decompositions}Let $G$ be an abelian group and let $\Gamma$ be a totally ordered abelian group. For a pro-$G$-graded Huber ring $(\mathcal{A}, \mathcal{A}_{\gr})$ and a continuous valuation \begin{equation*}v: \mathcal{A}\to \Gamma\cup\{0\},\end{equation*}the following are equivalent: \begin{enumerate}[(1)]\item $v$ is topologically homogeneous. \item For every element $f\in \mathcal{A}$, \begin{equation*}v(f)=\max_{g\in G}v(f_{g}),\end{equation*}where \begin{equation*}f=\sum_{g\in G}f_{g}\end{equation*}is the topological homogeneous decomposition of $f$ in $\mathcal{A}$.\end{enumerate}\end{lemma}
Note that the maximum in the statement of the above lemma always exists by virtue of the continuity of $v$ and convergence of the series $\sum_{g\in G}f_{g}$ in $\mathcal{A}$.
\begin{proof}It is clear that (2) implies (1). Suppose that the continuous valuation $v$ satisfies (1) and let $f\in \mathcal{A}$ with topological homogeneous decomposition $f=\sum_{g\in G}f_{g}$. Let $(G_{n})_{n}$ be an increasing sequence of finite subsets of $\supp(f)=\{g\in G\mid f_{g}\neq0\}$ exhausting $\supp(f)$. Then $f_{n}=\sum_{g\in G_{n}}f_{g}\in \mathcal{A}_{\gr}$ defines a sequence $(f_{n})_{n}$ in $\mathcal{A}_{\gr}$ with $f_{n}\to f$. Applying (1), we see that $v(f_{n})=\max\{v(f_{g})\mid g\in G_{n}\}$ for every $n$. By the continuity of $v$, there exists an integer $N\geq1$ such that $v(f)=v(f_{n})$ for all $n\geq N$. We conclude by the continuity of $v$ that $v(f)=\max\{v(f_{g})\mid g\in G\}$, as desired. \end{proof}
\begin{lemma}\label{Topologically homogeneous valuations and morphisms}For any morphism of pro-$G$-graded Huber pairs \begin{equation*}(\varphi, \varphi_{\gr}): ((\mathcal{A}, \mathcal{A}_{\gr}), \mathcal{A}^{+})\to ((\mathcal{B}, \mathcal{B}_{\gr}), \mathcal{B}^{+}),\end{equation*}the induced continuous map \begin{equation*}\Spa(\varphi): \Spa(\mathcal{B}, \mathcal{B}^{+})\to \Spa(\mathcal{A}, \mathcal{A}^{+})\end{equation*}sends topologically homogeneous valuations to topologically homogeneous valuations.\end{lemma}
\begin{proof}Let $f\in \mathcal{A}$ with topological homogeneous decomposition $f=\sum_{g\in G}f_{g}$. By Lemma~\ref{Homogeneous decompositions and morphisms}, $\varphi(f)=\sum_{g\in G}\varphi(f_{g})$ is the topological homogeneous decomposition of $\varphi(f)$. Using that $v$ is topologically homogeneous, we obtain \begin{equation*}\Spa(\varphi)(v)(f)=v(\varphi(f))=\max\{\, v(\varphi(f_{g}))\mid g\in G\,\}=\max\{\, \Spa(\varphi)(v)(f_{g})\mid g\in G\,\},\end{equation*}i.e., $\Spa(\varphi)(v)$ is topologically homogeneous. \end{proof}
We conclude with a general recipe for producing (topologically) homogeneous valuations. Various special cases of this statement appear to be well-known, but we could not locate a reference in the literature for the precise statement we needed.  
\begin{lemma}\label{Gauss Lemma}Let $G$ and $\Gamma$ be totally ordered abelian groups. Let $\mathcal{A}=\bigoplus_{g\in G}\mathcal{A}_{g}$ be a graded ring and let \begin{equation*}v: \bigcup_{g\in G}\mathcal{A}_{g}\to \Gamma\cup\{0\}\end{equation*}be a map from the set of homogeneous elements of $\mathcal{A}$ such that \begin{enumerate}[(1)]\item for every pair of homogeneous elements $f, f^{\prime}\in \mathcal{A}_{g}$ of the same degree $g$, we have \begin{equation*}v(f+f^{\prime})\leq \max\{v(f), v(f^{\prime})\},\end{equation*}and \item for every pair of homogeneous elements $f\in \mathcal{A}_{g}$, $f^{\prime}\in \mathcal{A}_{g^{\prime}}$ of possibly different degrees $g, g^{\prime}\in G$, \begin{equation*}v(ff^{\prime})=v(f)v(f^{\prime}).\end{equation*}\end{enumerate}For an arbitrary element $f\in \mathcal{A}$ with homogeneous decomposition $f=\sum_{g\in G}f_{g}$, set \begin{equation*}v(f)=\max\{\, v(f_{g})\mid g\in G\,\}.\end{equation*}Then $v: \mathcal{A}\to \Gamma\cup\{0\}$ is a homogeneous valuation on $\mathcal{A}$. \end{lemma}
\begin{proof}It is clear that the triangle inequality for all elements follows from the triangle inequality for homogeneous elements by looking at homogeneous decompositions. For $h, h^{\prime}\in \mathcal{A}$, with homogeneous decompositions $h=\sum_{g\in G}h_{g}$ and $h^{\prime}=\sum_{g\in G}h_{g}^{\prime}$, the product $hh^{\prime}$ has homogeneous components \begin{equation*}(hh^{\prime})_{g}=\sum_{g_{1}+g_{2}=g}h_{g_{1}}h_{g_{2}}, g\in G.\end{equation*}Hence, for every $g\in G$, \begin{equation*}v((hh^{\prime})_{g})\leq \max_{g_{1}+g_{2}=g}v(h_{g_{1}}h_{g_{2}}^{\prime})=\max_{g_{1}+g_{2}=g}v(h_{g_{1}})v(h_{g_{2}}^{\prime})\leq v(h)v(h^{\prime}).\end{equation*}To prove the converse inequality, let \begin{equation*}g=\min\{g\in G\mid v(h_{g})=v(h)\}\end{equation*}and \begin{equation*}g^{\prime}=\min\{g\in G\mid v(h_{g}^{\prime})=v(h^{\prime})\}.\end{equation*}For any $g_1, g_2\in G$ with $g_{1}+g_{2}=g+g^{\prime}$ and $(g_1, g_2)\neq (g, g^{\prime})$, we must have either $g_{1}<g$ or $g_{2}<g^{\prime}$. By the definition of $g$ and $g^{\prime}$, this implies $v(h_{g_{1}})<v(h)$ or $v(h_{g_{2}}^{\prime})<v(h^{\prime})$, whence \begin{equation*}v(h_{g_{1}}h_{g_{2}}^{\prime})<v(h)v(h^{\prime}).\end{equation*}It follows that \begin{equation*}v((hh^{\prime})_{g+g^{\prime}})=v(h_{g}h_{g^{\prime}}^{\prime})=v(h_{g})v(h_{g^{\prime}}^{\prime})=v(h)v(h^{\prime}).\end{equation*}In particular, $v(h)v(h^{\prime})\leq v(hh^{\prime})$, as desired. \end{proof}
\begin{rmk}\label{Remark on the Gauss Lemma}Note that the above lemma can be seen as a generalization of the Gauss Lemma. Indeed, a (rank-$1$) valuation $v$ on a field $K$ necessarily gives a map on homogeneous elements of the graded ring $K[X]$ satisfying the assumptions of Lemma~\ref{Gauss Lemma} by setting $v(aX^{i})=v(a)$ for all $a\in K$, $i\geq0$. The conclusion of Lemma~\ref{Gauss Lemma} then specializes to the statement that the Gauss norm on $K[X]$ is a valuation, i.e., the classical Gauss Lemma.\end{rmk}   

\section{The adic $\Proj$ construction}\label{sec:adic Proj}

We again fix an abelian group $G$, which now almost always will be assumed totally ordered (and, in particular, torsion-free). Most of the time we will restrict our attention to the two cases $G=\mathbb{Z}$ and $G=\mathbb{Z}[1/p]$. 

In this section, we finally introduce our adic space analog of the $\Proj$ construction from algebraic geometry. Before we proceed, let us recall the classical construction in a form which makes the motivation for our construction in adic geometry the most apparent. Recall that, in algebraic geometry, one attaches to a positively graded ring \begin{equation*}A=\bigoplus_{d\in\mathbb{Z}_{\geq0}}A_{d}\end{equation*}a topological space $\Proj(A)$ which is defined as the subspace of $\Spec(A)$ consisting of homogeneous prime ideals of $A$ which do not contain the irrelevant prime ideal \begin{equation*}A_{+}=\bigoplus_{d>0}A_{d}\end{equation*}of $A$. One then equips $\Proj(A)$ with the structure of a scheme as follows. Choosing a family of homogeneous elements of positive degree $(f_{s})_{s\in S}$ of $A$ with \begin{equation*}\sqrt{(f_{s}\mid s\in S)_{A}}=A_{+},\end{equation*}one considers the open subsets \begin{equation*}D_{+}(f_{s})=\{\, \mathfrak{p}\in \Proj(A)\mid f_{s}\not\in\mathfrak{p}\,\}\end{equation*}and notices that each of these is homeomorphic to the subspace of $\Spec(A[f_{s}^{-1}])$ consisting of the prime ideals of the localization $A[f_{s}^{-1}]$ which are homogeneous with respect to the natural $\mathbb{Z}$-grading on $A[f_{s}^{-1}]$, the homeomorphism being the map $\mathfrak{p}\mapsto \mathfrak{p}A[f_{s}^{-1}]$. We can denote by $\Proj(A[f^{-1}])$ the set of homogeneous prime ideals of the $\mathbb{Z}$-graded ring $A[f^{-1}]$ for a homogeneous element of positive degree $f\in A$ (note that this is not a special case of the $\Proj$ defined above since $A[f^{-1}]$ is not positively graded). One then observes that, for any homogeneous element $f$ in a graded ring $A=\bigoplus_{d}A_{d}$, the restriction map \begin{equation*}\Proj(A[f^{-1}])\to \Spec(A[f^{-1}]_{0})\end{equation*}is a homeomorphism, and one defines the scheme structure on $\Proj(A)$ by gluing the affine open subsets \begin{equation*}\Proj(A[f_{s}^{-1}])\cong \Spec(A[f_{s}^{-1}]_{0})\end{equation*}along their intersections. Our construction of the projective adic spectrum of a pro-graded Huber pair will be entirely analogous to this procedure, with the place of homogeneous prime ideals now occupied by topologically homogeneous valuations on a pro-graded Huber ring.

We begin by introducing a notion of irrelevant ideal for pro-$G_{\geq0}$-graded Huber rings.  
\begin{mydef}[Irrelevant ideal]\label{Irrelevant ideal}Let $G$ be a totally ordered abelian group and let $(\mathcal{A}, \mathcal{A}_{\gr})$ be a pro-$G_{\geq0}$-graded Huber ring. The irrelevant ideal of $\mathcal{A}$ is the (topological) closure inside $\mathcal{A}$ of the irrelevant ideal \begin{equation*}\bigoplus_{g>0}\mathcal{A}_{\gr, g}\end{equation*}of the $G_{\geq0}$-graded ring $\mathcal{A}_{\gr}$. We denote the irrelevant ideal by $(\mathcal{A}, \mathcal{A}_{\gr})_{+}$. \end{mydef}
In particular, we obtain a notion of irrelevant ideal $(R, R_{\gr})_{+}$ for any pro-$G$-graded adic ring $(R, R_{\gr})$ when $G$ is a totally ordered abelian group. 

Recall that for every uniform Huber pair we have the continuous specialization map \begin{equation*}\spc: \Spa(\mathcal{A}, \mathcal{A}^{+})\to \Spf(\mathcal{A}^{+})\end{equation*}which sends a continuous valuation $v\in \Spa(\mathcal{A}, \mathcal{A}^{+})$ to its center \begin{equation*}\spc(v)=\{\, f\in \mathcal{A}^{+}\mid v(f)<1\,\}\end{equation*}on $\mathcal{A}^{+}$. If $\mathcal{A}$ is a Tate ring, this map is closed and surjective (this follows, for example, from \cite{BhattNotes}, Theorem~8.1.2). By abuse of notation, we will denote the center on $\mathcal{A}^{+}$ of an element $v\in \Spa(\mathcal{A}, \mathcal{A}^{+})$ of the adic spectrum by $\spc(v)$ even if the Huber pair $(\mathcal{A}, \mathcal{A}^{+})$ is not necessarily uniform. We can now define the $\Proj$ of a pro-graded Huber pair. 
\begin{mydef}[The projective adic spectrum, as a topological space]\label{Adic Proj, part 1}Let $G$ be a totally ordered abelian group and let $((\mathcal{A}, \mathcal{A}_{\gr}), \mathcal{A}^{+})$ be a pro-$G_{\geq0}$-graded Huber pair. The \textit{projective adic spectrum} of $((\mathcal{A}, \mathcal{A}_{\gr}), \mathcal{A}^{+})$, as a topological space, is the subset \begin{equation*}\Proj(\mathcal{A}, \mathcal{A}_{\gr}, \mathcal{A}^{+})\end{equation*}of $\Spa(\mathcal{A}, \mathcal{A}^{+})$ consisting of topologically homogeneous continuous valuations $v\in \Spa(\mathcal{A}, \mathcal{A}^{+})$ whose center $\spc(v)$ on $\mathcal{A}^{+}$ does not contain $(\mathcal{A}, \mathcal{A}_{\gr})_{+}\cap \mathcal{A}^{+}$, with the topology given by unions of subsets of the form \begin{equation*}\Proj(\mathcal{A}, \mathcal{A}_{\gr}, \mathcal{A}^{+})(\frac{f_1,\dots, f_n}{f_{n+1}})=\{\, v\in \Proj(\mathcal{A}, \mathcal{A}_{\gr}, \mathcal{A}^{+})\mid v(f_{i})\leq v(f_{n+1})\neq 0, i=1,\dots, n\,\},\end{equation*}where $f_1,\dots, f_{n}, f_{n+1}$ are elements of $\mathcal{A}$ generating an open ideal of $\mathcal{A}$ and $f_{n+1}$ is a homogeneous element of $\mathcal{A}_{\gr}$. 

We call subsets of $\Proj(\mathcal{A}, \mathcal{A}_{\gr}, \mathcal{A}^{+})$ of the above form \textit{homogeneous rational subsets}. \end{mydef}
Note that every homogeneous rational subset is the restriction to $\Proj(\mathcal{A}, \mathcal{A}_{\gr}, \mathcal{A}^{+})$ of a rational subset of $\Spa(\mathcal{A}, \mathcal{A}^{+})$, but, a priori, we do not know whether the restriction of every rational subset is a homogeneous rational subset. Thus the topology on $\Proj(\mathcal{A}, \mathcal{A}_{\gr}, \mathcal{A}^{+})$ is a priori coarser than the subspace topology induced from $\Spa(\mathcal{A}, \mathcal{A}^{+})$. The topology on $\Proj$ can also be described using a more restricted family of homogeneous rational subsets. 
\begin{lemma}\label{Topology on Proj}The topology on $\Proj(\mathcal{A}, \mathcal{A}_{\gr}, \mathcal{A}^{+})$ is equivalently given by unions of homogeneous rational subsets \begin{equation*}\Proj(\mathcal{A}, \mathcal{A}_{\gr}, \mathcal{A}^{+})(\frac{f_1,\dots, f_n}{f_{n+1}})\end{equation*}where all of the elements $f_{1}, \dots, f_{n}, f_{n+1}$ are homogeneous. \end{lemma}
\begin{proof}It suffices to prove that every homogeneous rational subset is open with respect to the topology whose open sets are unions of subsets as in the statement of the lemma. Let $f_1,\dots, f_n, f_{n+1}$ be any family of elements of $\mathcal{A}$ generating an open ideal of $\mathcal{A}$ such that $f_{n+1}$ is a homogeneous element of $\mathcal{A}_{\gr}$ of some degree $g\in G$. Using \cite{Huber1}, Lemma~3.10, we can approximate $f_1,\dots, f_n$ by elements in the dense subring $\mathcal{A}_{\gr}$ without changing the rational subset, so we may assume that $f_1,\dots, f_n$ are elements of $\mathcal{A}_{\gr}$, with homogeneous decompositions $f_{i}=\sum_{g\in G}f_{i, g}$, $i=1,\dots, n$. In this case, \begin{equation*}\Proj(\mathcal{A}, \mathcal{A}_{\gr}, \mathcal{A}^{+})(\frac{f_{1},\dots, f_{n}}{f_{n+1}})=\bigcap_{g\in \bigcup_{i=1}^{n}\supp(f_{i})}\Proj(\mathcal{A}, \mathcal{A}_{\gr}, \mathcal{A}^{+})(\frac{f_{1, g},\dots, f_{n, g}}{f_{n+1}}).\end{equation*}Indeed, since every valuation $v$ in the adic projective spectrum is topologically homogeneous, the condition $v(f_{i})\leq v(f_{n+1})\neq0$ is equivalent to the condition $v(f_{i, g})\leq v(f_{n+1})$ for all $g\in G$ (or, equivalently, for all $g$ in the support of $f_{i}$). This finishes the proof since a finite intersection of open subsets is open in any topological space. \end{proof}
In the cases $G=\mathbb{Z}$ and $G=\mathbb{Z}[1/p]$, we can endow the adic projective spectrum of a ``sufficiently sheafy'' pro-graded Huber pair with the structure of an adic space. This is done by a gluing procedure, similar to that for the usual $\Proj$ construction in algebraic geometry.
\begin{mydef}\label{Fake adic Proj}Let $G$ be a totally ordered abelian group and let $((\mathcal{A}, \mathcal{A}_{\gr}), \mathcal{A}^{+})$ be a pro-$G_{\geq0}$-graded Huber pair. For a homogeneous element $f\in \mathcal{A}_{\gr}^{+}$, define the topological space \begin{equation*}\Proj(\mathcal{A}\langle f^{-1}\rangle, \mathcal{A}\langle f^{-1}\rangle_{\gr}, \mathcal{A}\langle f^{-1}\rangle^{+})\end{equation*}to be the subspace of $\Spa(\mathcal{A}\langle f^{-1}\rangle, \mathcal{A}\langle f^{-1}\rangle^{+})$ consisting of those $v\in \Spa(\mathcal{A}\langle f^{-1}\rangle, \mathcal{A}\langle f^{-1}\rangle^{+})$ which are topologically homogeneous on $(\mathcal{A}\langle f^{-1}\rangle, \mathcal{A}\langle f^{-1}\rangle_{\gr})$.\end{mydef}
Note that Definition~\ref{Fake adic Proj} is not a special case of Definition~\ref{Adic Proj, part 1}, since the pro-$G$-graded rational localization $(\mathcal{A}\langle f^{-1}\rangle, \mathcal{A}\langle f^{-1}\rangle_{\gr})$ is not a pro-$G_{\geq0}$-graded Huber ring unless $f$ is of degree $0$.  
\begin{lemma}\label{Topological space}Let $((\mathcal{A}, \mathcal{A}_{\gr}), \mathcal{A}^{+})$ be a pro-$G_{\geq0}$-graded Huber pair, for some totally ordered abelian group $G$, and let $(f_{s})_{s\in S}$ be a family of homogeneous elements of $\mathcal{A}^{+}$ such that $(f_{s})_{s\in S}$ generate the ideal $(\mathcal{A}, \mathcal{A}_{\gr})\cap \mathcal{A}^{+}$ of $\mathcal{A}^{+}$ up to radical and (topological) closure. Then the topological space $\Proj(\mathcal{A}, \mathcal{A}_{\gr}, \mathcal{A}^{+})$ is covered by the homogeneous rational subsets $\Proj(\mathcal{A}, \mathcal{A}_{\gr}, \mathcal{A}^{+})(\frac{1}{f_{s}})$. \end{lemma}
\begin{proof}Follows from the definitions.\end{proof}
\begin{lemma}\label{Gluing -1}Suppose that $G$ is a totally ordered abelian group. Let $((\mathcal{A}, \mathcal{A}_{\gr}), \mathcal{A}^{+})$ be a pro-$G_{\geq0}$-graded Huber pair, let $f\in \mathcal{A}_{\gr}^{+}=\mathcal{A}_{\gr}\cap \mathcal{A}^{+}$ be a homogeneous element of positive degree and let \begin{equation*}(\varphi, \varphi_{\gr}): ((\mathcal{A}, \mathcal{A}_{\gr}), \mathcal{A}^{+})\to ((\mathcal{A}\langle f^{-1}\rangle, \mathcal{A}\langle f^{-1}\rangle_{\gr}), \mathcal{A}\langle f^{-1}\rangle^{+})\end{equation*}be the canonical morphism of pro-$G$-graded Huber pairs. Then the homeomorphism \begin{equation*}\Spa(\varphi): \Spa(\mathcal{A}\langle f^{-1}\rangle, \mathcal{A}\langle f^{-1}\rangle^{+})\tilde{\to} \Spa(\mathcal{A}, \mathcal{A}^{+})(\frac{1}{f})\end{equation*}induced by $\varphi$ restricts to a homeomorphism \begin{equation*}\Proj(\mathcal{A}\langle f^{-1}\rangle, \mathcal{A}\langle f^{-1}\rangle_{\gr}, \mathcal{A}\langle f^{-1}\rangle^{+})\tilde{\to}\Proj(\mathcal{A}, \mathcal{A}_{\gr}, \mathcal{A}^{+})(\frac{1}{f}).\end{equation*}\end{lemma}
\begin{proof}By Lemma~\ref{Topologically homogeneous valuations and morphisms}, $\Spa(\varphi)(v)$ is topologically homogeneous for every \begin{equation*}v\in \Proj(\mathcal{A}\langle f^{-1}\rangle, \mathcal{A}\langle f^{-1}\rangle_{\gr}, \mathcal{A}\langle f^{-1}\rangle^{+}).\end{equation*}Moreover, no element of $\Spa(\mathcal{A}\langle f^{-1}\rangle, \mathcal{A}\langle f^{-1}\rangle^{+})$ or of $\Proj(\mathcal{A}, \mathcal{A}_{\gr}, \mathcal{A}^{+})(\frac{1}{f})$ can contain all homogeneous elements of $\mathcal{A}^{+}$ of positive degree in its center since $f$ is a homogeneous element of positive degree. This shows that $\Spa(\varphi)$ restricts to a map \begin{equation*}\Proj(\mathcal{A}\langle f^{-1}\rangle, \mathcal{A}\langle f^{-1}\rangle_{\gr}, \mathcal{A}\langle f^{-1}\rangle^{+})\to \Proj(\mathcal{A}, \mathcal{A}_{\gr}, \mathcal{A}^{+})(\frac{1}{f}).\end{equation*}Since $\Spa(\varphi)$ is a homeomorphism, its restriction is also a homeomorphism onto its image.

It remains to show that if $v\in \Spa(\mathcal{A}\langle f^{-1}\rangle, \mathcal{A}\langle f^{-1}\rangle^{+})$ and the restriction $\Spa(\varphi)(v)$ of $v$ is topologically homogeneous, then $v$ must itself be topologically homogeneous. Let $h\in \mathcal{A}\langle f^{-1}\rangle_{\gr}$ with topological homogeneous decomposition $h=\sum_{g\in G}h_{g}$. We want to show that $v(h)=\max\{\, v(h_{g})\mid g\in G\,\}$. After approximating $h$ with elements of the dense subring $\mathcal{A}_{\gr}[\frac{1}{f}]$ and using Lemma~\ref{Grading and topology}, we may assume that $h\in \mathcal{A}_{\gr}[\frac{1}{f}]$. Choose $k\geq0$ such that $f^{k}h\in \mathcal{A}_{\gr}$. Since $f\in \mathcal{A}_{\gr}$ is homogeneous, $f^{k}h=\sum_{g\in G}f^{k}h_{g}$ is the homogeneous decomposition of $f^{k}h$ in $\mathcal{A}_{\gr}$. Since $v$ restricts to the topologically homogeneous valuation $\Spa(\varphi)(v)$ on $\mathcal{A}$ and since $v(f)=1$, we get \begin{equation*}v(h)=v(f^{k}h)=\max\{\, v(f^{k}h_{g})\mid g\in G\,\}=\max\{\, v(h_{g})\mid g\in G\,\},\end{equation*}showing that $v$ is topologically homogeneous, as desired. \end{proof}
\begin{lemma}\label{Gluing 0}Suppose that either $G=\mathbb{Z}$ or $G=\mathbb{Z}[1/p]$. Let $((\mathcal{A}, \mathcal{A}_{\gr}), \mathcal{A}^{+})$ be a pro-$G_{\geq0}$-graded Huber pair and let $f\in \mathcal{A}_{\gr}^{+}=\mathcal{A}_{\gr}\cap \mathcal{A}^{+}$ be a homogeneous element of degree~$1$. If $G=\mathbb{Z}[1/p]$, assume that $f$ admits a compatible system of $p$-power roots in $\mathcal{A}_{\gr}^{+}$. Then the composition \begin{equation*}\Proj(\mathcal{A}, \mathcal{A}_{\gr}, \mathcal{A}^{+})(\frac{1}{f})\tilde{\to} \Proj(\mathcal{A}\langle f^{-1}\rangle, \mathcal{A}\langle f^{-1}\rangle_{\gr}, \mathcal{A}\langle f^{-1}\rangle^{+}) \to \Spa(\mathcal{A}\langle f^{-1}\rangle_{\gr, 0}, \mathcal{A}\langle f^{-1}\rangle_{\gr, 0}^{+}),\end{equation*}where the first map is the homeomorphism from Lemma~\ref{Gluing -1}, where $\mathcal{A}\langle f^{-1}\rangle_{\gr, 0}^{+}$ is the integral closure of $\mathcal{A}^{+}\langle f^{-1}\rangle_{\gr, 0}$ in $\mathcal{A}\langle f^{-1}\rangle_{\gr, 0}$ and where the second map is the natural restriction map, is a homeomorphism.\end{lemma}
\begin{proof}The restriction map is clearly continuous, so it suffices to prove that it is bijective with a continuous inverse. We construct the inverse map $\psi$. Let $v\in \Spa(\mathcal{A}\langle f^{-1}\rangle_{\gr, 0}, \mathcal{A}\langle f^{-1}\rangle_{\gr, 0}^{+})$, regarded as a continuous valuation \begin{equation*}v: \mathcal{A}\langle f^{-1}\rangle_{\gr,0}\to \Gamma\cup\{0\},\end{equation*}for some totally ordered abelian group $\Gamma$. If $G=\mathbb{Z}[1/p]$, we fix a compatible system of $p$-power roots $(f^{1/p^{n}})_{n\geq0} $ of $f$ in $\mathcal{A}_{\gr}^{+}$ and, for any $d=\frac{k}{p^{n}}\in \mathbb{Z}[1/p]_{\geq0}$ (with $k, n\in \mathbb{Z}_{\geq0}$), set $f^{d}=(f^{1/p^{n}})^{k}$ (this does not depend on the representation of $d$ as a fraction since the $p$-power roots are compatible with each other). To define the continuous, topologically homogeneous valuation $\psi(v)$ on $\mathcal{A}\langle f^{-1}\rangle$ (with $\psi(v)(h)\leq1$ for all $h\in \mathcal{A}\langle f^{-1}\rangle^{+}$), it suffices to define it on the dense subring $\mathcal{A}_{\gr}[\frac{1}{f}]$. 

Let $h\in \mathcal{A}_{\gr}[\frac{1}{f}]$ be homogeneous of degree $d\in G$. Then \begin{equation*}\frac{h}{f^{d}}\in \mathcal{A}_{\gr}[\frac{1}{f}]\subseteq \mathcal{A}\langle f^{-1}\rangle_{\gr}\end{equation*}is homogeneous of degree $0$. Set \begin{equation*}\psi(v)(h)=v(\frac{h}{f^{d}}).\end{equation*}This defines a map $\psi(v)$ from the set of homogeneous elements of $\mathcal{A}_{\gr}[\frac{1}{f}]$ to $\Gamma\cup\{0\}$. For any two homogeneous elements $h, h^{\prime}$ of the same degree $d$, we have \begin{equation*}\psi(v)(h+h^{\prime})=v(\frac{h+h^{\prime}}{f^{d}})\leq \max\{v(\frac{h}{f^{d}}), v(\frac{h^{\prime}}{f^{d}})\}=\max\{\psi(v)(h), \psi(v)(h^{\prime})\}.\end{equation*}For any two homogeneous elements $h, h^{\prime}\in \mathcal{A}_{\gr}[\frac{1}{f}]$ of possibly different degrees $d_{1}, d_{2}$, we have \begin{equation*}\psi(v)(hh^{\prime})=v(\frac{hh^{\prime}}{f^{d_{1}+d_{2}}})=v(\frac{h}{f^{d_{1}}})v(\frac{h^{\prime}}{f^{d_{2}}})=\psi(v)(h)\psi(v)(h^{\prime}).\end{equation*}By Lemma~\ref{Gauss Lemma}, this entails that $\psi(v)$ extends to a homogeneous valuation on the $G$-graded ring $\mathcal{A}_{\gr}[\frac{1}{f}]$ by setting \begin{equation*}\psi(v)(h)=\max\{\psi(v)(h_{g})\mid g\in G\}.\end{equation*}If $(h_{n})_{n}$ is a sequence in $\mathcal{A}_{\gr}[\frac{1}{f}]$ converging to some $h\in \mathcal{A}_{\gr}[\frac{1}{f}]$, then by Lemma~\ref{Grading and topology}, $h_{n, g}\to h_{g}$ for all $g\in G$ and then \begin{equation*}\psi(v)(h_{n, g}-h_{g})=v(\frac{h_{n, g}-h_{g}}{f^{g}})\to 0\end{equation*}by the continuity of $v$. Hence the homogeneous valuation $\psi(v)$ on $\mathcal{A}_{\gr}[\frac{1}{f}]$ is continuous, and we extend it by continuity to a topologically homogeneous continuous valuation on the completion $\mathcal{A}\langle f^{-1}\rangle$, which we again denote by $\psi(v)$. 

To prove that $\psi(v)(h)\leq1$ for all $h\in \mathcal{A}\langle f^{-1}\rangle^{+}$, it suffices to verify that $\psi(v)(h)\leq 1$ for all $h\in \mathcal{A}_{\gr}^{+}[\frac{1}{f}]$. Let $h\in \mathcal{A}_{\gr}^{+}[\frac{1}{f}]$. By Lemma~\ref{Lemma on rational localizations} and the assumption that $\mathcal{A}_{\gr}^{+}$ is a $G$-graded subring of $\mathcal{A}_{\gr}$, we know that the homogeneous components $h_{g}$ of $h$ also belong to $\mathcal{A}_{\gr}^{+}[\frac{1}{f}]$. Then, for every $g\in G$, \begin{equation*}\frac{h_{g}}{f^{g}}\in \mathcal{A}_{\gr}^{+}[\frac{1}{f}]_{0}\subseteq \mathcal{A}^{+}\langle f^{-1}\rangle_{\gr, 0}\subseteq \mathcal{A}\langle f^{-1}\rangle_{\gr, 0}^{+}\end{equation*}and thus \begin{equation*}\psi(v)(h_{g})=v(\frac{h_{g}}{f^{g}})\leq 1.\end{equation*}Since this holds for all $g\in G$, it follows that $\psi(v)(h)\leq1$. 

To prove that the map $\psi$ is continuous, by Lemma~\ref{Topology on Proj}, it suffices to prove that the pre-image of (the intersection with $\Proj(\mathcal{A}, \mathcal{A}_{\gr}, \mathcal{A}^{+})(\frac{1}{f})$ of) a homogeneous rational subset \begin{equation*}U=\Proj(\mathcal{A}, \mathcal{A}_{\gr}, \mathcal{A}^{+})(\frac{f_1,\dots, f_n}{f_{n+1}})\end{equation*}of $\Proj(\mathcal{A}, \mathcal{A}_{\gr}, \mathcal{A}^{+})$, where all of $f_1, \dots, f_{n}, f_{n+1}$ are homogeneous elements of $\mathcal{A}_{\gr}$, is an open subset of $\Spa(\mathcal{A}\langle f^{-1}\rangle_{\gr,0}, \mathcal{A}\langle f^{-1}\rangle_{\gr,0}^{+})$. To this end, let $f_1,\dots, f_n, f_{n+1}\in \mathcal{A}_{\gr}$ be homogeneous elements generating an open ideal in $\mathcal{A}$. Then, since all $f_{i}$ are homogeneous, the pre-image under $\psi$ of the open subset \begin{equation*}U\cap \Proj(\mathcal{A}\langle f^{-1}\rangle, \mathcal{A}\langle f^{-1}\rangle_{\gr}, \mathcal{A}\langle f^{-1}\rangle^{+})\end{equation*}is the open subset \begin{equation*}\{\, v\in \Spa(\mathcal{A}\langle f^{-1}\rangle_{\gr,0}, \mathcal{A}\langle f^{-1}\rangle_{\gr,0}^{+})\mid v(\frac{f_{i}}{f^{\deg(f_{i})}})\leq v(\frac{f_{n+1}}{f^{\deg(f_{n+1})}})\neq0, i=1,\dots, n+1\,\}.\end{equation*}This shows that the map $\psi$ is continuous.  

Finally, we have to show that the map \begin{equation*}\psi: \Spa(\mathcal{A}\langle f^{-1}\rangle_{\gr,0}, \mathcal{A}\langle f^{-1}\rangle_{\gr,0}^{+})\to \Proj(\mathcal{A}\langle f^{-1}\rangle, \mathcal{A}\langle f^{-1}\rangle_{\gr}, \mathcal{A}\langle f^{-1}\rangle^{+})\end{equation*}is inverse to the restriction map $\varphi$. Let $v\in \Spa(\mathcal{A}\langle f^{-1}\rangle_{\gr,0}, \mathcal{A}\langle f^{-1}\rangle_{\gr,0}^{+})$. To prove that $\varphi(\psi(v))=v$, it suffices to verify that the two continuous valuations are equal on all elements of the dense subring $\mathcal{A}_{\gr}[\frac{1}{f}]_{0}$ of $\mathcal{A}\langle f^{-1}\rangle_{\gr,0}$. Every element of $\mathcal{A}_{\gr}[\frac{1}{f}]_{0}$ can be written as $\frac{h}{f^{\deg(h)}}$ for some homogeneous element $h\in \mathcal{A}_{\gr}$. But \begin{equation*}\varphi(\psi(v))(\frac{h}{f^{\deg(h)}})=v(\frac{h}{f^{\deg(h)}})\end{equation*}by the definition of $\psi(v)$. Similarly, for a valuation $w\in \Proj(\mathcal{A}\langle f^{-1}\rangle, \mathcal{A}\langle f^{-1}\rangle_{\gr}, \mathcal{A}\langle f^{-1}\rangle^{+})$, it suffices to prove that $\psi(\varphi(w))$ coincides with $w$ on the dense subring $\mathcal{A}_{\gr}[\frac{1}{f}]$. Moreover, since $w$ and $\psi(\varphi(w))$ are topologically homogeneous, it suffices to show that they coincide on homogeneous elements of this $G$-graded ring. For $h\in \mathcal{A}_{\gr}[\frac{1}{f}]_{g}$, we verify: \begin{equation*}\psi(\varphi(w))(h)=\varphi(w)(\frac{h}{f^{g}}).\end{equation*}But $w(f)=1$, so $\varphi(w)(\frac{h}{f^{g}})=w(h)$. \end{proof}
\begin{cor}\label{Topological space 2}Let $G=\mathbb{Z}$ or $G=\mathbb{Z}[1/p]$ and let $((\mathcal{A}, \mathcal{A}_{\gr}), \mathcal{A}^{+})$ be a pro-$G_{\geq0}$-graded Huber pair. Suppose that there exist homogeneous elements of degree~$1$ \begin{equation*}f_1,\dots, f_{n}\in \mathcal{A}_{\gr}^{+}\end{equation*}which generate the irrelevant ideal of $(\mathcal{A}, \mathcal{A}_{\gr})$ up to topological closure and radical and if $G=\mathbb{Z}[1/p]$, suppose that $f_1,\dots, f_n$ admit compatible systems of $p$-power roots in $\mathcal{A}_{\gr}^{+}$. Then $\Proj(\mathcal{A}, \mathcal{A}_{\gr}, \mathcal{A}^{+})$ is quasi-compact.\end{cor}
\begin{proof}By Lemma~\ref{Topological space}, $\Proj(\mathcal{A}, \mathcal{A}_{\gr}, \mathcal{A}^{+})$ is covered by the open subspaces $\Proj(\mathcal{A}, \mathcal{A}_{\gr}, \mathcal{A}^{+})(\frac{1}{f_{i}})$, $i=1,\dots, n$. By Lemma~\ref{Gluing 0}, each of these open subspaces is homeomorphic to the adic spectrum of a Huber pair, hence quasi-compact. \end{proof}
\begin{lemma}\label{Gluing 1}Set assumptions as in Lemma~\ref{Gluing 0}. Let $f_1, f_2\in \mathcal{A}_{\gr}^{+}=\mathcal{A}_{\gr}\cap \mathcal{A}^{+}$ be homogeneous elements of degree~$1$ and, in case $G=\mathbb{Z}[1/p]$, suppose that $f_{1}$, $f_{2}$ admit compatible systems of $p$-power roots in $\mathcal{A}_{\gr}^{+}$. Then each of the maps \begin{equation*}\Proj(\mathcal{A}, \mathcal{A}_{\gr}, \mathcal{A}^{+})(\frac{1}{f_i})\to \Spa(\mathcal{A}\langle f_{i}^{-1}\rangle_{\gr, 0}, \mathcal{A}\langle f_{i}^{-1}\rangle_{\gr, 0}^{+}), i=1, 2,\end{equation*}from Lemma~\ref{Gluing 0} restricts to the same map \begin{equation*}\Proj(\mathcal{A}, \mathcal{A}_{\gr}, \mathcal{A}^{+})(\frac{1}{f_1})\cap \Proj(\mathcal{A}, \mathcal{A}_{\gr}, \mathcal{A}^{+})(\frac{1}{f_2})\to \Spa(\mathcal{A}\langle (f_{1}f_{2})^{-1}\rangle_{\gr, 0}, \mathcal{A}\langle (f_{1}f_{2})^{-1}\rangle_{\gr, 0}^{+})\end{equation*}and this map is a homeomorphism.\end{lemma}
\begin{proof}We clearly have \begin{equation*}\Proj(\mathcal{A}, \mathcal{A}_{\gr}, \mathcal{A}^{+})(\frac{1}{f_{1}})\cap \Proj(\mathcal{A}, \mathcal{A}_{\gr}, \mathcal{A}^{+})(\frac{1}{f_2})=\Proj(\mathcal{A}, \mathcal{A}_{\gr}, \mathcal{A}^{+})(\frac{1}{f_{1}f_{2}}).\end{equation*}Hence the assertion immediately follows from Lemma \ref{Gluing 0}.\end{proof}
\begin{lemma}\label{Gluing 1, part 2}Let $\mathcal{A}$ be any Huber ring and let $f_1, f_2\in \mathcal{A}^{\circ}$. Let $d_1, d_2\in \mathbb{Z}_{\geq0}$. Let \begin{equation*}f=\frac{f_{2}^{d_{1}}}{f_{1}^{d_{2}}}\in \mathcal{A}[f_{1}^{-1}].\end{equation*}Then \begin{equation*}\mathcal{A}[f_{1}^{-1}][f^{-1}]=\mathcal{A}[(f_{1}f_{2})^{-1}]\end{equation*}as Huber rings and, in particular,
\begin{equation*}\mathcal{A}\langle f_{1}^{-1}\rangle\langle f^{-1}\rangle=\mathcal{A}\langle (f_{1}f_{2})^{-1}\rangle.\end{equation*}\end{lemma}
\begin{proof}The second equality of Huber rings follows from the first by taking completions on both sides. It is clear that $f$ is invertible in $\mathcal{A}[(f_{1}f_{2})^{-1}]$ and that $f_{2}$ is invertible in $\mathcal{A}[f_{1}^{-1}][f^{-1}]$, with inverse $f^{-1}f_{1}^{-d_{2}}f_{2}^{d_{1}-1}$, so the first equality holds for the underlying abstract rings. By \cite{Huber1}, (1.2), to prove equality of Huber rings, it remains to show that $f^{-1}$ is power-bounded in $\mathcal{A}[f_{1}^{-1}][f_{2}^{-1}]$ and $f_{2}^{-1}$ is power-bounded in $\mathcal{A}[f_{1}^{-1}][f^{-1}]$. If $A$ is a ring of definition of $\mathcal{A}$, then $A[(f_{1}f_{2})^{-1}]$ is a ring of definition of $\mathcal{A}[(f_{1}f_{2})^{-1}]$ and $A[f_{1}^{-1}][f^{-1}]$ is a ring of definition of $\mathcal{A}[f_{1}^{-1}][f^{-1}]$, so, in particular, \begin{equation*}f_{2}^{-1}\in \mathcal{A}[(f_{1}f_{2})^{-1}]^{\circ}\end{equation*}and \begin{equation*}f_{1}^{-1}, f^{-1}\in \mathcal{A}[f_{1}^{-1}][f^{-1}]^{\circ}\end{equation*}Since we also assumed that $f_{1}, f_{2}\in \mathcal{A}^{\circ}$, we conclude that \begin{equation*}f^{-1}=(f_{2}^{-1})^{d_{1}}f_{1}^{d_{2}}\in \mathcal{A}[(f_{1}f_{2})^{-1}]^{\circ}\end{equation*}and \begin{equation*}f_{2}^{-1}=f^{-1}f_{1}^{-d_{2}}f_{2}^{d_{1}-1}\in \mathcal{A}[f_{1}^{-1}][f^{-1}]^{\circ},\end{equation*}as desired.\end{proof}
\begin{lemma}\label{Gluing 1, part 3}Let $G=\mathbb{Z}$ or $G=\mathbb{Z}[1/p]$. Let $(\mathcal{A}, \mathcal{A}_{\gr})$ be any pro-$G$-graded Huber ring. Let $f_{1}, f_{2}\in \mathcal{A}^{\circ}$ be two power-bounded elements which are homogeneous of positive degree and let $d_{1}, d_{2}\in \mathbb{Z}_{>0}$ such that $d_{2}\deg(f_{1})=d_{1}\deg(f_{2})$. Let \begin{equation*}f=\frac{f_{2}^{d_{1}}}{f_{1}^{d_{2}}}\in \mathcal{A}[f_{1}^{-1}],\end{equation*}a homogeneous element of degree $0$ in $\mathcal{A}[f_{1}^{-1}] $. Then \begin{equation*}\mathcal{A}[(f_{1}f_{2})^{-1}]_{\gr,0}=\mathcal{A}[f_{1}^{-1}]_{\gr,0}[f^{-1}]\end{equation*}as Huber rings. In particular, \begin{equation*}\mathcal{A}\langle (f_{1}f_{2})^{-1}\rangle_{\gr,0}=\mathcal{A}\langle f_{1}^{-1}\rangle_{\gr,0}\langle f^{-1}\rangle. \end{equation*}\end{lemma}
\begin{proof}The second equality follows from the first by taking completions. By Lemma~\ref{Gluing 1, part 2}, \begin{equation*}\mathcal{A}[f_{1}^{-1}][f^{-1}]=\mathcal{A}[(f_{1}f_{2})^{-1}]\end{equation*}as Huber rings, so, in particular, \begin{equation*}\mathcal{A}[f_{1}^{-1}][f^{-1}]_{\gr,0}=\mathcal{A}[(f_{1}f_{2})^{-1}]_{\gr,0}.\end{equation*}But $f\in \mathcal{A}[f_{1}^{-1}]$ is a homogeneous element of degree $0$, so \begin{equation*}\mathcal{A}[f_{1}^{-1}][f^{-1}]_{\gr,0}=\mathcal{A}[f_{1}^{-1}]_{\gr,0}\langle f^{-1}\rangle\end{equation*}as Huber rings by Lemma~\ref{Pro-graded rational localizations, special case}. The assertion of the lemma follows.\end{proof}
\begin{lemma}\label{Gluing 1, part 4}Set assumptions as in Lemma~\ref{Gluing 0}. Let $(f_{i})_{i\in S} \subseteq \mathcal{A}_{\gr}^{+}$ be a family of homogeneous elements of degree~$1$ and, in the case $G=\mathbb{Z}[1/p]$, assume that each $f_{i}$ admits a compatible system of $p$-power roots in $\mathcal{A}_{\gr}^{+}$. Endow \begin{equation*}\Proj(\mathcal{A}, \mathcal{A}_{\gr}, \mathcal{A})(\frac{1}{f_{i}})\end{equation*}with the structure of an adic space by identifying it with $\Spa(\mathcal{A}\langle f_{i}^{-1}\rangle_{\gr,0}, \mathcal{A}\langle f_{i}^{-1}\rangle_{\gr,0}^{+})$ via the homeomorphism from Lemma~\ref{Gluing 1}. For $i, j\in S$, $i\neq j$, let $U_{ij}$ be the adic space whose underlying topological space is \begin{equation*}\vert U_{ij}\vert=\Proj(\mathcal{A}, \mathcal{A}_{\gr}, \mathcal{A})(\frac{1}{f_{1}})\cap \Proj(\mathcal{A}, \mathcal{A}_{\gr}, \mathcal{A})(\frac{1}{f_{2}})\end{equation*}and whose structure sheaf is the restriction of the structure sheaf of $\Spa(\mathcal{A}\langle f_{j}^{-1}\rangle_{\gr,0}, \mathcal{A}\langle f_{j}^{-1}\rangle_{\gr,0}^{+})$ along the inclusion map \begin{equation*}\vert U_{ij}\vert\hookrightarrow \Spa(\mathcal{A}\langle f_{j}^{-1}\rangle_{\gr,0}, \mathcal{A}\langle f_{j}^{-1}\rangle_{\gr,0}^{+}).\end{equation*}Set \begin{equation*}f_{ij}=\frac{f_{i}}{f_{j}}\in \mathcal{A}\langle f_{j}^{-1}\rangle_{\gr,0}.\end{equation*}Then $U_{ij}$ is the rational subspace \begin{equation*}U_{ij}=\Spa(\mathcal{A}\langle f_{j}^{-1}\rangle_{\gr,0}\langle f_{ij}^{-1}\rangle, \mathcal{A}\langle f_{j}^{-1}\rangle_{\gr,0}\langle f_{ij}^{-1}\rangle^{+})\end{equation*}and there are isomorphisms of adic spaces\begin{equation*}\varphi_{ij}: U_{ij}\cong \Spa(\mathcal{A}\langle (f_{i}f_{j})^{-1}\rangle_{\gr,0}, \mathcal{A}\langle (f_{i}f_{j})^{-1}\rangle_{\gr,0}^{+})\cong U_{ji}\end{equation*}which satisfy the cocycle condition on triple intersections $U_{ji}\cap U_{jk}$. \end{lemma} 
\begin{proof}By construction of the homeomorphism from Lemma~\ref{Gluing 1}, the open subset $\vert U_{ij}\vert$ of $\Spa(\mathcal{A}\langle f_{j}^{-1}\rangle_{\gr,0}, \mathcal{A}\langle f_{j}^{-1}\rangle_{\gr,0}^{+})$ is precisely the rational subset defined by the inequality $\vert f_{ij}\vert\geq1$. This implies the first assertion of the lemma. The isomorphisms $\varphi_{ij}$ are then obtained from Lemma~\ref{Gluing 1, part 3}, and it only remains to show that these isomorphisms satisfy the cocycle condition on triple overlaps. But, for any $i, j, k$, isomorphisms $\varphi_{ij}$, $\varphi_{jk}$ and $\varphi_{ik}$ are induced by equalities of subrings of $\mathcal{A}[(f_{i}f_{j}f_{k})^{-1}]$ \begin{equation*}\mathcal{A}[f_{j}^{-1}]_{\gr,0}[f_{ij}^{-1}]=\mathcal{A}[(f_{i}f_{j})^{-1}]_{\gr,0}=\mathcal{A}[f_{i}^{-1}]_{\gr,0}[f_{ji}^{-1}],\end{equation*}\begin{equation*}\mathcal{A}[f_{k}^{-1}]_{\gr,0}[f_{kj}^{-1}]=\mathcal{A}[(f_{j}f_{k})^{-1}]_{\gr,0}=\mathcal{A}[f_{j}^{-1}]_{\gr,0}[f_{jk}^{-1}]\end{equation*}and\begin{equation*}\mathcal{A}[f_{k}^{-1}]_{\gr,0}[f_{ki}^{-1}]=\mathcal{A}[(f_{i}f_{k})^{-1}]_{\gr,0}=\mathcal{A}[f_{i}^{-1}]_{\gr,0}[f_{ik}^{-1}],\end{equation*}which are also equalities of Huber rings by Lemma~\ref{Gluing 1, part 3}, so the cocycle condition is satisfied. \end{proof}
We introduce an analog of the notion of a sheafy Huber pair in the pro-graded context. 
\begin{mydef}[Graded-sheafy pro-graded Huber pair]\label{Sheafy}Let $G=\mathbb{Z}$ or $G=\mathbb{Z}[1/p]$. We call a pro-$G_{\geq0}$-graded Huber pair $(\mathcal{A}, \mathcal{A}_{\gr}, \mathcal{A}^{+})$ graded-sheafy if \begin{enumerate}[(1)]\item $\mathcal{A}$ is sheafy, and \item there exists a family of homogeneous elements of degree~$1$  \begin{equation*}(f_{s})_{s\in S}\subseteq \mathcal{A}_{\gr}^{+},\end{equation*}generating the irrelevant ideal $(\mathcal{A}, \mathcal{A}_{\gr})_{+}$ up to radicals and topological closure and, in the case $G=\mathbb{Z}[1/p]$, admitting compatible systems of $p$-power roots in $\mathcal{A}_{\gr}^{+}$, such that, for every $s\in S$, the Huber pair \begin{equation*}(\mathcal{A}\langle f_{s}^{-1}\rangle_{\gr, 0}, \mathcal{A}\langle f_{s}^{-1}\rangle_{\gr, 0}^{+})\end{equation*}is sheafy. \end{enumerate}In this case, we call a family of homogeneous elements $(f_{s})_{s\in S}$ as above a defining family of homogeneous elements of degree~$1$ for the pro-$G_{\geq0}$-graded Huber pair $((\mathcal{A}, \mathcal{A}_{\gr}), \mathcal{A}^{+})$. \end{mydef}
\begin{mydef}[The projective adic spectrum, as an adic space]\label{Adic Proj, part 2}Suppose that either $G=\mathbb{Z}$ or $G=\mathbb{Z}[1/p]$. Let \begin{equation*}((\mathcal{A}, \mathcal{A}_{\gr}), \mathcal{A}^{+})\end{equation*}be a graded-sheafy pro-$G_{\geq0}$-graded Huber pair and let $(f_{s})_{s\in S}$ be a defining family of homogeneous elements of degree~$1$. Then the projective adic spectrum of $((\mathcal{A}, \mathcal{A}_{\gr}), \mathcal{A}^{+})$, denoted once again by \begin{equation*}\Proj(\mathcal{A}, \mathcal{A}_{\gr}, \mathcal{A}^{+}),\end{equation*}is the adic space whose underlying topological space is the topological space from Definition~\ref{Adic Proj, part 1} and whose structure sheaf is obtained by gluing the pullbacks, under the homeomorphisms from Lemma~\ref{Gluing 0}, of the structure sheaves of the adic spaces\begin{equation*}\Spa(\mathcal{A}\langle f_{s}^{-1}\rangle_{\gr, 0}, \mathcal{A}\langle f_{s}^{-1}\rangle_{\gr, 0}^{+})\end{equation*}along the isomorphisms $\varphi_{ij}$ of Lemma~\ref{Gluing 1, part 4}.\end{mydef}
\begin{rmk}By virtue of Lemma~\ref{Gluing 1} and Lemma~\ref{Gluing 1, part 4}, this definition does not depend on the choice of the defining family $(f_{s})_{s\in S}$.\end{rmk}
\begin{example}[Adic projective space]\label{Adic projective space}If $K$ is a nonarchimedean field and $G=\mathbb{Z}$, consider the pro-$\mathbb{Z}_{\geq0}$-graded Tate Huber pair \begin{equation*}((T_{n+1, K}, (T_{n+1, K})_{\gr}), T_{n+1, K}^{\circ}),\end{equation*}where $(T_{n+1, K}, (T_{n+1, K})_{\gr})$ is the pro-graded Tate algebra from Example~\ref{Pro-graded Tate algebra}. The variables $X_{1},\dots, X_{n+1}$ are homogeneous elements of degree~$1$ generating the irrelevant ideal and, for every $i=1,\dots, n+1$, \begin{equation*}T_{n+1, K}\langle X_{i}^{-1}\rangle_{\gr, 0}\cong T_{n, K}\end{equation*}via the map sending $\frac{X_{1}}{X_{i}}, \dots, \frac{X_{i-1}}{X_{i}}, \frac{X_{i+1}}{X_{i}}, \dots, \frac{X_{n+1}}{X_{i}}$ to $Y_{1},\dots, Y_{i-1}, Y_{i},\dots, Y_{n}$. In particular, we see that $((T_{n+1, K}, (T_{n+1, K})_{\gr}), T_{n+1, K}^{\circ})$ is graded-sheafy. Moreover, \begin{equation*}\Proj(T_{n+1, K}, (T_{n+1, K})_{\gr}, T_{n, K}^{\circ})=\mathbb{P}_{K}^{n,\ad},\end{equation*}the usual adic projective space over $K$ constructed by $n+1$ copies of the $n$-dimensional closed unit ball along their common annuli. \end{example}
\begin{example}[Perfectoid projective space]If $K$ is a perfectoid field and $G=\mathbb{Z}[1/p]$, we can analogously consider the pro-$\mathbb{Z}[1/p]_{\geq0}$-graded Tate Huber pair \begin{equation*}((T_{n+1, K}^{\perfd}, (T_{n+1, K}^{\perfd})_{\gr}), T_{n+1, K}^{\perfd\circ}),\end{equation*}where $(T_{n+1, K}^{\perfd}, (T_{n+1, K}^{\perfd})_{\gr})$ is the pro-graded perfectoid Tate algebra from Example~\ref{Pro-graded perfectoid Tate algebra}. The variables $X_1,\dots, X_n$ are homogeneous elements of degree~$1$ which admit compatible systems of $p$-power roots and which generate the irrelevant ideal up to radical and topological closure. Moreover, for every $i=1,\dots, n+1$,\begin{equation*}T_{n+1, K}^{\perfd}\langle X_{i}^{-1}\rangle_{\gr, 0}\cong T_{n, K}^{\perfd},\end{equation*}once again via the map taking \begin{equation*}\frac{X_{1}^{1/p^{m}}}{X_{i}^{1/p^{m}}},\dots, \frac{X_{i-1}^{1/p^{m}}}{X_{i}^{1/p^{m}}}, \frac{X_{i+1}^{1/p^{m}}}{X_{i}^{1/p^{m}}},\dots, \frac{X_{n+1}^{1/p^{m}}}{X_{i}^{1/p^{m}}}\end{equation*}to $Y_{1}^{1/p^{m}},\dots, Y_{i}^{1/p^{m}}, Y_{i+1}^{1/p^{m}},\dots, Y_{n}^{1/p^{m}}$ for all $m\geq0$. Then \begin{equation*}\Proj(T_{n+1, K}^{\perfd}, (T_{n+1, K}^{\perfd})_{\gr}, T_{n+1, K}^{\perfd\circ})=\mathbb{P}_{K}^{n,\perfd},\end{equation*}the perfectoid projective space over $K$ constructed by gluing $n+1$ copies of the perfectoid $n$-dimensional closed unit ball along their common (perfectoid) annuli.\end{example}

\section{Relation to formal schemes}\label{sec:formal}

The adic $\Proj$ construction also has an obvious cousin in the world of formal schemes. Note that an open ideal of a pro-$G$-graded adic ring is topologically homogeneous if and only if it is homogeneous. A similar notion of (formal) projective spectrum for pro-$G$-graded adic rings was also introduced by Ishizuka and Yoshikawa in \cite{Ishizuka-Yoshikawa26}, Appendix~A.2.   
\begin{mydef}[Formal projective spectrum, as a topological space]\label{Formal Proj, part 1}For a pro-$G_{\geq0}$-graded adic ring $(R, R_{\gr})$, where $G$ is a totally ordered abelian group, let \begin{equation*}\Proj(R, R_{\gr})\end{equation*}be the subspace of $\Spf(R)$ consisting of (topologically) homogeneous open prime ideals $\mathfrak{p}$ not containing the irrelevant ideal $(R, R_{\gr})_{+}$ of $(R, R_{\gr})$.\end{mydef}
As in the adic case, this recovers the usual formal projective space over an adic ring.
\begin{lemma}\label{Gluing 0, formal case}Suppose that $G=\mathbb{Z}$ or $G=\mathbb{Z}[1/p]$. Let $(R, R_{\gr})$ be a pro-$G_{\geq0}$-graded adic ring and let $f\in R_{\gr}$ be a homogeneous element of degree~$1$. In the case $G=\mathbb{Z}[1/p]$, assume that $f$ admits a compatible system of $p$-power roots in $R_{\gr}$. Then the map \begin{equation*}D_{+}(f)=\{\, \mathfrak{p}\in \Proj(R, R_{\gr})\mid f\not\in \mathfrak{p}\,\}\to \Spf(R\langle f^{-1}\rangle_{\gr, 0}), \mathfrak{p}\mapsto \mathfrak{p}R\langle f^{-1}\rangle\cap R\langle f^{-1}\rangle_{\gr, 0},\end{equation*}is a homeomorphism.\end{lemma}
\begin{proof}For an open prime ideal $\mathfrak{p}$ of $R\langle f^{-1}\rangle_{\gr, 0}$ and $g\in G$, we define \begin{equation*}\mathfrak{p}_{g}=\{\, h\in R\langle f^{-1}\rangle_{\gr, g}\mid \frac{h}{f^{g}}\in \mathfrak{p}\,\}.\end{equation*}Consider the additive subgroup $\bigoplus_{g\in G}\mathfrak{p}_{g}$ of $R\langle f^{-1}\rangle_{\gr}$. We claim that this is an open homogeneous ideal of $R\langle f^{-1}\rangle_{\gr}$, that its intersection with $R_{\gr}$ is an open homogeneous prime ideal and that the map \begin{equation*}\psi: \Spf(R\langle f^{-1}\rangle_{\gr, 0})\to D_{+}(f), \mathfrak{p}\mapsto (\bigoplus_{g\in G}\mathfrak{p}_{g})R\langle f^{-1}\rangle\cap R,\end{equation*}is an inverse to the map \begin{equation*}\varphi: \mathfrak{q}\mapsto \mathfrak{q}R\langle f^{-1}\rangle\cap R\langle f^{-1}\rangle_{\gr, 0}.\end{equation*}.

We first prove that $\bigoplus_{g\in G}\mathfrak{p}_{g}$ is an ideal. If $h\in \bigoplus_{g\in G}\mathfrak{p}_{g}$ and $a\in R\langle f^{-1}\rangle_{\gr}$, with homogeneous decompositions $h=\sum_{g\in G}h_{g}$ and $a=\sum_{g\in G}a_{g}$, then the homogeneous decomposition of $ah$ is \begin{equation*}ah=\sum_{g\in G}\sum_{g_{1}+g_{2}=g}a_{g_{1}}h_{g_{2}}.\end{equation*}Since $h\in \bigoplus_{g\in G}\mathfrak{p}_{g}$, we have $\frac{h_{g_{2}}}{f^{g_{2}}}\in \mathfrak{p}$, so \begin{equation*}\frac{a_{g_{1}}h_{g_{2}}}{f^{g_{1}+g_{2}}}=\frac{a_{g_{1}}}{f^{g_{1}}}\frac{h_{g_{2}}}{f^{g_{2}}}\in \mathfrak{p}\end{equation*}since $\mathfrak{p}$ is an ideal. Thus $a_{g_{1}}h_{g_{2}}\in \mathfrak{p}_{g_{1}+g_{2}}$ for every $g_{1}, g_{2}\in G$ and thus $ah\in \bigoplus_{g\in G}\mathfrak{p}_{g}$.

This shows that $\bigoplus_{g\in G}\mathfrak{p}_{g}$ is an ideal of $R\langle f^{-1}\rangle_{\gr}$. By definition, this ideal is homogeneous. Next, we show that $\bigoplus_{g\in G}\mathfrak{p}_{g}$ is open. Let $I$ be a finitely generated, homogeneous ideal of definition of $R_{\gr}$. Since $\mathfrak{p}$ is open in $R\langle f^{-1}\rangle_{\gr, 0}$, there exists an integer $N\geq1$ such that \begin{equation*}I^{N}R\langle f^{-1}\rangle\cap R\langle f^{-1}\rangle_{\gr, 0}\subseteq \mathfrak{p}.\end{equation*}For any homogeneous element $a\in I$ of degree $g\in G$, we have \begin{equation*}(\frac{a}{f^{g}})^{N}\in I^{N}R\langle f^{-1}\rangle\cap R\langle f^{-1}\rangle_{\gr, 0}\subseteq \mathfrak{p}.\end{equation*}Since $\mathfrak{p}$ is prime, this implies $\frac{a}{f^{g}}\in \mathfrak{p}$, so $a\in \mathfrak{p}_{g}$. Since $I$ is generated by homogeneous elements, this shows that $\bigoplus_{g\in G}\mathfrak{p}_{g}$ contains $IR\langle f^{-1}\rangle_{\gr}$, i.e., it is open. It follows that $(\bigoplus_{g\in G}\mathfrak{p}_{g})\cap R_{\gr}$ is an open homogeneous ideal of $R_{\gr}$ (the map $R_{\gr}\to R\langle f^{-1}\rangle_{\gr}$ being a continuous homomorphism of $G$-graded rings). 

Next, we show that the ideal $(\bigoplus_{g\in G}\mathfrak{p}_{g})\cap R_{\gr}$ is prime. Let $f_{1}, f_{2}\in R_{\gr}$ with homogeneous decompositions $f_{1}=\sum_{g\in G}f_{1, g}$, $f_{2}=\sum_{g\in G}f_{2, g}$ such that \begin{equation*}f_{1}, f_{2}\not\in \bigoplus_{g\in G}\mathfrak{p}_{g}.\end{equation*}Let \begin{equation*}g_{i}=\min\{g\in G\mid f_{i, g}\not\in \mathfrak{p}_{g}\}, i=1, 2.\end{equation*}
To prove that $f_{1}f_{2}\not\in \bigoplus_{g\in G}\mathfrak{p}_{g}$, we prove that $(f_{1}f_{2})_{g_{1}+g_{2}}\not\in \mathfrak{p}_{g_{1}+g_{2}}$. Write \begin{equation*}(f_{1}f_{2})_{g_{1}+g_{2}}=\sum_{g+g^{\prime}=g_{1}+g_{2}}f_{1,g}f_{2, g^{\prime}}.\end{equation*}For every $(g, g^{\prime})\neq (g_1, g_2)$ with $g+g^{\prime}=g_{1}+g_{2}$, then $g<g_{1}$ or $g^{\prime}<g_{2}$, so \begin{equation*}f_{1, g}f_{2, g^{\prime}}\in \mathfrak{p}_{g+g^{\prime}}=\mathfrak{p}_{g_{1}+g_{2}}\end{equation*}by our choice of $g_{1}$ and $g_{2}$. Hence if $(f_{1}f_{2})_{g_{1}+g_{2}}\in \bigoplus_{g\in G}\mathfrak{p}_{g}$, then $f_{1, g_{1}}f_{2, g_{2}}\in \mathfrak{p}_{g_{1}+g_{2}}$, contradicting the definition of $g_{1}$ and $g_{2}$. It follows that $(f_{1}f_{2})_{g_{1}+g_{2}}\not\in \mathfrak{p}_{g_{1}+g_{2}}$ and hence $f_{1}f_{2}\not\in \bigoplus_{g\in G}\mathfrak{p}_{g}$. This shows that $(\bigoplus_{g\in G}\mathfrak{p}_{g})\cap R_{\gr}$ is a prime ideal of $R_{\gr}$, as desired.

Since $\bigoplus_{g\in G}\mathfrak{p}_{g}$ is an open ideal of $R\langle f^{-1}\rangle_{\gr}$, we have \begin{equation*}(\bigoplus_{g\in G}\mathfrak{p}_{g})R\langle f^{-1}\rangle\cap R\langle f^{-1}\rangle_{\gr}=\bigoplus_{g\in G}\mathfrak{p}_{g}.\end{equation*}By definition, \begin{equation*}(\bigoplus_{g\in G}\mathfrak{p}_{g})\cap R\langle f^{-1}\rangle_{\gr, 0}=\mathfrak{p}_{0}=\mathfrak{p},\end{equation*}so \begin{equation*}(\bigoplus_{g\in G}\mathfrak{p}_{g})R\langle f^{-1}\rangle\cap R\langle f^{-1}\rangle_{\gr, 0}=\mathfrak{p}.\end{equation*}Moreover, we also have \begin{equation*}(\bigoplus_{g\in G}\mathfrak{p}_{g})R\langle f^{-1}\rangle\cap R=((\bigoplus_{g\in G}\mathfrak{p}_{g})\cap R_{\gr})R,\end{equation*}as we have seen that $(\bigoplus_{g\in G}\mathfrak{p}_{g})\cap R_{\gr}$ is an open ideal of $R_{\gr}$. Since we have seen that the open homogeneous ideal $(\bigoplus_{g\in G}\mathfrak{p}_{g})\cap R_{\gr}$ of $R_{\gr}$ is prime, we conclude that $(\bigoplus_{g\in G}\mathfrak{p}_{g})R\langle f^{-1}\rangle\cap R$ is an open homogeneous prime ideal of $R$. It then follows that \begin{equation*}(\bigoplus_{g\in G}\mathfrak{p}_{g})R\langle f^{-1}\rangle=((\bigoplus_{g\in G}\mathfrak{p}_{g})R\langle f^{-1}\rangle\cap R)R\langle f^{-1}\rangle.\end{equation*}Putting everything together, we conclude that $\varphi(\psi(\mathfrak{p}))=\mathfrak{p}$ for all open prime ideals $\mathfrak{p}$ of $R\langle f^{-1}\rangle_{\gr,0}$. 

Conversely, let $\mathfrak{q}\in D_{+}(f)$ and consider the homogeneous open prime ideal $\psi(\varphi(\mathfrak{q}))$ of $R$. To prove that $\mathfrak{q}=\psi(\varphi(\mathfrak{q}))$, it suffices to prove that these open prime ideals are equal after restricting to the dense subring $R_{\gr}$ of $R$. We verify: \begin{align*}\psi(\varphi(\mathfrak{q}))=\psi(\mathfrak{q}R\langle f^{-1}\rangle\cap R\langle f^{-1}\rangle_{\gr,0})\cap R_{\gr} \\ =\bigoplus_{g\in G}(\mathfrak{q}R\langle f^{-1}\rangle\cap R\langle f^{-1}\rangle_{\gr,0})\cap R_{\gr} \\ =\{\, h=\sum_{g\in G}h_{g}\in R_{\gr}\mid \frac{h_{g}}{f^{g}}\in \mathfrak{q}R\langle f^{-1}\rangle\cap R\langle f^{-1}\rangle_{\gr,0}\,\} \\ =\{\, h=\sum_{g\in G}h_{g}\in R_{\gr}\mid h_{g}\in \mathfrak{q}R\langle f^{-1}\rangle \cap R_{\gr}\,\} \\ =\{\, h=\sum_{g\in G}h_{g}\in R_{\gr}\mid h_{g}\in \mathfrak{q}\cap R_{\gr}\,\}=\mathfrak{q}\cap R_{\gr},\end{align*}where we used at multiple steps that $\mathfrak{q}$ is open in $R$ and where we used at the last step that $\mathfrak{q}\cap R_{\gr}$ is homogeneous. Thus the maps $\varphi$ and $\psi$ are inverses, as claimed. 

It remains to prove that the map $\psi$ is continuous. But to this end, we can observe that none of our above arguments (except the argument showing that $\bigoplus_{g\in G}\mathfrak{p}_{g}$ is prime for an open prime ideal $\mathfrak{p}$ of $R\langle f^{-1}\rangle_{\gr,0}$) used the fact that $\mathfrak{p}$ and $\mathfrak{q}$ are prime ideals of $R\langle f^{-1}\rangle_{\gr,0}$ and $R$, respectively, only that they are open ideals and that $\mathfrak{q}$ is homogeneous. Hence the mutually inverse bijections $\varphi$ and $\psi$ extend to inclusion-preserving, mutually inverse bijections between the set of open homogeneous ideals of $R$ not containing $f$ and the set of open ideals of $R\langle f^{-1}\rangle_{\gr,0}$. Consequently, $\varphi$ and $\psi$ are mutually inverse homeomorphisms. \end{proof}
\begin{lemma}\label{Gluing 1, formal case}Suppose that either $G=\mathbb{Z}$ or $G=\mathbb{Z}[1/p]$. Let $(R, R_{\gr})$ be a pro-$G_{\geq0}$-graded adic ring. Let $f_1, f_2\in R_{\gr}$ be homogeneous elements of degree~$1$. In the case $G=\mathbb{Z}[1/p]$, assume that $f_{1}$ and $f_{2}$ admit compatible systems of $p$-power roots in $R_{\gr}$. Then each of the maps \begin{equation*}D_{+}(f_{i})\to \Spf(R\langle f^{-1}\rangle_{\gr, 0}), i=1, 2,\end{equation*}from Lemma~\ref{Gluing 0, formal case} restrict to the same map \begin{equation*}D_{+}(f_{1})\cap D_{+}(f_{2})=D_{+}(f_{1}f_{2})\to \Spf(R\langle (f_{1}f_{2})^{-1}\rangle_{\gr, 0}).\end{equation*}and this map is a homeomorphism.\end{lemma}
\begin{proof}Follows from Lemma~\ref{Gluing 0, formal case}.\end{proof}
\begin{lemma}\label{Gluing 1, formal case, part 2}Suppose that $G=\mathbb{Z}$ or $G=\mathbb{Z}[1/p]$. Let $(R, R_{\gr})$ be a pro-$G_{\geq0}$-graded adic ring. Let $(f_{i})_{i\in S} \subseteq R_{\gr}$ be a family of homogeneous elements of degree~$1$. In the case $G=\mathbb{Z}[1/p]$, suppose that each $f_{i}$ admits a compatible system of $p$-power roots in $R_{\gr}$. Endow each \begin{equation*}D_{+}(f_{i})=\{\, \mathfrak{p}\in \Proj(R, R_{\gr})\mid f_{i}\not\in \mathfrak{p}\,\}\end{equation*}with the structure of a formal scheme by identifying it with $\Spf(R\langle f_{i}^{-1}\rangle_{\gr,0})$ via the homeomorphism from Lemma~\ref{Gluing 1, formal case}. For $i, j\in S$, $i\neq j$, let $\mathfrak{U}_{ij}$ be the formal scheme whose underlying topological space is \begin{equation*}\vert \mathfrak{U}_{ij}\vert=D_{+}(f_{i})\cap D_{+}(f_{j})=D_{+}(f_{i}f_{j})\end{equation*}and whose structure sheaf is the restriction of the structure sheaf of $\Spf(R\langle f_{j}^{-1}\rangle_{\gr,0})$ along the inclusion map \begin{equation*}\vert \mathfrak{U}_{ij}\vert\hookrightarrow \Spf(R\langle f_{j}^{-1}\rangle_{\gr,0}).\end{equation*}Set \begin{equation*}f_{ij}=\frac{f_{i}}{f_{j}}\in R\langle f_{j}^{-1}\rangle_{\gr,0}.\end{equation*}Then $\mathfrak{U}_{ij}$ is the open formal subscheme \begin{equation*}\mathfrak{U}_{ij}=\Spf(R\langle f_{j}^{-1}\rangle_{\gr,0}\langle f_{ij}^{-1}\rangle)\end{equation*}and there are isomorphisms of formal schemes\begin{equation*}\varphi_{ij}: \mathfrak{U}_{ij}\cong \Spf(R\langle (f_{i}f_{j})^{-1}\rangle_{\gr,0}\cong \mathfrak{U}_{ji}\end{equation*}which satisfy the cocycle condition on triple intersections $\mathfrak{U}_{ji}\cap \mathfrak{U}_{jk}$. \end{lemma} 
\begin{proof}For an open prime ideal $\mathfrak{p}$ of $R$ which does not contain $f_{j}$, we have $f_{i}\not\in \mathfrak{p}$ if and only \begin{equation*}f_{ij}\not\in \mathfrak{p}R\langle f_{j}^{-1}\rangle\cap R\langle f_{j}^{-1}\rangle_{\gr,0}.\end{equation*}
Hence the homeomorphism from Lemma~\ref{Gluing 1, formal case} identifies $\vert \mathfrak{U}_{ij}\vert$ with the open subset \begin{equation*}\Spf(R\langle f_{j}^{-1}\rangle_{\gr,0}\langle f_{ij}^{-1}\rangle)=\{\, \mathfrak{p}\in \Spf(R\langle f_{j}^{-1}\rangle_{\gr,0}\mid f_{ij}\not\in \mathfrak{p}\,\}\end{equation*}of $\Spf(R\langle f_{j}^{-1}\rangle_{\gr,0})$. The isomorphisms $\varphi_{ij}$ are then obtained from Lemma~\ref{Gluing 1, part 3}. The cocycle condition on triple intersections follows from equalities of subrings of $R[(f_{i}f_{j}f_{k})^{-1}]$ based on Lemma~\ref{Gluing 1, part 3}, directly analogous to those displayed at the end of the proof of Lemma~\ref{Gluing 1, part 4}. \end{proof}
\begin{mydef}[The formal projective spectrum of a pro-graded adic ring]\label{Formal Proj, part 2}Suppose that either $G=\mathbb{Z}$ or $G=\mathbb{Z}[1/p]$. Let $(R, R_{\gr})$ be a pro-$G_{\geq0}$-graded adic ring and let $(f_{s})_{s\in S}$ be a family of homogeneous elements of degree~$1$ in $R_{\gr}$ generating the irrelevant ideal of $(R, R_{\gr})$ up to radical and topological closure and, in the case $G=\mathbb{Z}[1/p]$, admitting compatible systems of $p$-power roots in $R_{\gr}$. 

We endow $\Proj(R, R_{\gr})$ with the structure of a formal scheme by identifying the open subsets $D_{+}(f_{s})$, $s\in S$, with affine formal schemes $\Spf(R\langle f_{s}^{-1}\rangle_{\gr, 0})$, $s\in S$, via the homeomorphism from Lemma~\ref{Gluing 0, formal case} and gluing their intersections along the isomorphisms $\varphi_{ij}$ from Lemma~\ref{Gluing 1, formal case, part 2}. We call the resulting formal scheme $\Proj(R, R_{\gr})$ the formal projective spectrum of the pro-$G$-graded adic ring $(R, R_{\gr})$. \end{mydef}
\begin{example}[Formal projective space]\label{Formal projective space}Let $A$ be a fixed complete adic ring with finitely generated ideal of definition $I$ and with trivial grading and define the pro-graded adic ring $(A\langle X_1,\dots, X_{n+1}\rangle, (A\langle X_1,\dots, X_{n+1}\rangle)_{\gr})$ as in Remark~\ref{Convergent power series, general case}. Then \begin{equation*}\Proj(A\langle X_{1},\dots, X_{n+1}\rangle, (A\langle X_{1},\dots, X_{n+1}\rangle)_{\gr})=\widehat{\mathbb{P}}_{A}^{n},\end{equation*}the $I$-adic formal completion of the usual (algebraic) projective space $\mathbb{P}_{A}^{n}$ over $A$. This is readily seen directly from the description of the formal scheme structure in Definition~\ref{Formal Proj, part 2}. There is also an obvious ``perfectoid'' analog of this example, which we omit spelling out here.\end{example}
We can consider the specialization map relating the adic and formal $\Proj$ constructions to each other. Note that if $((\mathcal{A}, \mathcal{A}_{\gr}), \mathcal{A}^{+})$ is a uniform pro-$G_{\geq0}$-graded Huber pair, then, by the uniformity assumption, $\mathcal{A}^{+}$, with the subspace topology, is an adic ring and, by the definition of a pro-$G$-graded Huber pair, $(\mathcal{A}^{+}, \mathcal{A}_{\gr}^{+})$ is a pro-$G_{\geq0}$-graded adic ring. In fact, $(\mathcal{A}^{+}, \mathcal{A}_{\gr}^{+})$ is a pro-$G$-graded Huber ring of definition of $(\mathcal{A}, \mathcal{A}_{\gr})$. In particular, we can speak of the formal projective spectrum of $(\mathcal{A}^{+}, \mathcal{A}_{\gr}^{+})$. 
\begin{lemma}\label{Specialization map, definition}Suppose that $G=\mathbb{Z}$ or $G=\mathbb{Z}[1/p]$ and that $((\mathcal{A}, \mathcal{A}_{\gr}), \mathcal{A}^{+})$ is a uniform graded-sheafy pro-$G_{\geq0}$-graded Tate Huber pair. Then the specialization map \begin{equation*}\Spa(\mathcal{A}, \mathcal{A}^{+})\to \Spf(\mathcal{A}^{+}), v\mapsto \spc(v)=\{\, f\in \mathcal{A}^{+}\mid v(f)<1\,\},\end{equation*}restricts to a continuous map \begin{equation*}\spc: \Proj(\mathcal{A}, \mathcal{A}_{\gr}, \mathcal{A}^{+})\to \Proj(\mathcal{A}^{+}, \mathcal{A}_{\gr}^{+}).\end{equation*}\end{lemma}
\begin{proof}If a continuous valuation $v\in \Spa(\mathcal{A}, \mathcal{A}^{+})$ is topologically homogeneous in the sense of Definition~\ref{Topologically homogeneous valuations}, then its center on $\mathcal{A}^{+}$ is a homogeneous prime ideal of $\mathcal{A}^{+}$. Moreover, a basis for the topology on $\Proj(\mathcal{A}^{+}, \mathcal{A}_{\gr}^{+})$ is given by basic open subsets $D_{+}(f)$ with $f\in \mathcal{A}_{\gr}^{+}$ homogeneous, and pre-images of such basic open subsets under $\spc$ are the homogeneous rational subsets $\Proj(\mathcal{A}, \mathcal{A}_{\gr}, \mathcal{A}^{+})(\frac{1}{f})$, which are open in $\Proj(\mathcal{A}, \mathcal{A}_{\gr}, \mathcal{A}^{+})$. \end{proof}
\begin{mydef}\label{Specialization map, actual definition}We call this map the specialization map on $\Proj(\mathcal{A}, \mathcal{A}_{\gr}, \mathcal{A}^{+})$.\end{mydef} 
As in our previous work \cite{Dine25-2}, Definition~2.16 and Lemma~2.14, we call a formal scheme $\mathfrak{X}$ with a torsion-free, adic morphism $\mathfrak{X}\to \Spf(\mathbb{Z}[[\varpi]])$ locally rig-sheafy (over $\mathbb{Z}[[\varpi]]$) if it is covered by open affines $\Spf(A)$, with ideal of definition generated by $\varpi$, such that the Tate ring $A[\varpi^{-1}]$, with pair of definition $(A, \varpi)$, is sheafy. It is clear that in the situation of Lemma~\ref{Specialization map, definition}, $\Proj(\mathcal{A}^{+}, \mathcal{A}_{\gr}^{+})$ is locally rig-sheafy over $\mathbb{Z}[[\varpi]]$, where $\varpi$ is a pseudo-uniformizer of the Tate Huber pair $(\mathcal{A}, \mathcal{A}^{+})$. For any locally rig-sheafy, $\varpi$-torsion-free adic formal $\mathbb{Z}[[\varpi]]$-scheme \begin{equation*}\mathfrak{X}\to \Spf(\mathbb{Z}[[\varpi]]),\end{equation*}one can define its adic analytic generic fiber $\mathfrak{X}_{\eta}^{\ad}$, an analytic adic space, equipped with a specialization map \begin{equation*}\spc_{\mathfrak{X}_{\eta}^{\ad}, \mathfrak{X}}: \mathfrak{X}_{\eta}^{\ad}\to \mathfrak{X},\end{equation*}which is a morphism of locally $v$-ringed spaces, see \cite{Dine25-2}, Lemma~3.8 and Definition~3.9. If $\mathfrak{X}=\Spf(A)$ is affine, where $A$ is a complete adic ring with ideal of definition generated by the non-zero-divisor $\varpi$, then \begin{equation*}\mathfrak{X}_{\eta}^{\ad}=\Spa(\mathcal{A}, \mathcal{A}^{+}),\end{equation*}with $\mathcal{A}=A[\varpi^{-1}]$ and with $\mathcal{A}^{+}$ the integral closure of $A$ inside $\mathcal{A}$. In the setting of rigid-analytic geometry (when the formal schemes in question are topologically of finite type over a rank-$1$ valuation ring), these constructions are classical and due to Raynaud~\cite{RaynaudReport} (in a slightly more general setup, they are due to Berthelot~\cite{Berthelot96}, \S0.2.6). We now relate these constructions to the specialization map introduced in Definition~\ref{Specialization map, definition}.   
\begin{thm}\label{Two specialization maps}Let the assumptions be as in Lemma~\ref{Specialization map, definition}. Then \begin{equation*}\Proj(\mathcal{A}, \mathcal{A}_{\gr}, \mathcal{A}^{+})=\Proj(\mathcal{A}^{+}, \mathcal{A}_{\gr}^{+})_{\eta}^{\ad}\end{equation*}and the specialization map \begin{equation*}\spc: \Proj(\mathcal{A}, \mathcal{A}_{\gr}, \mathcal{A}^{+})\to \Proj(\mathcal{A}^{+}, \mathcal{A}_{\gr}^{+})\end{equation*}coincides with the specialization map defined in \cite{Dine25-2}, Definition~3.9. \end{thm}
\begin{proof}Let $(f_{s})_{s\in S}$ be a defining family of homogeneous elements of degree~$1$, i.e., the elements $f_{s}\in \mathcal{A}_{\gr}^{+}$ generate the irrelevant ideal up to radical and topological closure (and admit compatible systems of $p$-power roots if $G=\mathbb{Z}[1/p]$). Consider the affine open cover \begin{equation*}D_{+}(f_{s})\simeq \Spf(\mathcal{A}^{+}\langle f_{s}^{-1}\rangle_{\gr, 0}), s\in S,\end{equation*}obtained from the definition of the formal scheme structure on $\Proj(\mathcal{A}^{+}, \mathcal{A}_{\gr}^{+})$. Note that the generic fiber of $\Spf(\mathcal{A}^{+}\langle f_{s}^{-1}\rangle_{\gr, 0})$ is precisely $\Spa(\mathcal{A}\langle f_{s}^{-1}\rangle_{\gr, 0}, \mathcal{A}\langle f_{s}^{-1}\rangle_{\gr, 0}^{+})$. Therefore, for every $s\in S$, we obtain a square of morphisms of locally $v$-ringed spaces \begin{center}\begin{tikzcd}\Proj(\mathcal{A}, \mathcal{A}_{\gr}, \mathcal{A}^{+})(\frac{1}{f_{s}})\arrow{d}{\simeq} \arrow{r}{\spc} & D_{+}(f_{s}) \arrow{d}{\simeq} \\ \Spa(\mathcal{A}\langle f_{s}^{-1}\rangle_{\gr, 0}, \mathcal{A}\langle f_{s}^{-1}\rangle_{\gr, 0}^{+}) \arrow{r}{\spc} & \Spf(\mathcal{A}\langle f_{s}^{-1}\rangle_{\gr, 0}^{+}),\end{tikzcd}\end{center}where the upper horizontal arrow is (the restriction to $\Proj(\mathcal{A}, \mathcal{A}_{\gr}, \mathcal{A}^{+})(\frac{1}{f_{s}})$) of the specialization map from Lemma~\ref{Specialization map, definition}, the lower horizontal arrow is the specialization map on the generic fiber as in \cite{Dine25-2}, Definition~3.9, and where the vertical arrows are the maps used to glue the adic projective spectrum and the formal projective spectrum, respectively. To finish the proof of the theorem, it suffices to check that this square is commutative. To this end, let $v\in \Proj(\mathcal{A}, \mathcal{A}_{\gr}, \mathcal{A}^{+})$ with $v(f_{s})=1$ and consider the image \begin{equation*}\spc(v)\mathcal{A}^{+}\langle f^{-1}\rangle\cap \mathcal{A}\langle f^{-1}\rangle_{\gr, 0}^{+}\end{equation*}of $v$ under the composite map \begin{equation*}\Proj(\mathcal{A}, \mathcal{A}_{\gr}, \mathcal{A}^{+})(\frac{1}{f_{s}})\to D_{+}(f_{s})\simeq \Spf(\mathcal{A}\langle f^{-1}\rangle_{\gr, 0}^{+}).\end{equation*}Since $v(f)=1$, we know that $v$ extends uniquely to a continuous valuation \begin{equation*}v'\in \Spa(\mathcal{A}\langle f^{-1}\rangle, \mathcal{A}\langle f^{-1}\rangle^{+})\end{equation*}by setting \begin{equation*}v'(\frac{a}{f^{n}})=v(a), a\in \mathcal{A}, n\geq1.\end{equation*}Then \begin{equation*}\spc(v)\mathcal{A}^{+}\langle f^{-1}\rangle=\spc(v')\end{equation*}as prime ideals of $\mathcal{A}\langle f^{-1}\rangle^{+}=\mathcal{A}^{+}\langle f^{-1}\rangle$ and thus \begin{equation*}\spc(v)\mathcal{A}^{+}\langle f^{-1}\rangle\cap \mathcal{A}\langle f^{-1}\rangle_{\gr, 0}^{+}=\spc(v'\vert_{\mathcal{A}\langle f^{-1}\rangle_{\gr, 0}}).\end{equation*}But by Lemma~\ref{Gluing -1}, $v'$ is topologically homogeneous and its restriction to $\mathcal{A}\langle f^{-1}\rangle_{\gr, 0}$ must coincide with the image of $v$ under the homeomorphism \begin{equation*}\Proj(\mathcal{A}, \mathcal{A}_{\gr}, \mathcal{A}^{+})(\frac{1}{f_{s}})\tilde{\to} \Spa(\mathcal{A}\langle f^{-1}\rangle_{\gr, 0}, \mathcal{A}\langle f^{-1}\rangle_{\gr, 0}^{+}).\end{equation*}This concludes the proof. \end{proof}
As an immediate consequence of Theorem~\ref{Two specialization maps}, we obtain the following inequality between the Krull dimensions of the locally spectral spaces $\Proj(\mathcal{A}, \mathcal{A}_{\gr}, \mathcal{A}^{+})$ and $\Proj(\mathcal{A}^{+}, \mathcal{A}_{\gr}^{+})$. 
\begin{cor}\label{Dimension estimate for the Proj}Set all assumptions as in Lemma~\ref{Specialization map, definition}. Then \begin{equation*}\dim(\Proj(\mathcal{A}, \mathcal{A}_{\gr}, \mathcal{A}^{+}))\geq \dim(\Proj(\mathcal{A}^{+}, \mathcal{A}_{\gr}^{+})).\end{equation*}\end{cor}
\begin{proof}The specialization map \begin{equation*}\spc_{\mathfrak{X}_{\eta}^{\ad},\mathfrak{X}}: \mathfrak{X}_{\eta}^{\ad}\to \mathfrak{X}\end{equation*}from the adic analytic generic fiber of a locally rig-sheafy $\varpi$-torsion-free adic formal $\mathbb{Z}[[\varpi]]$-scheme is always surjective and closed, by \cite{Dine25-2}, Corollary~4.16 and Corollary~4.20. Hence the assertion follows from Theorem~\ref{Two specialization maps}. \end{proof}

\section{Relation to projective varieties and GAGA}\label{sec:GAGA}

In this section, we want to describe adic analytifications of all projective varieties over a nonarchimedean field $K$ as instances of our adic $\Proj$ construction. Hence we specialize our discussion to pro-graded Huber pairs defined by the pro-graded affinoid algebras of Definition~\ref{Pro-graded affinoid algebras, definition}. 
\begin{lemma}\label{Pro-graded affinoids are graded-sheafy}The pro-$\mathbb{Z}_{\geq0}$-graded Tate Huber pair $((\mathcal{A}, \mathcal{A}_{\gr}), \mathcal{A}^{\circ})$ defined by any pro-graded affinoid $K$-algebra $(\mathcal{A}, \mathcal{A}_{\gr})$ is graded-sheafy. In particular, we can speak of the projective adic spectrum $\Proj(\mathcal{A}, \mathcal{A}_{\gr}, \mathcal{A}^{\circ})$ of $((\mathcal{A}, \mathcal{A}_{\gr}), \mathcal{A}^{\circ})$.\end{lemma}
\begin{proof}Let \begin{equation*}(\varphi, \varphi_{\gr}): (T_{n, K}, (T_{n, K})_{\gr})\to (\mathcal{A}, \mathcal{A}_{\gr})\end{equation*}be a morphism of pro-graded Huber rings such that $\varphi$ is surjective. If $n=0$, then $\mathcal{A}$ is a nonarchimedean field and there is nothing to prove. Assume that $n\geq1$. By Lemma~\ref{Surjective morphisms of pro-graded rings}, $\varphi_{\gr}$ is also surjective. Let $f_1,\dots, f_n$ be the images in $\mathcal{A}_{\gr}$ of the variables $X_1,\dots, X_n$. Since $\varphi_{\gr}$ is a surjective homomorphism of graded rings, the irrelevant ideal $\bigoplus_{g>0}\mathcal{A}_{\gr, g}$ is generated by the degree-$1$ homogeneous elements $f_1,\dots, f_n$. Since $\mathcal{A}$ is clearly sheafy, being an affinoid algebra, we only have to prove that $\mathcal{A}\langle f_{i}^{-1}\rangle_{\gr, 0}$ is sheafy for all $i=1,\dots, n$. Fix some $i\in\{1,\dots, n\}$. By the universal property of pro-graded rational localizations (Lemma~\ref{Universal property of pro-graded rational localization}), the graded map $(T_{n, K})_{\gr}\to \mathcal{A}_{\gr}\to \mathcal{A}\langle f_{i}^{-1}\rangle_{\gr}$ factors through a canonical graded map \begin{equation*}T_{n, K}\langle X_{i}^{-1}\rangle_{\gr}\to \mathcal{A}\langle f_{i}^{-1}\rangle_{\gr}\end{equation*}which is given by the restriction to $T_{n, K}\langle X_{i}^{-1}\rangle_{\gr}$ of the canonical map \begin{equation*}T_{n, K}\langle X_{i}^{-1}\rangle\to \mathcal{A}\langle f_{i}^{-1}\rangle.\end{equation*}Since $\varphi: T_{n, K}\to \mathcal{A}$ is strict surjective, its complete base extension \begin{equation*}T_{n, K}\langle X_{i}^{-1}\rangle\to T_{n, K}\langle X_{i}^{-1}\rangle\widehat{\otimes}_{T_{n, K}}\mathcal{A}=\mathcal{A}\langle f_{i}^{-1}\rangle\end{equation*}is also strict surjective, by \cite{BGR}, Proposition~2.1.8/6. Hence, the graded ring map $T_{n, K}\langle X_{i}^{-1}\rangle_{\gr}\to \mathcal{A}\langle f_{i}^{-1}\rangle_{\gr}$ is surjective, by Lemma~\ref{Surjective morphisms of pro-graded rings}. In particular, it restricts to a surjective ring map \begin{equation*}T_{n, K}\langle X_{i}^{-1}\rangle_{\gr, 0}\to \mathcal{A}\langle f_{i}^{-1}\rangle_{\gr, 0}.\end{equation*}But $T_{n, K}\langle X_{i}^{-1}\rangle_{\gr, 0}\cong T_{n-1, K}$, as in Example~\ref{Adic projective space}. This exhibits $\mathcal{A}\langle f_{i}^{-1}\rangle_{\gr, 0}$ as an affinoid $K$-algebra in the sense of Tate. In particular, this shows that $\mathcal{A}\langle f_{i}^{-1}\rangle_{\gr, 0}$ is sheafy. \end{proof}  
Note that for every pro-graded affinoid $K$-algebra $(\mathcal{A}, \mathcal{A}_{\gr})$, the underlying Banach $K$-algebra $\mathcal{A}$ is topologically generated over $K$ by finitely many homogeneous elements of degree $1$, namely, by the images of the variables $X_1,\dots, X_n$, for any map $(\varphi, \varphi_{\gr})$ as in the definition. In particular, by Corollary~\ref{Topological space 2}, the projective adic spectrum of any pro-graded affinoid algebra is quasi-compact. The following theorem is one of the highlights of this section.
\begin{thm}\label{Dimension theorem}Let $(\mathcal{A}, \mathcal{A}_{\gr})$ be a pro-graded affinoid $K$-algebra over some nonarchimedean field $K$. Then \begin{equation*}\dim(\Proj(\mathcal{A}, \mathcal{A}_{\gr}, \mathcal{A}^{\circ}))=\dim(\Spa(\mathcal{A}, \mathcal{A}^{\circ}))-1=\dim(\mathcal{A})-1.\end{equation*}\end{thm}
We first record the following well-known fact from the rigid-analytic geometry whose proof we include here for the reader's convenience. 
\begin{lemma}\label{Krull dimension of affinoid algebras}Let $\mathcal{A}$ be an affinoid algebra in the sense of Tate over a nonarchimedean field $K$ and let $\widetilde{\mathcal{A}}=\mathcal{A}^{\circ}/\mathcal{A}^{\circ\circ}$. Then $\dim(\mathcal{A})=\dim(\widetilde{\mathcal{A}})$.\end{lemma}
\begin{proof}Let $d$ be the Krull dimension of $\mathcal{A}$. By Noether normalization for affinoid algebras there exists a finite isometric homomorphism $\iota: T_{d, K} \hookrightarrow \mathcal{A}$. Since $\iota$ is isometric, it induces an injective homomorphism $\widetilde{\iota}: \widetilde{K}[X_{1}, \dots, X_{d}]=\widetilde{T_{d}}\hookrightarrow \widetilde{\mathcal{A}}$. By \cite{BGR}, Theorem~6.3.4/2, $\widetilde{\iota}$ is again a finite homomorphism. Hence $\dim(\widetilde{\mathcal{A}})=\dim(\widetilde{K}[X_{1}, \dots, X_{d}])=d$.\end{proof}
\begin{proof}[Proof of Theorem~\ref{Dimension theorem}]The fact that \begin{equation*}\dim(\Spa(\mathcal{A}, \mathcal{A}^{\circ}))=\dim(\mathcal{A})\end{equation*}for any affinoid $K$-algebra $\mathcal{A}$ was proved by Huber in \cite{Huber3}, Lemma~1.8.6.

We first reduce the other assertion of the theorem to the case when $\mathcal{A}$ is reduced and hence uniform. Note that the dimension of $\Proj(\mathcal{A}, \mathcal{A}_{\gr}, \mathcal{A}^{+})$ is the maximum of the dimensions of \begin{equation*}U_{i}=\Spa(\mathcal{A}\langle f_{i}^{-1}\rangle_{\gr,0}, \mathcal{A}\langle f_{i}^{-1}\rangle_{\gr,0}^{+}),\end{equation*}where $f_1,\dots, f_n$ are homogeneous elements of degree~$1$ in $\mathcal{A}_{\gr}^{+}$ generating the irrelevant ideal of $(\mathcal{A}, \mathcal{A}_{\gr})$. Indeed, if $v$ is any closed point of $\Proj(\mathcal{A}, \mathcal{A}_{\gr}, \mathcal{A}^{+})$, then it belongs to one of the open subspaces $U_{i}$ and then all generizations of $v$ are also contained in $U_{i}$. Now consider the nilradical $\Nil(\mathcal{A})$ of $\mathcal{A}$. Then the extension of $\Nil(\mathcal{A})$ to $\mathcal{A}\langle f_{i}^{-1}\rangle$ is still a nilpotent ideal of $\mathcal{A}\langle f_{i}^{-1}\rangle$, so the map $\mathcal{A}\langle f_{i}^{-1}\rangle\to (\mathcal{A}/\Nil(\mathcal{A}))\langle f_{i}^{-1}\rangle$, where, by abuse of notation, we again denote the image of $f_{i}$ in $\mathcal{A}/\Nil(\mathcal{A})$, has nilpotent kernel. In particular, the map \begin{equation*}\mathcal{A}\langle f_{i}^{-1}\rangle_{\gr,0}\to (\mathcal{A}/\Nil(\mathcal{A}))\langle f_{i}^{-1}\rangle_{\gr,0}\end{equation*}has nilpotent kernel. It follows that the dimension of each $U_{i}$ does not change upon replacing $\mathcal{A}$ with its reduced quotient $\mathcal{A}/\Nil(\mathcal{A})$, and hence the dimension of $\Proj(\mathcal{A}, \mathcal{A}_{\gr}, \mathcal{A}^{\circ})$ also does not change. 

Now assume that $\mathcal{A}$ is uniform and thus the open subring $\mathcal{A}^{\circ}$ of $\mathcal{A}$ is a complete adic ring, so that we can apply Corollary~\ref{Dimension estimate for the Proj}. By Lemma~\ref{Krull dimension of affinoid algebras}, $\dim(\mathcal{A})=\dim(\widetilde{\mathcal{A}})$, where $\widetilde{\mathcal{A}}=\mathcal{A}^{\circ}/\mathcal{A}^{\circ\circ}$. By Noether normalization for affinoid algebras and \cite{BGR}, Theorem~6.3.4/2, $\widetilde{\mathcal{A}}$ is a finite-type $\widetilde{K}$-algebra. The ideal $\mathcal{A}^{\circ\circ}$ of $\mathcal{A}^{\circ}$ is homogeneous since it is the radical of the ideal of $\mathcal{A}^{\circ}$ generated by a pseudo-uniformizer $\varpi$ of $K$. Regard $\widetilde{\mathcal{A}}$ as a finitely generated $\mathbb{Z}_{\geq0}$-graded $\widetilde{K}$-algebra, generated by finitely many homogeneous elements of positive degree, with the grading induced from the grading on $\mathcal{A}_{\gr}^{\circ}$. Then \begin{equation*}\dim(\Proj(\widetilde{\mathcal{A}}))=\dim(\widetilde{\mathcal{A}})-1=\dim(\mathcal{A})-1,\end{equation*}for example, by \cite{Hartshorne}, Ch.~I, Ex.~2.6 (for the case of finitely generated graded algebras generated in degree~$1$) together with \cite{EGAII}, Proposition~2.4.7. On the other hand, $\Proj(\widetilde{\mathcal{A}})$ has the same underlying topological space as $\Proj(\mathcal{A}^{\circ}, \mathcal{A}_{\gr}^{\circ})$, so $\dim(\Proj(\mathcal{A}^{\circ}, \mathcal{A}_{\gr}^{\circ}))=\dim(\mathcal{A})-1$. By Corollary~\ref{Dimension estimate for the Proj}, this entails that \begin{equation*}\dim(\Proj(\mathcal{A}, \mathcal{A}_{\gr}, \mathcal{A}^{\circ}))\geq \dim(\mathcal{A})-1.\end{equation*}On the other hand, by \cite{Dine25-2}, Corollary~4.20, the map \begin{equation*}\spc: \Spa(\mathcal{A}, \mathcal{A}^{\circ})\to \Spf(\mathcal{A}^{\circ})\end{equation*}is closed and, in particular, takes closed points to closed points. Thus the center on $\mathcal{A}^{\circ}$ of any closed point of $\Spa(\mathcal{A}, \mathcal{A}^{\circ})$ is a maximal ideal of $\mathcal{A}^{\circ}$ and corresponds to a maximal ideal of the finite-type $\widetilde{K}$-algebra $\widetilde{\mathcal{A}}$. But no maximal ideal of the finitely generated graded $\widetilde{K}$-algebra $\widetilde{\mathcal{A}}$ can define a point of $\Proj(\widetilde{\mathcal{A}})$. In this way, we see that no closed point of $\Spa(\mathcal{A}, \mathcal{A}^{\circ})$ maps to $\Proj(\mathcal{A}^{\circ}, \mathcal{A}_{\gr}^{\circ})$, and thus, that no closed point of $\Spa(\mathcal{A}, \mathcal{A}^{\circ})$ belongs to $\Proj(\mathcal{A}, \mathcal{A}_{\gr}, \mathcal{A}^{\circ})$. It follows that \begin{equation*}\dim(\Proj(\mathcal{A}, \mathcal{A}_{\gr}, \mathcal{A}^{\circ}))\leq \dim(\Spa(\mathcal{A}, \mathcal{A}^{\circ}))-1,\end{equation*}and this concludes the proof.  \end{proof}
We conclude this section by discussing the relationship between pro-graded affinoid algebras and projective varieties over a nonarchimedean field. 
\begin{lemma}\label{Restriction of topologically homogeneous ideals}If $I$ is a closed topologically homogeneous ideal of a pro-$G$-graded Huber ring $(\mathcal{A}, \mathcal{A}_{\gr})$, then for any $G$-graded subring $\mathcal{A}_{\gr}^{\prime}$ of $\mathcal{A}_{\gr}$ the restriction $I\cap \mathcal{A}_{\gr}^{\prime}$ is a homogeneous ideal of $\mathcal{A}_{\gr}^{\prime}$.\end{lemma}
\begin{proof}Follows from Lemma~\ref{Closed topologically homogeneous ideals}. \end{proof} 
\begin{lemma}\label{Kernels of topologically homogeneous valuations}Let $v$ be a topologically homogeneous continuous valuation of a pro-$G$-graded Huber ring $(\mathcal{A}, \mathcal{A}_{\gr})$. Then $\ker(v)$ is a topologically homogeneous prime ideal of $\mathcal{A}$, and for every $G$-graded subring $\mathcal{A}_{\gr}^{\prime}$ of $\mathcal{A}_{\gr}$, the kernel of $v\vert_{\mathcal{A}_{\gr}^{\prime}}$ is a homogeneous prime ideal of $\mathcal{A}_{\gr}^{\prime}$.\end{lemma}
\begin{proof}The second assertion follows from the first by Lemma~\ref{Restriction of topologically homogeneous ideals}. The first assertion follows directly from the definitions. Indeed, let $f\in \ker(v)$ and let $f=\sum_{g\in G}f_{g}$ be its topological homogeneous decomposition in $\mathcal{A}$. Since $v$ is topologically homogeneous, we have $0=v(f)=\max\{\, v(f_{g})\mid g\in G\,\}$, so $v(f_{g})=0$ for all $g\in G$, so $f_{g}\in \ker(v)$ for all $g\in G$.  \end{proof}  
\begin{mydef}[Kernel map]\label{Kernel map, definition}Let $K$ be a nonarchimedean field. We define the kernel map \begin{equation*}\pi: \mathbb{P}_{K}^{n,\ad}=\Proj(T_{n+1, K}, (T_{n+1, K})_{\gr}, T_{n+1, K}^{\circ})\to \mathbb{P}_{K}^{n}=\Proj(K[X_1,\dots, X_{n+1}])\end{equation*}to be the map sending a topologically homogeneous continuous valuation $v$ on $T_{n+1, K}$ to the kernel of its restriction to $K[X_1,\dots, X_{n+1}]$, which is a homogeneous prime ideal by Lemma~\ref{Kernels of topologically homogeneous valuations}.\end{mydef}
\begin{lemma}\label{Analytification}The kernel map from Definition~\ref{Kernel map, definition} coincides with the adic analytification map $\pi_{\ad}: \mathbb{P}_{K}^{n,\ad}\to \mathbb{P}_{K}^{n}$.\end{lemma}
\begin{proof}It suffices to prove that the two maps coincide when restricted to any one of the affinoid opens \begin{equation*}\Proj(T_{n+1, K}, (T_{n+1, K})_{\gr}, T_{n+1, K}^{\circ})(\frac{1}{X_{i}}), i=1,\dots, n+1.\end{equation*}The (adic) analytification map $\pi_{\ad}$ is defined by covering $\mathbb{P}_{K}^{n,\ad}$ with (adic) analytifications $D_{+}(X_{i})^{\ad}$ of the affine spaces \begin{equation*}D_{+}(X_{i})=\{X_{i}\neq 0\}\subseteq \mathbb{P}_{K}^{n}\end{equation*}and on each $D_{+}(X_{i})^{\ad}$ the map is given by sending a valuation on the corresponding polynomial ring $K[Y_1,\dots, Y_n]$ to its kernel.

Fix some $i\in\{1,\dots, n+1\}$ and consider the isomorphisms \begin{equation*}K[X_{1},\dots, X_{n+1}][X_{i}^{-1}]_{0}\cong K[Y_1,\dots, Y_n]\end{equation*}and \begin{equation*}K\langle X_{1},\dots, X_{n+1}\rangle\langle X_{i}^{-1}\rangle_{\gr, 0}\cong K\langle Y_{1},\dots, Y_{n}\rangle\end{equation*}given by $\frac{X_{j}}{X_{i}}\mapsto Y_{j}$ for $j<i$ and $\frac{X_{j}}{X_{i}}\mapsto Y_{j-1}$ for $j>i$. Regard $K[Y_1,\dots, Y_n]$ and $K\langle Y_1,\dots, Y_n\rangle$ as subrings of $K[X_1,\dots, X_{n+1}][X_{i}^{-1}]$ and $K\langle X_1,\dots, X_{n+1}\rangle\langle X_{i}^{-1}\rangle$ via these isomorphisms. For a continuous, topologically homogeneous valuation \begin{equation*}v\in \Proj(T_{n+1, K}, (T_{n+1, K})_{\gr}, T_{n+1, K}^{\circ})(\frac{1}{X_{i}}),\end{equation*}let $v'$ denote the canonical extension of $v$ to a continuous, topologically homogeneous valuation on $T_{n+1, K}\langle X_{i}^{-1}\rangle$, given by \begin{equation*}v'(\frac{F}{X_{i}^{k}})=v(F)\end{equation*}for all $F\in T_{n+1, K}$, $k\geq1$. We have to check that the diagram \begin{center}\begin{tikzcd}\Proj(T_{n+1, K}, (T_{n+1, K})_{\gr}, T_{n+1, K}^{\circ})(\frac{1}{X_{i}}) \arrow{r}{\raisebox{10pt}{$v\mapsto \ker(v\vert_{K[X_1,\dots, X_{n+1}]})$}} \arrow{d}{v\mapsto v'\vert_{K\langle Y_1,\dots, Y_n\rangle}} & D_{+}(X_{i})=\{X_{i}\neq 0\}\subseteq \mathbb{P}_{K}^{n} \arrow{d}{\mathfrak{p}\mapsto \mathfrak{p}K[X_1,\dots, X_{n+1}][X_{i}^{-1}]\cap K[Y_1,\dots, Y_n]} \\ \mathbb{B}_{K}^{n}=\Spa(K\langle Y_1,\dots, Y_n\rangle, K^{\circ}\langle Y_1,\dots, Y_n\rangle) \arrow[hook]{d}{w\mapsto w\vert_{K[Y_1,\dots, Y_n]}} & \Spec(K[Y_1,\dots, Y_n]) \\ \mathbb{A}_{K}^{n,\ad} \arrow[swap]{ur}{w\mapsto \ker(w)} \end{tikzcd}\end{center}is commutative. In other words, we have to show that for every topologically homogeneous $v\in \Spa(T_{n+1, K}, T_{n+1, K}^{\circ})$ with $v(X_{i})=1$, we have \begin{equation*}\ker(v\vert_{K[X_1,\dots, X_{n+1}]})K[X_1,\dots, X_{n+1}][X_{i}^{-1}]=\ker(v')\cap K[Y_1,\dots, Y_n].\end{equation*}To check this, let $f\in \ker(v\vert_{K[X_1,\dots, X_{n+1}]})K[X_1,\dots, X_{n+1}][X_{i}^{-1}]\cap K[Y_1,\dots, Y_n]$ for a valuation $v$ as above and write $f=\frac{F}{X_{i}^{k}}$ with $F\in \ker(v)$, $k\geq0$. Then \begin{equation*}f=\frac{F}{X_{i}^{k}}\in \ker(v'\vert_{K[Y_1,\dots, Y_n]}),\end{equation*}proving one of the two desired inclusions. Conversely, if $f=\frac{F}{X_{i}^{k}}\in \ker(v')\cap K[Y_1,\dots, Y_n]$, then, by the definition of $v'$, we have $v(F)=0$ and thus $f\in \ker(v)$.  \end{proof}
For a pro-$G$-graded Huber pair $((\mathcal{A}, \mathcal{A}_{\gr}), \mathcal{A}^{+})$ satisfying the assumptions of Definition~\ref{Adic Proj, part 2} and a closed topologically homogeneous ideal $I$ of $\mathcal{A}$, we set \begin{equation*}\mathcal{V}_{+}(I)=\{\, v\in \Proj(\mathcal{A}, \mathcal{A}_{\gr}, \mathcal{A}^{+})\mid \ker(v)\supseteq I\,\}.\end{equation*}We also use the analogous notation for the closed subset of $\Proj(R)$ defined by a homogeneous ideal $I$ of $R$, when $R$ is a usual graded ring. We can now describe Zariski-closed subsets of adic projective space in terms of homogeneous ideals.
\begin{prop}\label{Zariski-closed subsets and homogeneous ideals}Let $K$ be a nonarchimedean field, let $\mathcal{I}$ be a coherent ideal sheaf on $\mathbb{P}_{K}^{n,\ad}$ defining the Zariski-closed subset \begin{equation*}\mathcal{V}(\mathcal{I})=\{\, x\in \mathbb{P}_{K}^{n,\ad}\mid \mathcal{I}_{x}\neq \mathcal{O}_{\mathbb{P}_{K}^{n,\ad},x}\,\}\subseteq \mathbb{P}_{K}^{n,\ad}\end{equation*}and let $\mathcal{J}=\mathcal{I}_{\alg}$ be the coherent ideal sheaf on $\mathbb{P}_{K}^{n}$ corresponding to $\mathcal{I}$ by GAGA. Set \begin{equation*}\Gamma_{\ast}(\mathcal{J})=\bigoplus_{m\in\mathbb{Z}}\Gamma(\mathbb{P}_{K}^{n}, \mathcal{J}(m)),\end{equation*}a homogeneous ideal of $K[X_1,\dots, X_{n+1}]$. Then \begin{equation*}I=\Gamma_{\ast}(\mathcal{J})T_{n+1, K}\end{equation*}is a homogeneous ideal of $(T_{n+1, K}, (T_{n+1, K})_{\gr})$ satisfying \begin{equation*}\mathcal{V}_{+}(I)=\mathcal{V}(\mathcal{I})\end{equation*}as subsets of $\mathbb{P}_{K}^{n,\ad}$.\end{prop}
\begin{proof}Let $\pi=\pi_{\ad}: \mathbb{P}_{K}^{n,\ad}\to \mathbb{P}_{K}^{n}$ be the adic analytification map. By Lemma~\ref{Analytification}, we see that \begin{equation*}\mathcal{V}_{+}(I)=\pi^{-1}(\mathcal{V}(\Gamma_{\ast}(\mathcal{J})))=\pi^{-1}(\{\, \mathfrak{p}\in \mathbb{P}_{K}^{n}=\Proj(K[X_1,\dots, X_{n+1}])\mid \mathfrak{p}\supseteq \Gamma_{\ast}(\mathcal{J})\,\}).\end{equation*}By \cite{Hartshorne}, Ch.~II, proof of Corollary 5.16(a), \begin{equation*}\mathcal{V}_{+}(\Gamma_{\ast}(\mathcal{J}))=\mathcal{V}(\mathcal{J})=\{\, x\in \mathbb{P}_{K}^{n}\mid \mathcal{J}_{x}=0\,\}.\end{equation*}But on the other hand, $\mathcal{I}=\pi^{\ast}\mathcal{J}$ and the map of locally ringed spaces $\pi$ is flat, so, for all $x\in \mathbb{P}_{K}^{n,\ad}$, \begin{equation*}\mathcal{I}_{x}=\mathcal{J}_{\pi(x)}\mathcal{O}_{\mathbb{P}_{K}^{n,\ad}, x}.\end{equation*}It follows that \begin{equation*}\mathcal{V}_{+}(I)=\pi^{-1}(\mathcal{V}(\mathcal{J}))=\mathcal{V}(\mathcal{I}),\end{equation*}as claimed.  \end{proof}
\begin{lemma}\label{Surjection of pro-graded Huber pairs}Let $G$ be a totally ordered abelian group and let \begin{equation*}(\varphi, \varphi_{\gr}): ((\mathcal{A}, \mathcal{A}_{\gr}), \mathcal{A}^{+})\to ((\mathcal{B}, \mathcal{B}_{\gr}), \mathcal{B}^{+})\end{equation*}be a morphism of graded-sheafy pro-$G_{\geq0}$-graded Huber pairs such that $\varphi$ is surjective and $\mathcal{B}^{+}$ is the integral closure of $\varphi(\mathcal{A}^{+})$ inside $\mathcal{B}$. Then the continuous map \begin{equation*}\Spa(\varphi): \Spa(\mathcal{B}, \mathcal{B}^{+})\to \Spa(\mathcal{A}, \mathcal{A}^{+})\end{equation*}restricts to a homeomorphism \begin{equation*}\Proj(\mathcal{B}, \mathcal{B}_{\gr}, \mathcal{B}^{+})\tilde{\to} \mathcal{V}_{+}(I),\end{equation*}where $I=\ker(\varphi)$.\end{lemma}
\begin{proof}The map $\Spa(\varphi)$ is defined by $\Spa(\varphi)(v)(f)=v(\varphi(f))=v(f+I)$ for $f\in \mathcal{A}$. By Lemma~\ref{Topologically homogeneous valuations and morphisms}, $\Spa(\varphi)(v)$ is topologically homogeneous for every topologically homogeneous valuation $v$ in $\Spa(\mathcal{B}, \mathcal{B}^{+})$. Moreover, \begin{equation*}\mathcal{A}_{\gr}/(I\cap \mathcal{A}_{\gr})=\bigoplus_{g\in G}\mathcal{A}_{\gr, g}/(\mathcal{A}_{\gr, g}\cap I)\end{equation*}and, in particular, the homogeneous elements of positive degree in $\mathcal{A}_{\gr}/(I\cap \mathcal{A}_{\gr})$ are precisely the images of homogeneous elements of positive degree in $\mathcal{A}_{\gr}$. It follows that $v(f)<1$ for all homogeneous elements $f$ of positive degree in $\mathcal{A}_{\gr}^{+}/(I\cap \mathcal{A}_{\gr}^{+})$ if and only if $\Spa(\varphi)(v)$ is $<1$ on all homogeneous elements of positive degree in $\mathcal{A}_{\gr}^{+}$. But if $f\in \mathcal{B}_{\gr}^{+}$ is any homogeneous element of positive degree, then $f$ is integral over $\mathcal{A}^{+}/(I\cap \mathcal{A}^{+})$ by the hypothesis on $\varphi$. By looking at topological homogeneous decompositions of coefficients of such an equation of integral dependence, we see that $f$ must satisfy such an equation with coefficients being homogeneous elements of positive degree. Hence $v$ being $<1$ on all homogeneous elements of positive degree in $\mathcal{A}_{\gr}^{+}/(I\cap \mathcal{A}_{\gr}^{+})$ also implies $v(f)<1$. This shows that the center of $v$ on $\mathcal{B}^{+}$ contains $(\mathcal{B}, \mathcal{B}_{\gr})_{+}\cap \mathcal{B}^{+}$ if and only if the center of $\Spa(\varphi)(v)$ on $\mathcal{A}^{+}$ contains $(\mathcal{A}, \mathcal{A}_{\gr})_{+}\cap \mathcal{A}^{+}$. In particular, the map $\Spa(\varphi)$ restricts to a map \begin{equation*}\Proj(\mathcal{B}, \mathcal{B}_{\gr}, \mathcal{B}^{+})\to \Proj(\mathcal{A}, \mathcal{A}_{\gr}, \mathcal{A}^{+}).\end{equation*}Since $\Spa(\varphi)$ is a homeomorphism onto its image, this map is also a homeomorphism onto its image. It remains to show that the image of this map is precisely $\mathcal{V}_{+}(I)$, and for this it remains to show that $\Spa(\varphi)(v)$ being topologically homogeneous also implies that $v$ is topologically homogeneous. Suppose that $v$ is such that $\Spa(\varphi)(v)$ is topologically homogeneous. Let $\overline{f}\in \mathcal{B}$ with topological homogeneous decomposition $\overline{f}=\sum_{g\in G}\overline{f}_{g}$. Choose $f\in \mathcal{A}$ with $\varphi(f)=\overline{f}$, with topological homogeneous decomposition $f=\sum_{g\in G}f_{g}$ in $\mathcal{A}$. By Lemma~\ref{Homogeneous decompositions and morphisms}, $\varphi(f_{g})=\overline{f}_{g}$ for all $g\in G$. Then \begin{equation*}v(\overline{f})=\Spa(\varphi)(v)(f)=\max\{\, \Spa(\varphi)(v)(f_{g})\mid g\in G\,\}=\max\{\, v(\overline{f}_{g})\mid g\in G\,\},\end{equation*}as desired. \end{proof}
\begin{mydef}[Adic space structure on $\mathcal{V}_{+}(I)$]\label{Adic space structure on Zariski-closed subsets}Let $K$ be a nonarchimedean field, let $(\mathcal{A}, \mathcal{A}_{\gr})$ be a pro-graded affinoid $K$-algebra and let $I$ be a closed topologically homogeneous ideal of $\mathcal{A}$. We endow the closed subset \begin{equation*}\mathcal{V}_{+}(I)\subseteq \Proj(\mathcal{A}, \mathcal{A}_{\gr}, \mathcal{A}^{+})\end{equation*}with the structure of an adic space by identifying it with $\Proj(\mathcal{A}/I, \mathcal{A}_{\gr}/(I\cap \mathcal{A}_{\gr}), (\mathcal{A}/I)^{+})$ via the homeomorphism from Lemma~\ref{Surjection of pro-graded Huber pairs}\footnote{Note that the quotient $(\mathcal{A}/I, \mathcal{A}_{\gr}/(I\cap \mathcal{A}_{\gr}))$ is again a pro-graded affinoid $K$-algebra. Moreover, the ring of integral elements $(\mathcal{A}/I)^{+}$ is equal to $(\mathcal{A}/I)^{+}$ by \cite{BGR}, Proposition~6.3.4/1. Hence \begin{equation*}((\mathcal{A}/I, \mathcal{A}_{\gr}/(I\cap \mathcal{A}_{\gr})), (\mathcal{A}/I)^{+})\end{equation*}is a pro-graded Huber pair by Lemma~\ref{Pro-graded affinoid algebras define pro-graded Huber pairs} and Lemma~\ref{Pro-graded affinoids are graded-sheafy}.}. We call this the adic space structure on $\mathcal{V}_{+}(I)$ induced by $I$. \end{mydef}
\begin{prop}\label{Zariski-closed subsets and homogeneous ideals 2}Set assumptions and notation as in Proposition~\ref{Zariski-closed subsets and homogeneous ideals}. Endow $\mathcal{V}_{+}(I)$ with the adic space structure induced by $I$. Then $\mathcal{V}_{+}(I)=\mathcal{V}(\mathcal{I})$ as adic spaces.\end{prop}
\begin{proof}Let $f_1,\dots, f_{n+1}$ be the images of $X_1,\dots, X_{n+1}$ under the map $T_{n+1, K}\to T_{n+1, K}/I$. By definition, \begin{equation*}\mathcal{V}_{+}(I)=\Proj(T_{n+1, K}/I, (T_{n+1, K})_{\gr}/(I\cap (T_{n+1, K})_{\gr}), (T_{n+1, K}/I)^{+})\end{equation*}is obtained, as an adic space, by gluing the affinoid adic spaces \begin{equation*}\Spa((T_{n+1, K}/I)\langle f_{i}^{-1}\rangle_{\gr, 0}, (T_{n+1, K}/I)\langle f_{i}^{-1}\rangle_{\gr, 0}^{+})\end{equation*}along their intersections \begin{equation*}\Spa((T_{n+1, K}/I)\langle (f_{i}f_{j})^{-1}\rangle_{\gr, 0}, (T_{n+1, K}/I)\langle (f_{i}f_{j})^{-1}\rangle_{\gr, 0}^{+}).\end{equation*}Hence it suffices to prove that, for every $i=1,\dots, n+1$, the intersection of $\mathcal{V}(\mathcal{I})$ with \begin{equation*}\Proj(T_{n+1, K}, (T_{n+1, K})_{\gr}, T_{n+1, K}^{\circ})(\frac{1}{X_{i}})=\{\, v\in \mathbb{P}_{K}^{n,\ad}\mid v(X_{i})=1\,\}\end{equation*}is equal to $\Spa((T_{n+1, K}/I)\langle f_{i}^{-1}\rangle_{\gr, 0}, (T_{n+1, K}/I)\langle f_{i}^{-1}\rangle_{\gr, 0}^{+})$. Fix some $i\in\{1, \dots, n+1\,\}$. Note that \begin{equation*}(T_{n+1, K}/I)\langle f_{i}^{-1}\rangle=T_{n+1, K}\langle X_{i}^{-1}\rangle/IT_{n+1, K}\langle X_{i}^{-1}\rangle\end{equation*}and thus \begin{equation*}(T_{n+1, K}/I)\langle f_{i}^{-1}\rangle_{\gr, 0}=K\langle Y_1,\dots, Y_n\rangle/(IT_{n+1, K}\langle X_{i}^{-1}\rangle\cap K\langle Y_1,\dots, Y_n\rangle),\end{equation*}where $Y_{j}=\frac{X_{j}}{X_{i}}$ for $j<i$ and $Y_{j}=\frac{X_{j+1}}{X_{i}}$ for $j\geq i$. 

On the other hand, let $\mathcal{J}$ be the coherent ideal sheaf on the algebraic projective space $\mathbb{P}_{K}^{n}$ corresponding via GAGA to $\mathcal{I}$, so that $I=\Gamma_{\ast}(\mathcal{J})T_{n+1, K}$. Using \cite{Hartshorne}, Ch.~II, proof of Corollary 5.16(a), the intersection of $\mathcal{V}(\mathcal{J})$ with the affine open $\{X_{i}\neq 0\}$ of $\mathbb{P}_{K}^{n}$ is given by \begin{align*}\mathcal{V}(\mathcal{J})\cap \{X_{i}\neq 0\}=\mathcal{V}_{+}(\Gamma_{\ast}(\mathcal{J}))\cap \{X_{i}\neq 0\} \\ =\Spec(K[Y_1,\dots, Y_n]/(\Gamma_{\ast}(\mathcal{J})K[X_1,\dots, X_{n+1}][X_{i}^{-1}]\cap K[Y_1,\dots, Y_n])).\end{align*}But the pre-image in $\mathbb{P}_{K}^{n,\ad}$ of the latter affine open subset of $\mathbb{P}_{K}^{n}$ intersects \begin{equation*}\{\, v\in\mathbb{P}_{K}^{n,\ad}\mid v(X_{i})=1\,\}\end{equation*}precisely in the Zariski-closed adic subspace of $\Spa(K\langle Y_1,\dots, Y_n\rangle)$ defined by \begin{align*}K\langle Y_1,\dots, Y_n\rangle/(\Gamma_{\ast}(\mathcal{J})T_{n+1, K}\langle X_{i}^{-1}\rangle\cap K\langle Y_1,\dots, Y_n\rangle) \\ =K\langle Y_1,\dots, Y_n\rangle/(IT_{n+1, K}\langle X_{i}^{-1}\rangle\cap K\langle Y_1,\dots, Y_n\rangle).\end{align*}This concludes the proof of the proposition.\end{proof}
We can now finally establish a correspondence between pro-graded affinoid algebras and projective varieties over a nonarchimedean field $K$.
\begin{thm}\label{Graded GAGA}Let $K$ be a nonarchimedean field. For every closed subscheme \begin{equation*}Z\subseteq \mathbb{P}_{K}^{n}\end{equation*}with adic analytification $Z^{\ad}$ there exists a pro-graded affinoid $K$-algebra $(\mathcal{A}, \mathcal{A}_{\gr})$ such that \begin{equation*}Z^{\ad}=\Proj(\mathcal{A}, \mathcal{A}_{\gr}, \mathcal{A}^{\circ})\end{equation*}as adic spaces. Conversely, if $(\mathcal{A}, \mathcal{A}_{\gr})$ is a pro-graded affinoid $K$-algebra, then there exists a closed subscheme $Z\subseteq \mathbb{P}_{K}^{n}$ with $Z^{\ad}=\Proj(\mathcal{A}, \mathcal{A}_{\gr}, \mathcal{A}^{\circ})$. \end{thm}
\begin{proof}Let $Z\subseteq \mathbb{P}_{K}^{n}$ be a closed subscheme and let $\mathcal{I}$ be a coherent ideal sheaf on $\mathbb{P}_{K}^{n,\ad}$ such that $Z^{\ad}=\mathcal{V}(\mathcal{I})$ as adic spaces. Proposition~\ref{Zariski-closed subsets and homogeneous ideals 2} yields a homogeneous ideal $I$ of $T_{n+1, K}$ such that \begin{equation*}Z^{\ad}=\mathcal{V}_{+}(I)=\Proj(T_{n+1, K}/I, (T_{n+1, K})_{\gr}/(I\cap (T_{n+1, K})_{\gr}), (T_{n+1, K}/I)^{+})\end{equation*}as adic spaces. By \cite{BGR}, Proposition~6.3.4/1, $(T_{n+1, K}/I)^{+}=(T_{n+1, K}/I)^{\circ}$, so this proves one assertion of the theorem.

Conversely, let $(\mathcal{A}, \mathcal{A}_{\gr})$ be a pro-graded affinoid $K$-algebra, let \begin{equation*}(\varphi, \varphi_{\gr}): (T_{n+1, K}, (T_{n+1, K})_{\gr})\to (\mathcal{A}, \mathcal{A}_{\gr})\end{equation*}be a morphism of pro-graded Huber rings such that $\varphi$ is surjective and let $I=\ker(\varphi_{\gr})$, a homogeneous ideal of the graded ring $(T_{n+1, K})_{\gr}$. For every $i=1,\dots, n+1$, we define a coherent ideal sheaf $\mathcal{I}_{i}$ on the affinoid open subspace \begin{equation*}U_{i}=\Proj(T_{n+1, K}, (T_{n+1, K})_{\gr}, T_{n+1, K}^{\circ})(\frac{1}{X_{i}})=\Spa(T_{n+1, K}\langle X_{i}^{-1}\rangle_{\gr, 0}, T_{n+1, K}\langle X_{i}^{-1}\rangle_{\gr, 0}^{+})\end{equation*}by setting \begin{equation*}\mathcal{I}_{i}(V)=IT_{n+1, K}\langle X_{i}^{-1}\rangle_{\gr, 0}\otimes_{\mathcal{O}_{U_{i}}(U_{i})}\mathcal{O}_{U_{i}}(V)\end{equation*}for every rational subset $V$. On the intersections $U_{i}\cap U_{j}$, we have \begin{equation*}\mathcal{I}_{i}\vert_{U_{i}\cap U_{j}}=\mathcal{I}_{j}\vert_{U_{i}\cap U_{j}}\end{equation*}since \begin{equation*}\mathcal{I}_{i}(U_{i}\cap U_{j})=IT_{n+1, K}\langle (X_{i}X_{j})^{-1}\rangle_{\gr, 0}=\mathcal{I}_{j}(U_{i}\cap U_{j}).\end{equation*}Thus the sheaves $\mathcal{I}_{i}$ glue to a coherent ideal sheaf $\mathcal{I}$ on $\mathbb{P}_{K}^{n,\ad}$ with $\mathcal{V}_{+}(I)=\mathcal{V}(\mathcal{I})$. By GAGA, this implies that $\mathcal{V}_{+}(I)$ is equal to the adic analytification of a closed subscheme $Z$ of $\mathbb{P}_{K}^{n}$.  \end{proof}

\bibliographystyle{plain} 
\bibliography{Bib}

\begin{thebibliography}{10}

\bibitem{Berkovich95Drinfeld}
V.~Berkovich.
\newblock The automorphism group of the {D}rinfeld half-plane.
\newblock {\em C. R. Acad. Sci. Paris}, 321:1127--1132, 1995.

\bibitem{Berthelot96}
P.~Berthelot.
\newblock Cohomologie rigide et cohomologie rigide {\`a} support propres.
\newblock Preprint 96-03, Institut de Recherche Math. de Rennes, 1996.

\bibitem{BhattNotes}
B.~Bhatt.
\newblock Lecture notes for a class on perfectoid spaces.
\newblock Available at
  \url{https://www.math.ias.edu/~bhatt/teaching/mat679w17/lectures.pdf}, 2017.

\bibitem{BGR}
S.~Bosch, U.~G{\"u}ntzer, and R.~Remmert.
\newblock {\em Non-Archimedean analysis. A systematic approach to
  rigid-analytic geometry}, volume 261 of {\em Grundlehren der Math. Wiss.}
\newblock Springer-Verlag, Berlin, 1984.

\bibitem{Dine22}
D.~Dine.
\newblock Topological spectrum and perfectoid {T}ate rings.
\newblock {\em Algebra Number Theory}, 16(6):1463--1500, 2022.

\bibitem{Dine25-3}
D.~Dine.
\newblock Banach modules, almost mathematics and condensed mathematics.
\newblock {\em arXiv e-prints}, arXiv:2508.11268, 2025.

\bibitem{Dine25-2}
D.~Dine.
\newblock Formal models for relative adic spaces.
\newblock {\em arXiv e-prints}, arXiv:2507.11073, 2025.

\bibitem{Dine-Ishizuka23}
D.~Dine and R.~Ishizuka.
\newblock Tilting and untilting for ideals in perfectoid rings.
\newblock {\em Math. Z.}, 307(66), 2024.

\bibitem{EGAII}
A.~Grothendieck and J.~Dieudonn{\'e}.
\newblock {\'E}l{\'e}ments de g{\'e}om{\'e}trie alg{\'e}brique: Ii. {\'e}tude
  globale {\'e}l{\'e}mentaire de quelques classes de morphismes.
\newblock {\em Publ. Math. IH{\'E}S}, 8:5--222, 1961.

\bibitem{Hartshorne}
R.~Hartshorne.
\newblock {\em Algebraic {G}eometry}, volume~52 of {\em Graduate Texts in
  Mathematics}.
\newblock Springer-Verlag, New York-Heidelberg, 1977.

\bibitem{Henkel}
T.~Henkel.
\newblock An open mapping theorem for rings with a zero sequence of units.
\newblock arXiv:1407.5647v2, 2014.

\bibitem{Huber1}
R.~Huber.
\newblock Continuous valuations.
\newblock {\em Math. Z.}, 212:455--477, 1993.

\bibitem{Huber3}
R.~Huber.
\newblock {\em {\'E}tale cohomology of rigid analytic varieties and adic
  spaces}, volume E30 of {\em Aspects of Mathematics}.
\newblock Friedr. Vieweg \& Sohn, Braunschweig, 1996.

\bibitem{Ishizuka-Yoshikawa25}
R.~Ishizuka and S.~Yoshikawa.
\newblock Graded perfectoid rings.
\newblock {\em arXiv e-prints}, arXiv:2511.02322, 2025.

\bibitem{Ishizuka-Yoshikawa26}
R.~Ishizuka and S.~Yoshikawa.
\newblock Algebraization of absolute perfectoidization via section rings.
\newblock {\em arXiv e-prints}, arXiv:2604.02682, 2026.

\bibitem{Kedlaya17}
K.S. Kedlaya.
\newblock Sheaves, stacks and shtukas.
\newblock In {\em Perfectoid {S}paces: {L}ectures from the 2017 Arizona Winter
  School}, volume 242 of {\em Mathematical Surveys and Monographs}, pages
  58--205. Amer. Math. Soc., 2019.

\bibitem{RaynaudReport}
M.~Raynaud.
\newblock G{\'e}om{\'e}trie analytique rigide d'apr{\`e}s {T}ate, {K}iehl,...
\newblock {\em M{\'e}m. Soc. Math. Fr.}, 39-40:319--327, 1974.

\bibitem{Scholze}
P.~Scholze.
\newblock Perfectoid spaces.
\newblock {\em Publ. Math. IH{\'E}S}, 116:245--313, 2012.

\bibitem{Swanson-Huneke}
I.~Swanson and C.~Huneke.
\newblock {\em Integral {C}losure of {I}deals, {R}ings, and {M}odules}, volume
  336 of {\em London Mathematical Society Lecture Note Series}.
\newblock Cambridge University Press, 2006.

\end{thebibliography}
\end{document}